\documentclass[twoside,11pt]{article}

\usepackage{blindtext}
\usepackage{lastpage}
\usepackage{jmlr2e}

\usepackage[english]{babel}	
\usepackage[utf8x]{inputenc}		
\usepackage{enumerate}
\usepackage{textcomp}                
\usepackage{typearea}
\usepackage{pdfpages}
\usepackage[toc]{appendix}
\usepackage[disable]{todonotes}
\usepackage{subcaption}

\let\proof\relax

\usepackage{mathtools, amsthm, amsfonts, dsfont}
\usepackage{booktabs}

\usepackage{pgf, tikz}
\usepackage{bm}

\theoremstyle{definition}

\newtheorem*{assumption*}{Assumption}

\newtheorem*{condition*}{Condition}
\newtheorem{Example}{Example}
\theoremstyle{plain}
\newtheorem{Theorem}{Theorem}
\newtheorem{Proposition}{Proposition}
\newtheorem{Corollary}{Corollary}
\newtheorem{Lemma}{Lemma}
\theoremstyle{remark}
\newtheorem{Remark}{Remark}

\newenvironment{myproof}[1]{\paragraph{Proof  {#1}:}}{\hfill$\square$}

\usepackage{graphicx}
\usepackage{hyperref}
\allowdisplaybreaks

\newcommand{\N}{\mathbb{N}}
\newcommand{\Z}{\mathbb{Z}}

\newcommand{\R}{\mathbb{R}}

\renewcommand{\P}{\mathbb{P}}
\newcommand{\E}{\mathbb{E}}

\newcommand{\Var}{\operatorname{Var}}
\newcommand{\Cov}{\operatorname{Cov}}

\newcommand{\wt}[1]{\widetilde{#1}}

\newcommand{\calB}{\mathcal{B}}

\newcommand{\calD}{\mathcal{D}}

\newcommand{\calF}{\mathcal{F}}

\newcommand{\calL}{\mathcal{L}}

\newcommand{\calN}{\mathcal{N}}

\newcommand{\calQ}{\mathcal{Q}}
\newcommand{\calS}{\mathcal{S}}
\newcommand{\calT}{\mathcal{T}}
\newcommand{\calU}{\mathcal{U}}

\newcommand{\eins}{{\bm 1}}
\newcommand{\matA}{{\mathbf{A}}}
\newcommand{\matB}{{\mathbf{B}}}

\newcommand{\matM}{{\mathbf{M}}}

\newcommand{\matS}{{\mathbf{S}}}

\newcommand{\matV}{{\mathbf{V}}}
\newcommand{\matW}{{\mathbf{W}}}

\newcommand{\vecnull}{{\bm 0}}

\newcommand{\veca}{{\bm a}}

\newcommand{\vecf}{{\bm f}}

\newcommand{\vecG}{{\bm G}}
\newcommand{\vecH}{{\bm H}}

\newcommand{\vecv}{{\bm v}}

\newcommand{\vecw}{{\bm w}}
\newcommand{\vecx}{{\bm x}}
\newcommand{\vecX}{{\bm X}}
\newcommand{\vecy}{{\bm y}}
\newcommand{\vecY}{{\bm Y}}
\newcommand{\vecZ}{{\bm Z}}

\newcommand{\bfmu}{\bm \mu}
\newcommand{\bfbeta}{\bm\beta}

\newcommand{\bftheta}{{\bm\theta}}

\newcommand{\bfeta}{\bm\eta}

\newcommand{\bfeps}{\bm \epsilon}

\newcommand{\bfGamma}{\bm\Gamma}
\newcommand{\bfSigma}{\bm\Sigma}

\newcommand{\blau}{black}

\jmlrheading{27}{2026}{1-\pageref{LastPage}}{03/23; Revised 9/24}{1/26}{23-0274
}{Ansgar Steland}

\ShortHeadings{When to Retrain Deep Learners?}{Steland}
\firstpageno{1}

\begin{document}
	
\title{Online Detection of Changes in Moment--Based Projections: When to Retrain Deep Learners or Update Portfolios?}
	
	\author{\name Ansgar Steland \email steland@stochastik.rwth-aachen.de \\
		\addr Institute of Statistics and AI Center \\
		RWTH Aachen University \\
		52062 Aachen, Germany}
	
	\editor{Aur\'elion Garivier}
	
	\maketitle

\begin{abstract}
	Training deep learning neural networks often requires massive amounts of computational ressources. We propose
	to sequentially monitor network predictions to trigger retraining only if the predictions are no longer valid. This can reduce drastically computational costs  and opens a door to green deep learning. Our approach is based on the relationship to projected second moments monitoring, a problem also arising in other areas such as computational finance. Various open--end as well as closed--end monitoring rules are studied under mild assumptions on the training sample and the observations of the monitoring period. The results allow for high--dimensional non-stationary time series data and thus, especially, non--i.i.d. training data. Asymptotics is based on Gaussian approximations of projected partial sums allowing for an estimated projection vector. Estimation of projection vectors is studied both for classical non--$\ell_0$--sparsity as well as under sparsity. For the case that the optimal projection depends on the unknown covariance matrix, hard-- and soft--thresholded estimators are studied. The method is analyzed by simulations and supported by synthetic data experiments.
\end{abstract}


\section{Introduction} 

\color{\blau} Despite the success of fast algorithms such as ADAM, \color{black} training of deep learners requires massive amounts of computational resources, and when the data distribution changes, called  {\em concept drift}, the network predictions may be outdated and render the network useless. This asks for a monitoring approach to trigger retraining, if (and only if) the predictions are no longer valid.  \color{\blau} A crucial observation is that neural networks with linear output layer compute inner products between the features derived by the hidden layers and the connection weights of  that output layer yielding the prediction for some input $ \vecx $. If the distribution of the input samples changes and thus differs from the distribution in the learning sample used to train the net, these predictions will be outdated and may lead to biased net outputs. This suggests the following approach: When using the deep learner to compute predictions for input data, a monitoring algorithm is used to detect when the predictions are no longer compatible with the training sample thus indicating the need to retrain the net. Our solution monitors efficiently in an online fashion a second moment functionals associated to the predictions, which evaluates the features derived by the deep learner, and signals the need to retrain the network.

Such second moment functionals, even in a much simpler form, also arise in computational finance when determining optimal portfolios. Thus, before proceeding with our primary deep learning application, let us briefly discuss this example from finance: \color{black} The risk of a portfolio $ \vecv $ of $d$ risky assets with log returns $ \vecY $ is usually measured by the variance $ \vecv^\top \bfSigma \vecv $ of the portfolio return $ \vecv^\top \vecY $, where $ \bfSigma = \Var( \vecY )$ is the $d \times d $ covariance matrix. It is well known that the Markowitz--optimal portfolio explicitly depends on the structure of the first two moments and requires to compute the inverse $ \bfSigma^{-1} $. It should be updated if and only if its risk has changed, since for a large investment universe such an update leads to substantial computational costs and, even more important, substantial trading costs. 

\color{\blau} The problem to examine neural network predictions breaks down to a similar quadratic form. Anticipating the main message behing the detailed discussion in Section~\ref{Sec:Applications}, the basic idea is as follows:  \color{black} Recall that deep neural networks with $H$ hidden layers and linear output process inputs $ \vecX $ and linearly combine derived features $ \vecZ = f_H( \vecX ) $ by some network function $ f_H $ to determine the net output $  \bfbeta^\top \vecZ $ where $ \bfbeta $ denotes the vector of output weights. This applies regardless of the network topology. \color{\blau} Here, for simplicity of presentation, we consider a single output neuron.  \color{\blau} Since the network prediction for a new input $ \vecX^* $ 
is given by the projection $ \hat Y^* = \hat\bfbeta_m^\top f_H( \vecX^*; \hat\bftheta_m ) $, where $ (\hat\bftheta_m, \hat \bfbeta_m) $ are the trained network parameters regarded as estimators of the true but unknown optimal weights $ (\bftheta_0,\bfbeta_0) $,  it is clear that the net needs retraining using an updated training sample, only if the change is relevant for that projection. This can be assessed by monitoring the second moment $ \E( \bfbeta^\top f_H( \vecX^* ))^2 = \bfbeta^\top \matM^H \bfbeta $, which is a quadratic form with coefficient matrix $ \matM_n^H = \E ( f_H(\vecX^*;\bftheta_0) f_H(\vecX^*;\bftheta_0)^\top ) $ given by second moment matrix of the last hidden layer, which represents the features learned by the net. 
The practical algorithms discussed later will be based on associated quantities where unknowns are estimated from the training sample and compared to current data.
\color{black}

We approach this machine learning problem from a sequential statistical perspective, which allows to draw on what is known on sequential detection. \color{\blau} It is assumed that after having trained the deep learner, the trained net is used to generate predictions of a sequential input data stream. We use these predictions to monitor the stream, in order to derive a signal when the current data are no longer compatible with the training sample indicating that the data distribution has changed and the neural network needs to be retrained.  \color{black} Our solution employs the open-end monitoring framework usually addressed to  \cite{ChuStinchcombeWhite1996}, which incorporates such a training sample of size $m$ satisfying the no--change null hypothesis of interest and allows for open--end monitoring. Open--end means that the stream may be monitored inifinitely long. Since often neural networks are retrained on a coarse regular grid, we also consider the closed--end setting where monitoring stops latest at a fixed time horizon. The training sample is used to estimate unknowns and to compare current data of the monitoring period with stable training data. We study several detectors based on a  cumulated sum (CUSUM) statistic following \cite{HHKS2004} to monitor regression residuals and elaborated in \cite{SeqMonSqrdRes2010} to monitor squared residuals for a variance change of regression errors, as well as detectors which draw on classical results on boundary-crossing  probabilities for Brownian motion, specifically \cite{ErdoesKac1946} and \cite{RobbinsSiegmund1970}. 
A common basic feature is that the data from the monitoring period is successively cumulated, such that information increases as monitoring proceeds, and the unknown no--change second moment functional is replaced by a suitable estimator calculated from the training sample. The associated sequential stream of detector statistics is then compared with a boundary function to decide at each time instant whether or not a signal should be raised and the no-change null hypothesis be rejected. 

This paper contributes by establishing asymptotic theory for high--dimensional nonlinear and nonstationary time series of growing dimension, and we particularly address the case that the projection vector is estimated from the training sample. The results are based on strong Gaussian approximations of partial sums which extend previous results of \cite{MiesSteland2022} and \cite{StelandJMVA2020}. We show that any estimator with convergence rate $ O_\P(  \sqrt{\frac{r_d}{m} } ) $, $ r_d \ge 0 $, can be used, provided the dimension $d$ ensures $ \sqrt{d r_d} = O(m^\zeta)  $ for $ 0<\zeta<\frac12$, so that for $ \sqrt{m} $--consistent estimators the dimension $d$ can be almost as large as  $ \sqrt{m} $. For highdimensional settings where $ d $ is large compared to $m$ or even much larger, it is common to impose sparseness assumptions and employ sparse estimators. Examples include estimation of eigenvectors.  If  $\vecv $ is $ s$--sparse and the estimator $ \hat{\vecv}_m $ is concentrated on the support of $ \vecv $ with sufficiently large probability for large $m$, then the condition can be relaxed to $ \sqrt{s r_d} = O(m^\zeta)  $. Under sparsity, the dimension an grow almost exponentially, see our discussion in Section~\ref{Sec:Assumptions}. 

For the above two examples, Markowitz portfolios and deep learning networks, the optimal projection vector is a function of the unknown covariance resp. precision matrix. Thus, we study (soft--) thresholding estimates \color{\blau} and contribute consistency guarantees under the general assumptions imposed in this paper, thus generalizing the results of \cite{BickelLevina2008}. Especially, we contribute a non--probabilistic bound on the maximum norm of the difference between a thresholded matrix and its estimation target, which may be of independent interest. \color{black}

We elaborate on both applications discussed above, with a focus on neural net monitoring. Here, as a side--product of independent interest, we contribute an improved result on the Lipschitz constant of a (deep) neural network wth respect to the network parameters.

The organization of the rest of the paper is as follows. \color{\blau} Section~\ref{Sec:Methods} describes the framework and the proposed detection methods. \color{black} Asymptotic results \color{\blau} justifying the procedure and providing theoretical insights into its properties \color{black} are provided in Section~\ref{Sec:Asymptotics}. Section~\ref{Sec:EstimationOfV} studies (soft--) thresholding estimators of projections depending on the covariance matrix and the mean vector suitable for high--dimensional data.  Applications to deep learning and financial portfolio optimization are discussed in Section~\ref{Sec:Applications}, and the approach is illustrated by monitoring synthetic data \color{\blau} and investigated by simulations. \color{black} Proofs of the main results are provided in Section~\ref{Sec:Proofs}.

\section{Monitoring Framework}
 \label{Sec:Methods}

\color{\blau}
 In view of the fact that machine learning approaches for prediction are typically applied whenever new input data arrives, we consider an idealized setting where data vectors are sequentially observed, one after the other, and at each time point a $ \{ 0,1 \}$--valued detector is applied where $1$ stands for a signal indicating a change and the need for a reaction, such as retraining the deep learner, whereas $0$ means that we keep using the given neural network specification. \color{black} We  observe random vectors
\[
\vecY_1, \ldots, \vecY_m, \vecY_{m+1}, \cdots
\]
of dimension $d = d_m $, $m \in \N $, defined on a common probability space $ (\Omega, \calF, \P) $, where $ \vecY_1, \ldots, \vecY_m $ represents a training sample of size $m$, \color{\blau} and starting at time $ m+1 $ the observations $ \vecY_{m+i} $, $i \ge 1 $, arrive one after the other. Going beyond the classical assumption of independent observations, we consider below a general vector time series model allowing for serial correlations as arising in real--world applications, see below and Section~\ref{Sec:Assumptions} for details. 

\color{black} Our goal is to sequentially decide at each time $ n \ge 1 $ whether or not the second moment structure of the observations $ \vecY_{m+1}, \ldots, \vecY_{m+n} $ has changed in comparison to the training sample. We reduce dimensionality by projecting the data onto a vector $ \vecv \in \R^d $. Precisely, the detector proposed in this paper monitors the functional 
\[
\theta(\P_{\vecY_i}) = \vecv^\top [\E( \vecY_i \vecY_i^\top ) ]\vecv, \qquad i \ge 1,
\]
where $ \P_\vecY $ denotes the law of a random vector $ \vecY $. The  parameter $ \theta( \P_\vecY ) $ reacts to a change in the mean as well as to a change in the covariance matrix, if these changes are visible through the bilinear form induced by $ \vecv $ and do not cancel. There may be problems where methods are preferable which are able to explicitly decouple changes in the mean and the variances; here we refer to methods proposed by  \cite{RafaSteland2014}.  If $ \vecv$ is unknown, we replace it by an estimator calculated from the training sample as discussed in greater detail below, where we consider the classical non--$\ell_0$--sparse setting as well as $\ell_0$--sparse vectors and associated sparse estimators which typically yield much better convergence rates. 
 
\color{\blau} The sequential testing problem of interest is formalized as a at--most--one--change (AMOC), since our goal is to detect a change quickly to initiate measurements instead of finding and dating all change-points possibly present in the data. We test the no-change null hypothesis 
\[
	H_0 :  \vecv^\top \matM_i \vecv = \theta_0 \ \text{for all $i \ge 1$}
\] 
which asserts that the second moments -- viewed through the projection vector $ \vecv $ -- is constant and, especially, coincides with the second moments of the training sample, against 
the AMOC alternative hypothesis,
\begin{align*}
	H_A: &\ \text{there is $k^* \ge m \delta $ with $ \vecv^\top \matM_i \vecv = \theta_0 $, $ 1\le i < m + k^* $}, \\
	& \ \text{and $ \vecv^\top\matM_i \vecv = \theta_A$, $i \ge m+k^*$ },
\end{align*}
for some  $ \theta_A \not= \theta_0 $. Under $ H_A $ there is a change at time $k^*$ in terms of the second moments. In our treatment, $ \theta_0 $ and $ \theta_A $ may be unknown positive constants, namely the pre- and after--change values of the functional $ \vecv^\top \matM_i \vecv $,  and the change--point $k^* $ is assumed unknown and non--random. More general functional alternative hypotheses are considered in the next section. If the underlying vector time series is assumed to be centered, then $ \matM_t = \bfSigma_t $ and the above change--point problem considers a change in the variance of the projection $ \vecv^\top \vecY_t $. \color{black}

\color{\blau} Our time series model follows the widely used physical systems approach introduced by \cite{Wu2005}, which hosts many specific time series models as special cases, including autoregressive and moving average models, conditional heteroskedastic series and recursive stochastic models. \color{black} We assume that $ \vecY_t = (Y_{tj})_{j=1}^d $, $ t \ge 1 $, is a nonlinear and possibly nonstationary vector time series with mean vectors $ \bfmu_t = (\mu_{t1}, \ldots, \mu_{td})^\top $, second moment matrix $ \matM_t = \E( \vecY_t \vecY_t^\top ) $, covariance matrices $ \bfSigma_t = ( \Sigma_{t,ij} )_{i,j=1,\ldots, d} $ and coordinates given by
\[
  Y_{tj} = G_{tj}( \epsilon_t, \epsilon_{t-1}, \ldots), 
\]
for a measurable mapping $ \vecG_t = (G_{tj})_{j=1}^d : \R^\infty \to \R^d $ and i.i.d. $ U(0,1) $--distributed  random variables $ \epsilon_t  $. To formulate the physical weak dependence assumptions, let
\[
  \bfeps_t = (\epsilon_t, \epsilon_{t-1}, \ldots ), \qquad \bfeps_{t,j} = ( \epsilon_t, \ldots, \epsilon_{j+1}, \epsilon_j', \epsilon_{j-1}, \ldots )
\]
for $ j \le t $ and $ t \ge 1 $, where $ \{ \epsilon_t' \} $ is an independent i.i.d. sequence with $ \{ \epsilon_t' \} \stackrel{d}{=} \{ \epsilon_t \} $. For a vector $ \veca = (a_1, \ldots a_n)^\top \in \R^n $, $ \| \veca \|_r $  denotes the vector--$r$ norm, $ r \ge 1 $, and $ \| \veca \|_\infty = \max_{1 \le i \le n} |a_i| $ the maximum norm.

\subsection{Open-end monitoring} 
The open--end monitoring procedure essentially works as follows: At time  $ m+k $ of the monitoring period $ m+1, m+2, \ldots $ the current sum of  the most recent $k$ transformed observations is compared to an average calculated from the  training sample. To avoid that decisions are made based on too few data points of the monitoring period, the
 constraint $ k \ge m \delta $ for some (small) $ \delta > 0$ is imposed, so that monitoring starts at observation $ m(1+\delta) $ and continues either infinitely long or until a signal in favor of a change is given. 
 We focus on detectors using different strategies to threshold a suitable cumulated sum centered at the corresponding average in the learning sample. They are constructed such that the probability of a false alarm is controlled. Such a procedure has a significance level  comparable with the significance level of classical fixed sample Neyman--Pearson type hypothesis tests, which eases interpretation of results. But, clearly, the method differs from such a classical test in that it tries to come to the conclusion that the null hypothesis needs to be rejected as soon as the detector statistic provides enough evidence. This approach rules out classical methods studied in the field of statistical quality control which have the unpleasant property that in infinite monitoring an alarm is raised with probability one.

If we define the random matrices $ \matS_i = \vecY_i \vecY_i^\top $, then  $ \E( \matS_i ) = \matM_i $. Therefore, to detect a change \color{\blau}  at the $k$th time point, $m+k$, of the monitoring period, \color{black} we may consider for known $ \vecv $ the detector statistics
\[
   \left| \sum_{m < i \le m+k} \vecv^\top \matS_i \vecv - \frac{k}{m} \sum_{j=1}^m  \vecv^\top \matS_j \vecv \right|, \qquad k \ge \delta m.
\]
Since usually  $ \vecv $ is unknown, it is replaced by some estimator $ \hat{\vecv}_m $ computed from the training sample and safisfying Assumption A--iv (cf. next section), leading to the detector
\[
Q(m,k) = \left| \sum_{m < i \le m+k} \hat{\vecv}_m^\top \matS_i \hat{\vecv}_m - \frac{k}{m} \sum_{j=1}^m  \hat{\vecv}_m^\top \matS_j \hat{\vecv}_m \right|.
\]

All procedures signal at time $i \ge \delta m$ of the monitoring period and reject the null hypothesis of no change, if the trajectory of the standardized detector statistic, $ \{ Q(m,j) / \hat{\sigma}_{0,\infty} :  j \ge \delta m  \} $,  crosses a suitable multiple of a boundary function, $\{  c g(m,j) : j \ge \delta m \}$, at the $i$th time instant for the first time, where $ \hat{\sigma}_{0,\infty}^2  $ is an estimator of the asymptotic variance $ \sigma_{0,\infty}^2 $ discussed below. This means,
\begin{equation}
	\label{GammaDetector}
	\frac{Q(m,i)}{\hat{\sigma}_{0,\infty} g(m,i)} > c, \quad \text{but} \quad \frac{Q(m,j)}{\hat{\sigma}_{0,\infty} g(m,j)} \le c, \ m \delta \le j < i.
\end{equation}
$ c > 0 $ is a critical value, whose choice is discussed below. Note that the signal is given at the random stopping time 
\[
\tau_m = \inf\{ k \ge m\delta: Q(m,k) / \hat{\sigma}_{0,\infty} > c g(m,k) \},
\]
with the convention $ \inf \emptyset = \infty $. 

{\em Erd\"os-Kac detector:} The first procedure compares $ Q(m,k) $ with a multiple of the simple boundary function 
\[
  b(m, k) = m^{\frac{1}{2}} \left( \frac{m+k}{m} \right).
\]
As we will formally show, the probability that $ m^{-1/2} (1+k/m)^{-1}Q(m,k) / \sigma_{0,\infty} $ does not cross $ b(m,k) $  converges to $ \P( \sup_{\delta/(1+\delta)\le t \le 1} | B(t) | \le c ) $ for some Brownian motion $B$. Thus, comparing $ Q(m,k) $ with a multiple of the boundary $ b(m,k) $  asymptotically corresponds to comparing the absolute value of $B$ with the constant $c$. For small $ \delta $ one may select $c$ such that 
\[ \P\left( \sup_{t \in [0,1]} | B(t) | \le c\right) = \frac{4}{\pi} \sum_{i=0}^\infty \frac{(-1)^i}{2i+1} \exp\left(\frac{-(2i+1)^2 }{\pi^2/8c^2}\right)  = 1- \alpha,
\]
see \cite{ErdoesKac1946} and \cite[ch.~2]{ShorackWellner1986}.

{\em $\gamma$--detector:} The increasing variance of Brownian motion suggests to replace the latter constant by  a curved boundary function. The choice $ t^{1/2} $ would standardize $B$, but due to the law of the iterated logarithm the corresponding supremum is not well defined. Following \cite{HHKS2004} and generalizing the Erd\"os--Kac detector, the $ \gamma $--detector uses a multiple of the boundary
\[
  g(m, k) = m^{\frac{1}{2}} \left( \frac{m+k}{m} \right) \left( \frac{k}{m+k} \right)^\gamma 
\]
where $ \gamma $ is a tuning constant with $ 0  \le \gamma < 1/2 $.  For a late change $k^* / (m+k^*) \approx 1 $, so that the choice of $ \gamma $ has a negligible effect, but the closer the change to $m$, the smaller $ k^* / (m+k^*) $, so that large values of $ \gamma $ close to $ 1/2 $ are preferable for quick detection. 
Since monitoring may last forever and then analyzes infinitely many observations, we need to ensure that asymptotically the null hypothesis is never rejected with probability $ 1-\alpha $, if it is true. To achieve the latter,  we select $c$ such that
\[
\lim_{m \to \infty} \P\left( \sup_{k \ge \delta m} \frac{ Q(m,k)}{\hat{\sigma}_{0,\infty} g(m,k)} > c \right) = \alpha.
\]
Our theoretical results, see Section~\ref{Sec:Asymptotics}, justify to select $c$ as the (simulated) $1-\alpha $ quantile of the distribution function
$$
\mathcal{D}_\infty(x) = \P \left( \sup_{\frac{\delta}{1+\delta} \le t \le 1} \frac{|B(t)|}{t^\gamma} \le x \right), \qquad x \in \R, 
$$
where $ B(t) $ denotes a standard Brownian motion. Even for $ \delta = 0 $ explicit formulas for $ \calD_\infty(x) $ are not kown except $ \gamma = 0 $. 

{\em Robbins--Siegmund detector:} Lastly, we consider a   detector which does not require to simulate a critical value. It makes use of  the fact that  $|B(t)|$ crosses the boundary 
\[ 
	g_a(t) = \sqrt{(t+1)(a^2 + \log(t+1))}, \qquad a > 0,
\] 
for some $t > 0 $ with probability $e^{-a^2/2} $, i.e.,  $ \P( |B(t)| \ge g_a(t),  \exists t > 0 ) = e^{-a^2/2} $, see \cite{RobbinsSiegmund1970} and \cite{ChuStinchcombeWhite1996}. This result suggests to use $ a = \sqrt{2\log(\alpha)} $ to ensure a crossing probability of $ \alpha \in (0,1) $. For random processes converging weakly to a standard Brownian motion $B(t)$, $ t \ge 0 $, this boundary--crossing probability carries over. 
However, this result is not directly applicable to $ Q(m,k) $. But applying a suitable time transformation and interpreting the resulting rule in terms of the physical time scale $ k \in \N$ yields the following detection scheme: Raise a signal at time $ m+i $, $ i \ge \delta m$, if 
\begin{equation}
	\label{DetectorED}
	Q(m,i) > \hat \sigma_{0,\infty} g_a(m,i), \quad  g_a(m,i) :=
	\sqrt{m+i^*} \frac{m}{m-i^*} \sqrt{a^2+ \log \frac{m+i^*}{m} }
\end{equation}
for the first time, where $i^* = \lceil i \frac{m}{m+i} \rceil $.

The estimator $ \hat \sigma_{0,\infty} $ of the asymptotic variance $ \sigma_{0,\infty}^2 = \lim \text{Var}(  m^{-1/2} \sum_{i=1}^m \vecv^\top \matS_i \vecv ) $ can be obtained as follows. If the training data is an i.i.d. sample, use  the sample variance of $ \vecv^\top \matS_i \vecv $, $ 1 \le i \le m $. For time series data under Assumptions A--i -- A--iii of the next section we propose
\[
  \hat{\sigma}_{0,\infty}^2 = \frac{1}{m-b+1} \sum_{j=0}^{m-b} (S_j(b) - \bar{S}(b))^2,
\]
where $  S_j(b) = \frac{1}{\sqrt{b}} \sum_{i=j-b+1}^j \vecv^\top \matS_i \vecv $, $ 0 \le j \le m-b $, are subsampled scaled partial sums, $ \bar{S}(b) = \frac{1}{m-b+1} \sum_{j=0}^{m-b} S_j(b) $, and $b $ is a bandwidth parameter. Consistency of this estimator, which can be dated back to works of \cite{Carlstein1986} and \cite{PeligradShao1995},  has been shown in  \cite[Theorem~4.1]{MiesSteland2022}, provided Assumptions A--i and A--ii hold and the bandwidth is selected such that $ b = O(m^{\rho}) $ for some $ 0 < \rho < \frac{1}{2} $.

\subsection{Closed--end monitoring}

The above detection schemes can also be applied as closed--end schemes where monitoring stops when a prespecified time horizon parameterized as $ T m $, for some fixed $ T > 0 $, is reached. Let us briefly discuss the closed--end scheme for the $ \gamma $--detector. Here, one is interested in detecting a change point $ m\delta  \le k^* \le Tm $ and thus considers the stopping time
\[
  \tau_m^T = \min \{  m\delta \le k \le mT : Q(m,k) / \hat{\sigma}_{0,\infty} > c g(m,k) \}
\]
with the convention that $ \min \emptyset = \infty $ (no change detected). If the monitoring procedure examines the partial sums from $ \delta m $ to $ Tm $, one has  more freedom to choose a boundary function $ g(m,k) $ including case $ \gamma = 1/2 $. The decision threshold $ c $ is now selected as the (simulated)  $ 1- \alpha $ quantile of the distribution function
\[
  \mathcal{D}_T(x) =  \P\left( \sup_{t \in [\delta, T]} \frac{| B_{1}(t ) - t  B_{2}(1)|}{(1+t)\left(\frac{t}{1+t} \right)^\gamma} \le x \right), \qquad x \in \R,
\]
for two independent standard Brownian motions $ B_1 $ and $ B_2 $.

\subsection{Application to deep learners}

To monitor a deep learner to decide when it needs retraining, the above algorithms are applied with 
\[ \vecY_i \leftarrow  f_H( \vecX_i; \hat \bftheta_m ), \vecv \leftarrow \hat{\bfbeta}_m, \qquad \text{(detector $D$, monitoring of network predictions)}
\]
for $ i \ge 1 $, where $ (\hat{\bftheta}_m, \hat{\bfbeta}_m) $ are the trained connection weights of the deep learning. 

The above detector considers a single output neuron. For a network with $q$ outputs, one can apply a Bonferroni correction. Here, we monitor each output neuron separately, using a common critical value $ c $ corresponding to a nominal significance level $ \alpha / q $, and give a signal at time $ m +k $ if at least one of the $q$ detectors signals. If we denote by $ \tau_m(i)$ the time points when the $i$th detectors gives a signal, the overall detector signals at time $ \min( \tau_m(1), \ldots, \tau_m(q) )$.
\color{black}


\section{Asymptotics}
\label{Sec:Asymptotics}

\color{\blau}
The theoretical results justifying the proposed methods hold true under a general nonlinear time series model, whose assumptions are discussed in the detail in the next subsection. The results comprise Gaussian approximations of projected partial sums, on which the detectors rely and that is of independent interest, non--asymptotic bounds and convergence rates of estimators of means, second moments and covariance matrices, and the asymptotic theory of the proposed detectors. Here we provide the required asymptotic null distributions as well as asymptotics under local alternatives for the $ \gamma $--detector. The class of local alternatives consists of functional alternatives going beyond the classical case of jump-- or trend--like local alternatives.  \color{black} 

\subsection{Assumptions}
\label{Sec:Assumptions}

\color{\blau} We will now list and discuss all assumptions for the general nonlinear time series model $ \vecY_t = \vecG_t( \bfeps_t ) $, which substantially goes beyond typical settings assuming independent data. Besides classical moment assumptions, the assumptions limit the degree of dependence (A--i) and non-stationarity (A--ii), formalize the training sample (A--iii), and provide the required assumptions of proper estimators of the projection vector (A--iv). 
\color{black} 

\textbf{Assumption A--i} For some $ \beta > 2 $ and $ q > 4 $, $ \| \bfmu_t \|_\infty \le K $,
\begin{equation*}
	\label{A-i-a}
	\left( \E | G_{tj}( \bfeps_0 ) - \mu_{tj} |^q  \right)^{\frac{1}{q}} \le K,
\end{equation*}
\begin{equation*}
	\label{A-i-b}
	\left( \E | G_{tj}( \bfeps_t ) - G_{tj}( \bfeps_{t,t-\ell} ) |^q \right)^{\frac{1}{q}} \le C_1 \ell^{-\beta}, 
\end{equation*}
for $ 1 \le j \le  d $, $ \ell \ge 1 $ and $t \ge 1$, for constants $ K $ and $ C_1 $.

If $ \vecY_t = \vecY_t^{(m)}  $, $ t, m \ge 1 $, is an array indexed by  $m \in \N$, then the  constants $ K, C_1 $ arising in Assumption A--i  and all assumptions to follow are required to be universal, i.e., independent of  $d$ and $ m $.

The first condition of A--i requires existence of the central moments of order $4+\varepsilon $, $ \varepsilon >0 $ arbitrarily small. The second condition is a uniform coordinate--wise weak dependence assumption in the sense of Wu's physical dependence measure, so that the impact when replacing in $ Y_{tj} $ the lag $\ell$ innovation $ \epsilon_{t-\ell} $  by some i.i.d. copy is of magnitude $ \ell^{-\beta} $ in the $ L_q $--norm. 

\textbf{Assumption A--ii}  There exists some $ L > 0 $ such that the inequality
\[
V(\bar q ) = \max_{1 \le j \le d} \sum_{t=2}^n \left( \E | G_{tj}( \bfeps_0) - G_{t-1,j}( \bfeps_0 ) |^{\bar q} \right)^{\frac{1}{\bar q}} \le KL
\]
holds for $ \bar q = 4$, for all $ n \ge 1$. 

Assumption A--ii constrains the degree of nonstationarity by assuming that the
cumulated increments $ Y_{tj} - Y_{t-1,j} $ remain bounded when they are measured in the $ L_4 $--norm. This hold if the mapping $ t \mapsto G_{tj} $ is smooth in a Lipschitz sense. The following example serves as a \color{black} non--trivial \color{black} illustration and discusses the case of random curves, such as temperature in climate studies, which are discretely sampled at time instants $ t_k = \frac{k}{m} $, $ 1 \le k \le m $.

\begin{Example}  Suppose that the time series $ Y_{tj} $ is obtained by discretely sampling $d$ random curves $ G_j( u, \bfeps_t ) $,  $ u \in [0,\infty) $, $G_j : [0, \infty) \times \R^\infty \to \R $,  at equdistant time instants $ u \in \{ \frac{t}{m} : 1 \le t \le m\} $, such that $ Y_{tj} = G_{tj}( \bfeps_t ) = G_j( \frac{t}{m}, \bfeps_t ) $, $ 1 \le j \le d $. For example, $ G_j( \frac{t}{m}, \bfeps_t ) $ can be the temperature at day $j$ of year $t$, and the learning sample consists of $ m $ years of daily temperature data. Here, the influence of time, $t/m$, on the functional dependence of the randomness $ \bfeps_t $ on the temperature may be nonlinear.
	
	Let us assume that the functions $ G_j(u,\vecx) $ satisfy a uniform Lipschitz property
	\[
	|G_j(u_1,\vecx) - G_j(u_2,\vecx) | \le \operatorname{Lip}(G) |u_2-u_1| H(\vecx), \qquad \vecx \in \R^\infty, u_1, u_2 \in [0,\infty), 1 \le j \le d,
	\]
	for some Lipschitz constant $ \text{Lip}(G) $ and a function $H: \R^\infty \to \R $ with $ \E( H^4(\bfeps_0) ) < \infty $. Then
	\[
	\left( \E | G_{tj}( \bfeps_0) - G_{t-1,j}( \bfeps_0 ) |^4 \right)^{\frac{1}{4}}
	\le \operatorname{Lip}(G) \frac{\left( \E H^4(\bfeps_0) \right)^\frac{1}{4}}{m}, \qquad 1 \le j \le d, t \ge 1,
	\]
	such that A--ii holds with $ L = \text{Lip}(G) (\E H^4(\bfeps_0))^{1/4}/ K $.
\end{Example}

In our framework the first $m$ observations represent a training sample, also called  non--contamination period. In real world applications the training sample is sometimes carefully selected and assumed to satisfy much stronger assumptions than imposed here. Typical assumptions are that the training sample is a second--order stationary  vector time series, or that  it consists of i.i.d. random vectors of dimension $d$ with covariance matrix $ \bfSigma $, perhaps even assuming that $ \bfSigma^{-1/2} \vecY_i $ have i.i.d. (Gaussian) coordinates with finite fourth moments. These restrictive assumptions limit the applicability of results, since they are in practice often not fulfilled. They are included in the following assumption on the training sample. 

\textbf{Assumption A--iii} The training sample $ \vecY_t $, $ 1\le t \le m $, satisfies
the no--change (stability) hypothesis
\begin{equation}
	\label{equal-second-moments}
	\theta( \P_{\vecY_t} ) = \vecv^\top \matM_t  \vecv = \theta_0, \qquad  1\le t \le m,
\end{equation}
for a constant $ \theta_0 $,   there exists a measurable mapping $ \vecG_0 = ( G_{0j} )_{j=1}^d : \R^\infty \to \R^d $ such that for all $ 1 \le t \le m $,  $ 1 \le j, k \le d $ and $ h \ge 0 $
\begin{equation}
	\label{A-iii-b}
	\left|  \frac{\Cov( G_{tj}( \bfeps_t ), G_{tk}( \bfeps_{t+h} ) )}{\Cov( G_{0j}( \bfeps_t ), G_{0k}( \bfeps_{t+h} ) ) } - 1 \right| \le \frac{C_2}{m}.
\end{equation}
for some constant $C_2$, and  
\begin{equation}
	\label{A-iii-c}
	\inf_{i \ge 1 } \lambda_{\min}( \bfSigma_{i,\infty} ) > 0  \qquad \text{where} \qquad \bfSigma_{i,\infty}  = \sum_{h=-\infty}^\infty \Cov( \vecG_i( \bfeps_0 ) \vecG_i(\bfeps_h )^\top ). 
\end{equation}

Condition (\ref{equal-second-moments}) of A--iii requires that $ \theta( \P_{\vecY_t} ) $  attains a constant value in the training sample, which is required to be known, but it nevertheless allows for nonstationarity to some extent with constraints on the process mean and the covariance structure. 
Clearly, if $ \vecY_t $, $ 1 \le t \le m $, have mean zero, then $ \theta( \P_{\vecY_t} ) = \vecv^\top \bfSigma_t \vecv $ and \eqref{equal-second-moments} collapses to the null hypothesis that $ \Var( \vecv^\top \vecY_1 ) = \cdots = \Var( \vecv^\top \vecY_m) $. Further, a sufficient condition for constancy of $ \theta( \P_{\vecY_i}) $, $ 1 \le i \le m $, is to assume that
\[
\vecv^\top \bfSigma_t \vecv = \sigma_0^2, \qquad \vecv^\top \bfmu_t = m_0, \qquad 1 \le t \le m,
\]
for constants $ \sigma_0^2  \ge  0 $ and $ m_0 \in \R $.

Condition (\ref{A-iii-b}) of Assumption A--iii requires that the relative error when replacing any pseudo--lag $h$ cross--covariance of $ \vecY_t = \vecG_t( \bfeps_t ) $ with $ \vecG_t( \bfeps_{t+h} ) $ by the lag $h$ cross--covariance of $ \vecG_0( \bfeps_t ) $ is of the order $ O(\frac{1}{m} ) $. It limits the nonstationarity of the covariances but not necessarily those of other nuisance characteristics such as skewness or tail behavior. Condition (\ref{A-iii-b}) can be replaced by the assumption 
\begin{equation}
	\label{SuffA-iii}
	|\Cov( G_{tj}( \bfeps_t ), G_{tk}( \bfeps_{t+h} ) ) - \Cov( G_{0j}( \bfeps_t ), G_{0k}( \bfeps_{t+h} ) ) | \le \frac{k_h}{m}
\end{equation}
for some summable sequence $ k_h $, which can be formulated as
\begin{equation}
	\label{SuffA-iii-b}
	\|\Cov( \vecG_t(\bfeps_t), \vecG_t(\bfeps_{t+h} ) ) - 
	\Cov( \vecG_0( \bfeps_t), \vecG_0(\bfeps_{t+h}) ) \|_\infty 
	\le \frac{k_h}{m}.
\end{equation}
Condition (\ref{A-iii-c}) requires that the smallest eigenvalues of the matrices $ \bfSigma_{i,\infty} $ are bounded away from zero. Combined with (\ref{A-iii-b}) this will ensure that all summands of the cumulated sum arising in the CUSUM detector of interest have asymptotically the same scale, such that a well defined limiting law exists which is governed by a standard Brownian motion. If these assumption are omitted, then one could still construct Gaussian approximations using \cite{MiesSteland2022}, but then there may be no well defined asymptotic limit distribution of the detection scheme. Assumptions (\ref{A-iii-b}) and (\ref{A-iii-c}) ensure that, asymptotically, the procedure is governed by Brownian motion, and the asymptotic distribution is related to previous works on related detectors. 

The question arises whether the assumption on the second moments is compatible with nonstationarity and the assumed form of $ \vecY_t $.

\begin{Example} Consider the vector time series 
	\[
	  \vecY_t =  \sqrt{\frac{\nu_t-2}{\nu_t}} \tilde{\vecY}_t + e_t, \qquad t \ge 1,
	\]
	where $ e_t = \sum_{j=0}^\infty \rho^j \zeta_{t-j} $ for i.i.d. $\zeta_s \sim (0,1) $, $ s \in \Z $, and some $ \rho \in (-1,1) $, is a stationary autoregressive process of order 1 and $ \tilde \vecY_t \sim t( \bfmu, \bfSigma, \nu_t) $ with varying degrees of freedom specified as $ \nu_t =  \nu_0+ t $, $ t \ge 1$, for some $ \nu_0 > 4 $. Here, the multivariate $t$--distribution, $t(\bfmu, \Sigma, \nu) $, is defined as the law of $ \bfmu + \vecZ^* / \sqrt{Q/\nu} $ for independent random variables $ Q \sim \chi^2(\nu) $ and $ \vecZ^* \sim \calN( \vecnull, \bfSigma ) $. Its covariance matrix is $ \frac{\nu}{\nu-2} \bfSigma $.  Thus, $ \vecY_t $ has stationary first and second order moments. However, the higher order moments are nonstationary, since the mixed moments of a $ t(\vecnull,\bfSigma,\nu_0+t)$--distribution are given by
\[
\E( Y_{t1}^{n_1} \cdot \ldots \cdot  Y_{td_m}^{n_{d_m}}) =  [(t-1)/(t+1) ]^{k} 
(\nu_0+t)^{k/2} \frac{\Gamma( (\nu_0+t-k)/2)}{2^k\Gamma((\nu_0+t)/2)} \prod_{i=1}^{d_m} \frac{n_i!}{(n_i/2)!},  
\] 
for $  (n_1, \ldots, n_{d_m})^\top \in \N_0^{d_m} $, with $ k = \sum_{i=1}^{d_n} n_i $, if all $ n_i $ are even and $ \nu_0 + t > k $, and vanishes if at least one $ n_i $ is odd, see, e.g.,  \cite{Sutradhar1986}.
\end{Example}

More generally, the question arises whether a specific admissible moment structure for the whole vector process can be realized in the assumed form $ \vecG_t( \bfeps_t ) $. This question is related to two classical probems: Firstly, to the classical moment problem, which is not solved satisfactorily, especially in terms of algorithms, and, secondly, the problem to construct stochastic processes from prescribed marginals. However,  under regularity conditions addressing these two basic problems, such a construction is possible, even on the original probability space, provided it  carries uniform random variables $ \epsilon_t $, and this is a major motivation to study such nonlinear processes governed by i.i.d. innovations. For given real numbers $ \mu_{\bm n} $, $ \bm n \in \N_0^{k d_m} $, introduce the linear Riesz functional $ R$ on the space of polynomials in $kd_m$ variables via $ R( \bm x^{\bm n} )  = \mu_{\bm n} $  for monomials $  \bm x^{\bm n} = x_{11}^{n_{11}} \cdot \ldots \cdot x_{kd_m}^{n_{kd_m}} $, $ \bm n \in \N_0^{k d_m} $.

\begin{Proposition}
	\label{Construction}
	Suppose that $ \mu_{\bm t, \bm n}^{(m)}, \bm t \in \N^k, \bm n \in \N_0^{k d_m }, k \ge 1 $, are real numbers such that for all $ \bm t \in \N^k $, $ \mu_{\bm t, \bm n}^{(m)}, \bm n \in \N_0^{k d_m }, $ satisfy the Riesz--Haviland condition of  positive semidefiniteness of the Riesz functional, i.e, $ R(p) \ge 0 $ for all polynomials $ p$ with $ p(\bm x ) \ge 0 $ for all $ \bm x$, equivalently, 
	\begin{equation}
		\label{HamburgerCondition}
		\sum_{\bm n \in \N_0^{kd_m}} a_{\bm n} \mu_{\bm t, \bm n}^{(m)} \ge 0, \qquad \forall  \{ a_{\bm n} : \bm n \in \N_0^{kd_m} \} \subset \R.
	\end{equation}
	Then, for any $ 1 \le t_1 < \cdots < t_k $, $ k \in \N $, there exists a distribution $ Q_{t_1,\ldots, t_k}$ and a random  vector \[ (X_{t_11}, \ldots, X_{t_1d_m}, \ldots, X_{t_k1}, \ldots, X_{t_kd_m})^\top \sim Q_{t_1,\ldots, t_k} \] with the prescribed moments 
	$ \E( X_{t_11}^{n_{11}} \cdot \ldots \cdot X_{t_kd_m}^{n_{t_kd_m}}) = \mu_{(t_1, \ldots, t_k), (n_{11}, \ldots, n_{t_kd_m}) } $, for all $ \bm n = ( n_{11}, \ldots, n_{t_kd_m})^\top \in \N_0^{kd_m} $. If $ Q_{t_1, \ldots, t_k} $, $ t_1 < \cdots < t_k $, $ k \ge 1 $, satisfy Kolmogorov's consistency conditions, then one can construct, on a new probability space $ (\tilde \Omega, \tilde \calF, \tilde \P ) $, a $d_m $--dimensional    process, $ \tilde \vecX_t $, $ t \in \N $, whose finite dimensional distributions are given by the $ Q_{t_1, \ldots, t_k }$. Given i.i.d. $ \epsilon_t \sim U(0,1) $ defined on $ (\Omega, \calF, \P) $, one can construct, on $ (\Omega, \calF, \P ) $, a process $ \vecX_t $, $t  \ge 1 $,  of the form $ \vecX_t = \vecG_t( \bfeps_t ) $ for Borel functions $ \vecG_t $, with $ \{ \vecX_t : t \ge 1 \} \stackrel{d}{=} \{ \tilde\vecX_t : t \ge 1 \} $.
\end{Proposition}

Condition (\ref{HamburgerCondition}) cannot be weakened, since it is a necessary and sufficient condition to ensure that a probability measure with the given numbers as moments exists, see \cite{Haviland1936} and \cite{GhasemiEtAl2016} for recent extensions. The assumption $ \vecv^\top \matM_t \vecv = \theta_0 $, $ 1 \le t \le m $, puts a constraint only on the elements $ \mu_{(1, \ldots, m),\bm n} $ with $ \sum_{i=1}^{d_m} n_{ti} = 2 $ for $ 1 \le t \le m $, but the remaining  free parameters allow to specify nonstationary higher order moments. However, to best of our knowledge, specific algorithms to construct distributions with prescribed moments are constraint to dimension one (and hence to multivariate distributions with independent coordinates). For example, \cite{MeadPapaniolaou1984} and \cite{erdogmus2004minimax} study algorithms maximizing entropy and \cite{John2007} an approach using splines. Clearly, one can combine these methods with copulas, but this is certainly feasible in practice only to generate sequences of independent $d_m$--vectors. 
\color{black}

If $ \vecv $ is unknown, it is natural to estimate it from the training sample of size $m$. 
We consider the non--$ \ell_0 $--sparse case and {\em $s$--sparseness}. The latter means that the true projection $ \vecv = (v_1, \ldots, d)^\top $ has active set $ A = \{ j : v_j \not= 0 \} $ of cardinality $ s $.

%

\textbf{Assumption A--iv} Under $ s$--sparseness there is a sparse estimator $ \hat{\vecv}_m $ and  constants $ r_d \ge 1  $, such that , with $ \hat{A} = \{ j : \hat{v}_j \not= 0 \} $,
\[
 \P( \hat{A} \not = A )  =  O_\P\left( \sqrt{ \frac{s r_d}{m} } \right), \quad \| \hat{\vecv}_m - \vecv \|_2 = O_\P\left( \sqrt{ \frac{r_d}{m} } \right)
\]
and it holds 
$
\sqrt{s r_d} = O( m^\zeta ) 
$
for some $ 0 \le \zeta < \frac12 $. In the non--$ \ell_0$--sparse setting, the estimator $ \hat{\vecv}_m $ ensures the above convergence rate and it holds 
$
\sqrt{d r_d} = O( m^\zeta ) 
$.

This assumption includes the typical cases $ r_d = d $, $ r_d = d \log(d) $ and $ r_d = s \log(d) $. We discuss this issue and some special cases of interest in a separate section.
In the non--$ \ell_0 $--sparse setting, the dimension $d$ can grow at best almost at the rate $ \sqrt{m} $. This growth of the dimension can be attained when $ \vecv $ is $ \ell_1 $--sparse with coordinates satisfying a decay condition. 

Under $s$--sparseness the situation is much better. Especially, the constraint $ d = O( \sqrt{m} ) $ on the dimension is no longer required. If, for example, $ r_d = s \log d $, then the number of active variables can be almost of the order $ O( \sqrt{m} ) $. But for the more important case that the number $ s $ of relevant variables grows slowly or is bounded, the constraint on the dimension $d $ becomes very mild. For example, if $ r_d = s \log(d) $ with  $ s = O(1) $, the dimension can be of the order $ O( \exp(m^{2\eta} ) ) $, i.e., $d$ can grow almost exponentially in $m$. 

 
\subsection{Basic estimation under nonstationarity: Non--asymptotic guarantees}

In view of the general assumptions the following theorem provides non--asymptotic consistency guarantees (i.e., for all $m \ge 1 $) for sample moments. Let  
\[
	 \hat{\bfSigma}_m = \frac{1}{m} \sum_{t=1}^m (\vecY_t - \hat{\bfmu}_m)(\vecY_t - \hat{\bfmu}_m)^\top, \qquad  \hat{\matM}_m = \frac{1}{m} \sum_{t=1}^m \vecY_t \vecY_t^\top,
\] 
where $ \hat{\bfmu}_m = \frac{1}{m} \sum_{t=1}^m \bfmu_{t} $. Further put $ \bar\bfmu_m = \E( \hat{\bfmu}_m ),  \bar{M}_m = \E \hat{\matM}_m $ and $ \bar \bfSigma_m = \bar{M}_m - \bar{\bfmu}_m \bar{\bfmu}_m^\top $.  \color{\blau} In general, these population quantities are averages and depend on $m$, since we do not assume stationary data. For example, $ \bar \bfmu_m = \frac{1}{m} \sum_{i=1}^m \bfmu_i $, etc. Throughout, we indicate this by a bar--notation. \color{black} Denote the elements of $ \hat{\bfmu}_m $ by $ \hat{\mu}_{mj} $ and those of $ \hat{\bfSigma}_m $ by $ \hat{\Sigma}_{m,jk} $, etc. Let $ \zeta(z) = \sum_{j=1}^\infty j^{-z} $, $ z > 0 $, denote Riemann's zeta function. The following theorem summarizes basic consistency results of the estimators $ \hat{\bfSigma}_m, \hat{\bfmu}_m $ and $\hat{\matM}_m $. Threshold estimators of  covariance matrices are discussed in greater detail in the next section.

\begin{Theorem}
	\label{LemmaRatesMeanCovEst}
	Under Assumption A--i the sample mean $\hat{\bfmu}_m $, sample covariance matrix $ \hat{\bfSigma}_m  $ and sample second moments, $ \hat \matM_m $, are elementwise $ \sqrt{m} $--consistent in $ L_r $ and $ L_{\frac{r}{2}} $, respectively, for any $ 2 \le r \le q $, and hence in $ L_2 $, with
	\[
	\max_{1\le j \le d} \left( \E | \hat{\mu}_{mj} - \bar \mu_{mj} |^{r} \right)^\frac{1}{r} \le \frac{c \zeta(\beta)}{\sqrt{m}},
	\]
	\[
	  \max_{1 \le j,k \le d} \left(  \E | \hat{M}_{m,jk} - \bar M_{jk} |^\frac{r}{2} \right)^{\frac{2}{r}} \le \frac{c \zeta(\beta)}{\sqrt{m}},
	\]
	\[
	\max_{1\le j, k \le d} \left( \E | \hat{\Sigma}_{m,jk} - \bar \Sigma_{m,jk} |^{\frac{r}{2}} \right)^\frac{2}{r} \le \frac{c \zeta(\beta)}{\sqrt{m}},
	\]
	for some universal constant $ c$ (only depending on $q$) \color{\blau} and all $m \ge 1$. \color{black} 
\end{Theorem}

The above results provide the following convergence rates with respect to the maximum and Frobenius norm.

\begin{Corollary} 
\label{ConvRate}
	Under the assumptions of Theorem~\ref{LemmaRatesMeanCovEst} 
	\[
	  \| \hat{\bfmu}_m - \bar{\bfmu}_m \|_\infty = O_\P\left( \frac{d^\frac{2}{q}}{\sqrt{m}} \right), \ 
	  \| \hat{\bfSigma}_m - \bar{\bfSigma}_m \|_\infty = O_\P\left( \frac{d^\frac{4}{q}}{\sqrt{m}} \right), \  
	  \| \hat{\matM}_m - \bar{\matM}_m \|_\infty = O_\P\left( \frac{d^\frac{4}{q}}{\sqrt{m}} \right).  
	\]
	and
	\[
	\| \hat{\bfmu}_m - \bar{\bfmu}_m \|_2 = O_\P\left( \frac{d^{\frac{2-q}{q}}}{\sqrt{m}} \right), \ 
	\| \hat{\bfSigma}_m - \bar{\bfSigma}_m \|_F = O_\P\left( \frac{d^{2+\frac{4}{q}}}{\sqrt{m}} \right), \  
	\| \hat{\matM}_m - \bar{\matM}_m \|_F = O_\P\left( \frac{d^{2+\frac{4}{q}}}{\sqrt{m}} \right).  
\]	
\end{Corollary}

\subsection{A sequential Gaussian approximation}

The following theorem extends sequential Gaussian approximations for projected partial sums studied in \cite{MiesSteland2022} and \cite{MiesSteland2023} to estimated projections and provides approximations in terms of Brownian motion.  

\begin{Theorem} (Sequential Gaussian Approximation) \\
	\label{ThIP}
	Suppose that Assumptions A--i -- A--iii are satisfied and $ \| \vecv \|_1 \le C_\vecv $ for some constant $ C_\vecv $.  Let 	\[
	\xi(q,\beta) = \begin{cases}
		\frac{q-2}{6q-4}, & \beta \geq 3, \\
		\frac{(\beta-2)(q-2)}{(4\beta-6)q-4},& \frac{3+\frac{2}{q}}{1+\frac{2}{q}} < \beta < 3, \\
		\frac{1}{2}-\frac{1}{\beta}, & 2< \beta \leq \frac{3+\frac{2}{q}}{1+\frac{2}{q}}.\\
	\end{cases}     
	\]
	For some standard Brownian motion $ \{ B_t:  t \ge  0 \} $ the following assertions hold. 
	\begin{itemize}
		\item[(i)] 
		We have 
		\[ \max_{k \le m} \left| \sum_{i=1}^k \frac{ \vecv^\top (\matS_i - \matM_i) \vecv}{\sigma_{i,\infty}}  - B_k \right|  = O_\P\left( m^{\frac{1}{\nu}} \right) \]
		as well as
		\[ \max_{k \le m} \left| \sum_{i=1}^k \frac{ \vecv^\top (\matS_i - \matM_i) \vecv}{\sigma_{0,\infty}}  - B_k \right|  = O_\P\left( m^{\frac{1}{\nu}} \right) \]	
		with $ \nu = \frac{1}{\frac{1}{2} - \xi(q,\beta) + \varepsilon} > 2 $ for $ 0 < \varepsilon < \xi(q,\beta)  $,	where 
		\[
		\sigma_{i,\infty}^2 = \sum_{h \in \Z} \vecv^\top \Cov( \vecG_i( \bfeps_0 ),  \vecG_i( \bfeps_h ) ) \vecv, \ i \ge 0.
		\]
		\item[(ii)] If Assumption A--iv holds, then 
		\[ \max_{k \le m} \left| \sum_{i=1}^k \frac{ \hat{\vecv}_m^\top (\matS_i - \matM_i) \hat{\vecv}_m}{\sigma_{i,\infty}}  - B_k \right|  = O_\P\left( m^{\frac{1}{\nu}} \right) \]
		and
		\[ \max_{k \le m} \left| \sum_{i=1}^k \frac{ \hat{\vecv}_m^\top (\matS_i - \matM_i) \hat{\vecv}_m}{\sigma_{0,\infty}}  - B_k \right|  = O_\P\left( m^{\frac{1}{\nu}} \right) \]
		with $ \nu = \min\left( \frac{1}{\frac{1}{2} - \xi(q,\beta) + \varepsilon}, \frac{1}{\zeta + \varepsilon} \right) > 2 $ for some small $ \varepsilon > 0 $.
	\end{itemize}
\end{Theorem}

The first assertion of part (i) of Theorem~\ref{ThIP} refines \cite{MiesSteland2023}, where related Gaussian approximations  (not in terms of a Brownian motion) are derived under a slightly different assumption and used to study shrinkage covariance matrix estimators.

\subsection{Detector asymptotics under $H_0$ and local alternatives}

\color{\blau} The following theorem entails the detector's limiting distribution using the boundary function $ g(m,k) $, when no change is present. \color{black}

\begin{Theorem} 
\label{ThH0}
Suppose that Assumptions A--i -- A--iii are fulfilled and $ \| \vecv \|_1 \le C_\vecv $ for some constant $ C_\vecv $. If $ \vecv $ is estimated, additionally assume that A--iv holds, otherwise formally put $ \hat{\vecv}_m = \vecv$. Under the null hypothesis $ H_0 $ the detector statistic converges in distribution,
\[
  \lim_{m \to \infty} \P\left( \sup_{k \ge \delta m} \frac{Q(m,k)}{\hat{\sigma}_{0,\infty} g(m,k)} \le x \right) 
  = \P\left(  \sup_{\frac{\delta}{1+\delta} \le t \le 1} \frac{|B(t)|}{t^\gamma} \le x \right),
\qquad x \in \R,
\]
where $ B(t) $, $ 0 \le t < \infty $, is a standard Brownian motion. 
\end{Theorem}

\begin{Remark}
Note that the $ \calS_{\delta,\gamma} = \sup_{\frac{\delta}{1+\delta} \le t \le 1} \frac{|B(t)|}{t^\gamma} \to \calS_{0,\gamma} = \sup_{0 \le t \le 1} \frac{|B(t)|}{t^\gamma} $, $ \delta \to 0 $, so that for small $ \delta $ one can work with $ \calS_{0,\gamma} $. Explicit formulas for $ \calS_{0,\gamma} $ are, however, only known for $ \gamma = 0 $, but one can simulate quantiles or rely on quantiles provided in the literature.
\end{Remark}

For the Robbins--Siegmund detector, \eqref{DetectorED}, we have the following result.

\begin{Theorem} 
	\label{ErdoesKacDetector}
	Suppose that Assumptions A--i -- A--iii hold true. Then
	\[
	\lim_{m \to \infty}	\P\left( 
	Q(m,k) > \hat \sigma_{0,\infty}
	\sqrt{m+k^*} \frac{m}{m-k^*} \sqrt{a^2+ \log \frac{m+k^*}{m} }, \exists \delta m \le k < m
	\right) \uparrow \alpha,
	\]
	as $ \delta \downarrow 0 $, 
	where $k^* = \lceil k \frac{m}{m+k} \rceil $  and $ a = \sqrt{  2\log( 1/\alpha )} $. 
\end{Theorem}
\color{black}

Let us briefly discuss how these results are shown. The derivation of Theorem~\ref{ThH0} relies on the Gaussian approximation of Theorem~\ref{ThIP}. That result also yields the asymptotics of the related closed--end procedures using weighting functions of the form 
\[
  g(m,k) = g\left( \frac{k}{m} \right), \qquad g: [0,\infty) \to (0,\infty).
\]
One obtains, by the continuous mapping theorem, 
\[
  \max_{m \delta \le k \le mT} \frac{Q(m,k)}{ \hat{\sigma}_{0,\infty} g(k/m)} \stackrel{d}{\to} \sup_{\delta \le t \le T} \frac{|B_1(t) - t B_2(t) |}{g(t)}, 
\]
as $ m \to \infty $, where $ B_1, B_2 $ are two independent Brownian motions, and the fact that $ B_1(t) - t B_2(t) $ is equal in distribution to $ (1+t)B( t/(t+1) ) $ then leads to the result. The proof of Theorem~\ref{ErdoesKacDetector} applies a time transformation so that the latter process becomes a standard Brownian motion. Careful inspection and interpretation of the related  boundary crossing result for the time transformed $ Q(m,k) $ then yield the rule \eqref{DetectorED} and  Theorem~\ref{ErdoesKacDetector}. 

In the sequel we study asymptotic properties under alternatives and focus on the procedure \eqref{GammaDetector}. It is well known that the  $\gamma$--detector based on $ Q(m,k) $ has asymptotic power $1$ for a fixed alternative for the mean of the cumulated random variables. For very early changes $ k^* \ll \sqrt{m} $ the behaviour of such weighted CUSUMs for changes in the mean has been studied in \cite{AueHorvath2004} assuming independent errors. Here we want to study the behaviour under local functional alternatives which converge to $0$, if $m$ increases, and allow for a functional shape after the change point for the problem under investigation. For that purpose, we consider the closed--end monitoring setting where monitoring stops at the time horizon $ Tm $, for some fixed $ T >  0$. 

We formulate the local functional alternative model directly in terms of the moment functional and thus assume that for $ k \ge 1 $
\begin{equation}
\label{LocalFunctionalAlternative}
\vecv^\top \matM_{m + k} \vecv = \theta_0 + \frac{\Delta\left( \frac{k-k^*}{m} \right)}{\sqrt{m}} \eins_{k \ge k^*}
\end{equation}
for some integrable function $ \Delta : [0,T] \to \R $ which is right-continuous at $0$ and satisfies $ \Delta(0) \not= 0 $. The function $ \Delta(x) $ models the shape of the moment functional after the change-point $k^*$ which is assume to satisfy 
\[ 
  \frac{k^*}{m} \to \vartheta, \qquad m \to \infty,
\]
for some $ \vartheta \in (0, T) $.

\begin{Theorem} 
	\label{ThH1}
	Suppose that the assumptions of Theorem~\ref{ThH0} are fulfilled. If $ \Delta $ has bounded total variation on $ [0,T] $, then  under the sequence of local functional alternatives (\ref{LocalFunctionalAlternative})
	\[
	\lim_{m \to \infty} \P\left( \max_{\delta m \le k \le Tm} \frac{Q(m,k)}{\hat{\sigma}_{0,\infty} g(m,k)} \le x \right) 
	= \calD_T(x; \Delta) 
\]
where
\[
  \calD_T(x; \Delta ) = \P\left(  \sup_{t \in [\delta, T]} \frac{|B_1(t) - t B_2(1) + \eins_{t \ge \vartheta } \int_0^{t-\vartheta} \Delta(s) \, ds |}{(1+t)\left(\frac{t}{1+t} \right)^\gamma} \le x \right),  
	\]
	$x \in \R$, where $ B_1 $ and $B_2 $ are independent standard Brownian motions. 
\end{Theorem}


\section{Estimating projections depending on covariance and precision matrices}
\label{Sec:EstimationOfV}

As discussed in detail in the next section, in finance and deep learning one would like to select the projection by optimizing some criterion, which leads to an optimal projection, $ \vecv $, depending on the covariance or second moment matrix (or its inverse) and the mean. As a general and widely applicable approach we consider threshold estimators dating back to the works of \cite{BickelLevina2008} and \cite{RothmanLevinaZhu2009} for i.i.d. samples.  Here we show that the results can be generalized to nonlinear nonstationary vector time series. The results apply to covariance matrices as well as second moment matrices and their estimator. For brevity of presentation, we formulate the results for covariance matrices. Introduce the uniformity class
\[
\mathcal{U}( r, s_0, M, \varepsilon_0 ) = \left\{ \bfSigma : \Sigma_{ii} \le M, \sum_{j=1}^d |\Sigma_{ij}|^r \le s_0, 1 \le i \le d, \lambda_{\min}(\bfSigma) \ge \varepsilon_0 > 0 \right\},
\]
for some constant $ s_0 = s_0(d) $ and some $ 0 \le r < 1$. Matrices of this class have bounded variances and their spectrum is contained in $  [\varepsilon_0,  M^{1-r} c_0(d) ] $, since the general upper bound  $ \max_i \sum_j |\Sigma_{ij} | $ for the maximal eigenvalue  $\lambda_{\max}( \bfSigma ) $ is bounded by $ M^{1-r} s_0 $ for each $ \bfSigma \in \mathcal{U}( r, s_0, M, \varepsilon_0 )$. The thresholded estimator is defined by
\[
\hat{\bfSigma}_{th} = \calT_{t} \hat{\bfSigma}_m  = ( \hat{\Sigma}_{m,ij} \eins_{|\hat{\Sigma}_{m,ij} | \ge  t} )_{1 \le i\le d \atop 1 \le j \le d},
\]
for some threshold $ t \ge 0 $, where $ \hat{\bfSigma}_m $ is the sample covariance matrix. The operator $ \calT_t : \R^{d\times d} \to \R^{d \times d} $ performs a hard--thresholding operation. Soft--thresholding also nulls entries which are smaller in magnitude than $t$, but allows to shrink entries towards zero, if their magnitude is small but larger than the threshold. Any operator $ \calS_t : \R^{d \times d} \to \R^{d \times d} $ which operates element--wise, i.e., $ \calS_t \bfGamma = ( \calS_t( \Gamma_{ij} ) )_{ij} $,  $ \bfGamma = ( \Gamma_{ij} )_{ij} \in \R^{d \times d} $, and satisfies 
\begin{itemize}
	\item[(i)] $ | \calS_t( x ) | \le |x| $,
	\item[(ii)] $ \calS_t( x ) = 0 $ for $ |x| \le t $,
	\item[(iii)] $ |\calS_t( x) - x | \le t $,
\end{itemize} 
is called a soft--thresholding operator. Important special cases are hard--thresholding,  lasso soft--thresholding defined by
\[
  \calS_t^{l}(x) = \text{sign}(x) (|x|-t)_+,
\]
and SCAD thresholding, which is between hard--thresholding and lasso thresholding.
The resulting error when approximating some covariance matrix $ \bfSigma $ by a (soft--) thresholded version, $ \calS_t  \bfGamma  $, of an arbitrary  $d \times d $ matrix $ \bfGamma $ attains the following crucial bound with respect to the spectral norm $ \| \cdot \|_{op} $, when $ \bfSigma $ is constrained to $ \calU(r, s_0, M ) $.

\begin{Theorem}
	\label{ThBoundThresholdOp} Suppose that $ \bfSigma \in \calU( r, s_0, M) $. Then for any symmetric $d \times d$ matrix $ \bfGamma = ( \Gamma_{ij} )_{ij} $ the soft--threshold operator $ \calS_t $ satisfies the bound
	\[
	\|\calS_t {\bfGamma} - \bfSigma \|_{op} 
	\le 2 t^{1-r} s_0 + \| {\bfGamma} - \bfSigma \|_\infty
	\left( \#( |{\Gamma}_{ij} - \Sigma_{ij}| > (1-\gamma)t  ) + \frac{1}{(\gamma t)^r} s_0 + \frac{2}{t^r} s_0 \right)
	\]
	for any $ 0 < \gamma < 1 $. 
\end{Theorem}

The following theorem extends the result of \cite{BickelLevina2008} and \cite{RothmanLevinaZhu2009} for i.i.d. $ \vecY_1, \ldots, \vecY_m $ with uniformly bounded moments of order $q$ to the nonlinear nonstationary times series model satisfying Assumptions A--i and A--ii. 

\begin{Theorem} 
	\label{ThConsThresholdedCov}
	Under Assumptions A--i and A--ii the thresholded sample covariance matrix is consistent for $ \bar \bfSigma_m \in \mathcal{U}(r, s_0, M, \varepsilon_0 ) $ in the operator norm with
	\[
	\| \calS_t \hat{\bfSigma}_m - \bar\bfSigma_m \|_{op} = O_\P\left(\left( \frac{d^{\frac{4}{q}}}{\sqrt{m}} \right)^{1-r}\right) 
	\]
	if the threshold is selected as 
	\[
	t_{th} = C_{th} \frac{d^{\frac{4}{q}}}{\sqrt{m}}
	\]
	for some large enough constant $ C_{th} $ and $ d = o\left(m^{\frac{q}{8} } \right) $.
\end{Theorem}

\color{\blau}
A related problem is estimation of the precision matrix $ \bar \Phi_m = \bar \bfSigma_m^{-1} $. Here one may invert the thresholded estimator  $ \calS_{t_m} \hat{\bfSigma}_m  $ and thus define the estimator
\[
  \hat{\Phi}_m = ( \calS_{t_m} \hat{\bfSigma}_m  )^{-1}.
\]
This allows estimation of projection vectors which are smooth functions of the precision matrix $ \bar \Phi_m $ and the expectation $ \bar\bfmu_m $ as required when computing optimal portfolios.

\begin{Theorem} 
\label{ThConsPrecision}
	Under the assumptions of  Theorem~\ref{ThConsThresholdedCov} we have
\[
	\|\hat\Phi_m - \bar\Phi_m \|_{op}
	= O_\P\left(\left( \frac{d^{\frac{4}{q}}}{\sqrt{m}} \right)^{1-r}\right).
\]
A projection given by $ \vecv = \vecv_m = f( \bar \bfSigma_m^{-1}, \bar\bfmu_m ) $, for some Lipschitz continuous function $f$, can be estimated by
\[
\hat{\vecv}_m = f( (\calS_{t_m} \hat{\bfSigma}_m )^{-1}, \hat{\bfmu}_m ), 
\]
at the same convergence rate provided the estimator $ \hat{\bfmu}_m $ attains the same convergence rate. 
\end{Theorem}

\color{black}

It is worth discussing the rate of convergence obtained in the above results and compare it with the classical parametric rate of convergence $ O(1/\sqrt{m}) $. For small $r$ the convergence rate is almost $ O_\P\left( \frac{r_d}{\sqrt{m}} \right) $ with $ r_d = d^{\frac{4}{q}} $, and for $ q \to \infty $ one obtains almost the parametric rate $ O_\P(\frac{1}{\sqrt{m}} ) $.

\section{Applications and Simulations}
\label{Sec:Applications}

\subsection{Deep neural networks}

We consider  a (deep) artificial neural network with $H$ hidden layers for a $d$--dimensional input $ \vecx \in \calD \subset \R^{d} $ and a univariate output $ y $ defined recursively as  
\[
  \vecf_1(\vecx) = \sigma_1( \matW_1 \vecx ) , \qquad 
  \vecf_j( \vecx ) = \sigma_j( \matW_ j\vecf_{j-1}( \vecx ) ), \quad 2 \le j \le H,
\]
for weight matrices $ \matW_j \in \R^{n_j \times n_{j-1}} $  and (nonlinear) activation functions $ \sigma_j $, $ 1 \le j \le H $, where $ n_0 = d $ is the input dimension and $n_j \in \N $ is the width of the $j$th layer, $ 1 \le j \le H$. The linear output layer linearly combines the features generated by the $H$ layers and computes the output
\[
  y = f_{net}( \vecx ) =  \bfbeta^\top f_H( \vecx ),
\]
with $ \bfbeta \in \R^{n_H} $. 
For multivariate output of dimension $ n_o $  one takes a coefficient matrix $ \matB \in \R^{n_o \times n_H} $ and calculates $ \vecy = \matB f_H( \vecx  ) $. The network can be compactly written as
\[
  \vecy = \matB \sigma_H^{\matW_H} \circ \cdots \circ \sigma_1^{\matW_1}( \vecx ),
  \qquad \sigma_j^{\matW_j}(\vecx ) = \sigma_j( \matW_j \vecx ), \qquad \vecx \in \R^{n_{j-1}}.
\]
And if we put $ \sigma_{H+1}(x) = x $ and $ \matW_{H+1} = \matB $, then 
\[
  \vecy = \sigma_{H+1}^{\matW_{H+1}} \circ \cdots \circ \sigma_1^{\matW_1}( \vecx ).
\] 
For simplicity of presentation, however, we consider a univariate output in what follows.

The activation functions, $ \sigma_j : \R \to \R $, act elemenwise on vectors. Common choices of the activation functions, $ \sigma $, such as the ReLU given by $ \max(x,0) $, the sigmoid $ \frac{1}{1+e^{-x}} $ or the tanh activation $ \frac{e^x - e^{-x}}{e^x+e^{-x}} $, are Lipschitz continuous and satisfy $ \sigma(0) = 0 $. Indeed, Lipschitz continuity has beneficial effects for deep learning. Whereas the latter choices are differentiable, the widely used ReLU activation is not differentiable at $0$, where it has the subdifferential $ [0,1] $. However, one may approximate it by $ \sigma_k(x) = \frac{1}{2k} log( 1 + e^{2kx} ) $ for some large $k>0$. 

Denote by $ \bftheta = ( \matW_1, \ldots, \matW_H ) $ the parameters of the hidden layers; the full Euclidean parameter is then $ \bfeta = (\text{vec}( \bftheta), \bfbeta ) $. Ocassionally, we write $ \vecf_j( \vecx; \bftheta ) $ and $ f_{net}( \vecx; \bfbeta, \bftheta ) $ to stress the dependence on the weights.

The following result establishes the Lipschitz property of a deep neural network, which plays an important role in the generalization capabilities of a network, \cite{BartlettFosterTelgarsky2017}. The result is similar to \cite[Prop.~6]{Lederer2023}, but our proof is much simpler and the Lipschitz constant is typically  smaller (note that typically $ L_H \ll 1 $). For bounds on the Lipschitz constant for specific hidden layers and a discussion of related generalzation bounds see \cite{GoukEtAl2021}. Define for $ 1 \le k \le H $
\[ 
\bftheta( \wt{\matW}_k ) = ( \wt{\matW}_1, \ldots, \wt{\matW}_k, \matW_{k+1}, \ldots, \matW_H ), 
\] and $ \bftheta( \wt{\matW}_0 ) = (\matW_1, \ldots, \matW_H) $. 

\begin{Proposition}  
	\label{LipschitzDeepLearner}
	Consider a deep learning network given by $ \sigma_H^{\matW_H} \circ \cdots \circ \sigma_1^{\matW_1}( \vecx ) $ with  $ \rho_i $--Lipschitz continuous activation functions $ \sigma_i: \R \to \R $ satisfying $ \sigma_i(0) = 0 $, $ 1 \le i \le H $. Suppose that $ \calD \times \Theta $ is bounded or $ \sigma_H $ is bounded. Then 
	\[
	\| \vecf_j( \vecx; \bftheta(\wt{\matW}_k) ) - \vecf_j( \vecx; \bftheta( \matW_k ) ) \|_2 \le L_{jk}  \| \vecx \|_2 \| \wt{\matW}_k - \matW_k \|_{op},
	\]
	for $ 1 \le k \le H $, where $   L_{jk} = \prod_{i=1, i \not= k}^j \rho_i  \|\matW_i \|_{op}   \rho_k $.
	
	The neural network has the Lipschitz property
	\[
	\| \vecf_H( \vecx; {\bftheta} ) - \vecf_H( \vecx; \wt{\bftheta} ) \|_2 
	\le L_H \sqrt{H} \| \vecx \|_2 \| \wt{\bftheta} - \bftheta \|_F ,
	\]
	where $ L_H = \max_{1 \le k \le H} L_{Hk} $. 
\end{Proposition}

The deep learning regression model assumes that before a change--point, $ k^* $,
\[
   Y_t = \bfbeta_0^\top \vecf_H( \vecX_t; \bftheta_0 ) + e_t, \qquad t \le k^*, 
\]
for  true parameters $ \bfbeta_0, \bftheta_0 $ (see below), and a zero mean noise process $ \{ e_t \} $. 

Whereas there is a vast of literature on performance guarantees for neural networks, results about parameter--consistency are quite limited. Under ideal conditions, nonlinear least squares estimators of fixed--topology neural networks are consistent and asymptotically normal, \cite{WhiteHalbert1992}. We confine here to this setting, mainly because recent and ongoing  work on statistical guarantees for networks with multiple layers is not yet settled, although the results of \cite{LohWainwright2015} for penalized $M$--estimators under restricted strong convexity (RSC) are potentially applicable to wide classes of deep networks, see also the recent result about RSC for the least squares loss of \cite{BanerjeeEtAl2022}, where local minima in layer--wise $ \ell_2 $--balls around a random initialization are considered. Further, for the perceptron with additive i.i.d. noise \cite{pmlr-v48-yangc16} have shown that the $ \ell_1 $--penalized least squares estimator, $ \hat{\bfbeta}_m $, attains under suitable conditions with high probability the optimal rate of convergence  $O\left( \sqrt{\frac{s\log d}{m} } \right) $, where $ s $ is the number of active entries of the true parameter.
 
Let us now fix a multiple layered network as above and consider minimization of the
least squares loss in the training sample,
\[
   \ell_m(\bfeta) = \frac{1}{m} \sum_{i=1}^m \ell( (Y_i, \vecX_i); \bfeta), \qquad \ell((y,\vecx), \bfeta) = (y - f_{net}( \vecx; \bfeta ))^2, y \in \R, \vecx \in \R^d,
\] 
over a compact set $ \calB \times \Theta $ around a random initialization. 

In practice, an iterative numerical solver is applied which ideally converges to the  local minimum of the error landscape closest to the initialization. Common algorithms used in practice are classical backpropagation, stochastic gradient descent, \cite{RobbbinsMonro1951},  and ADAM, \cite{Kingma2015}, which is the state--of--art algorithm from a practical point of view due its superior running time. Let us assume that the numerical solver succeeds in finding a stationary point of $ \ell_m( \bfeta ) $ located in $ \calB \times \Theta $, and thus yields an estimator $ \hat \bfeta_m = ( \hat \bftheta_m, \hat{\bfbeta}_m ) $ of the true value $ \bfeta_0 = (\bftheta_0, \bfbeta_0 ) = \operatorname{argmin}_{ (\bfbeta, \bftheta)\in \calB \times \Theta} \E \ell_m(\bfbeta,\bftheta) $, assumed to be an interior point and unique. Suitable conditions for smooth activation functions such that weakly consistent solutions exist, a.s., can be formulated in terms of the estimating function $ G_m(\bfeta) = \nabla \ell_m(\bfeta) $, see \cite{JacodSoerensen2018}:
\begin{itemize}
	\item[(i)]		$ G_m(\bfeta_0) \stackrel{\P}{\to} 0$, $ m \to \infty $.
	\item[(ii)]  $ G $ and all $ G_m $, $m \ge 1 $, are a.s. differentiable on $  M = \calB \times \Theta $ and $ \nabla G_m $ converges uniformly in probability,
	\[
	 	\sup_{\bfeta \in M} \| \nabla G_m(\bfeta) - \nabla G(\bfeta) \|_2 \stackrel{\P}{\to} 0.
	\]
	\item[(iii)] The Hessian satisfies $ \nabla^2 G(\bfeta_0) > 0 $.
\end{itemize}
If, additionally, the estimating function satisfies a central limit theorem,
\[
  \frac{1}{\sqrt{m}} \sum_{i=1}^m \nabla \ell( (Y_i,\vecX_i), \bfeta_0)  \stackrel{d}{\to} N( \mathbf{0}, \matV( \bfeta_0 )), 
\]
as $ m \to \infty $, for some matrix $ \matV(\bfeta_0) >  0 $, then the estimator $ \hat{\bfeta}_m $ attains the convergence rate $ O_\P( \frac{1}{\sqrt{m}} ) $. If the training sample $ (Y_1, \vecX_1), \ldots, (Y_m,\vecX_m) $ forms a strictly stationary nonlinear vector time series of the form $ \vecH( \bfeps_t ) $ satisfying Assumptions A--i -- A--ii, then the above conditions (i) and (ii) typically hold, especially if $ Y_t $ and $ \vecX_t $ are bounded and the activation functions are $ C^2 $, since $  M = \calB \times \Theta $ is bounded, cf. Theorem~\ref{LemmaRatesMeanCovEst}. A central limit theorem for the estimating function $ G_m(  \bfeta ) = \nabla \ell_m( \bfeta ) $ holds as well, where the positive definiteness of the limiting covariance matrix depends on the neural network and the dependence structure of the training sample and thus needs to be assumed, cf. Theorem~\ref{ThIP}.

Denote $ \vecZ_i(\bfbeta) = \vecf_H( \vecX_i; \bftheta) $, $ 1 \le i \le m $. Observe that $ \bfbeta_0 $ is a solution of the normal equations of the output layer, 
\[ 
	\bfSigma_Z \bfbeta_0 = \E( \vecZ_1(\bftheta_0) Y_1 ), \qquad  \text{with\ } \vecZ_i(\bfbeta) = \vecf_H( \vecX_i; \bftheta).
\] 
It attains the explicit representation $ \bfbeta_0 = \bfSigma_Z^{-1} \E( \vecZ_1(\bftheta_0) Y_1 ) $ if $ \bfSigma_Z > 0 $. $\hat\bfbeta_m $ solves the corresponding sample normal equations, 
\[
\hat{\bfSigma}_{Zm}(\hat{\bftheta}_m)  \hat\bfbeta_m
=  \frac{1}{m} \sum_{i=1}^m \vecZ_i( \hat \bftheta_m ) Y_i
\]
where $ \hat{\bfSigma}_{Zm}( \hat{\bftheta}_m ) = \frac{1}{m} \sum_{i=1}^m \vecZ_i( \hat\bftheta_m ) \vecZ_i( \hat \bftheta_m)^\top $.  Thus,  under the assumptions discussed above, $ \hat{\bfbeta}_m = \bfbeta_0 + O_\P( \frac{1}{\sqrt{m}} ) $ with $ \bfbeta_0 = \bfSigma_Z( \bftheta_0 )^{-1} \E( \vecZ_1(\bftheta_0) Y_1 ) $. The optimal weights encode the second moment structure of the last hidden layer and its covariances with the responses.

For a new input $ \vecX^* $ the true network prediction is $ \bfbeta_0^\top \vecf_H( \vecX^*; \bftheta_0 ) $ and the prediction of the trained network is $ \hat{\bfbeta}_m^\top \vecf_H( \vecX^*; \hat{\bftheta}_m ) $. As long as the projection  $ \bfbeta_0^\top \vecf_H( \vecX^*; \bftheta_0 ) $ does not change, there is no need to retrain the network. But, contrary, if that projection changes, which will be indicated sooner or later when cumulating information contained in the stream $f_H( \vecX_{m+i}; \hat{\bftheta}_m )$, $ i \ge 1 $, the network predictions may no longer be valid. This leads to the  detector, $ D $, with 
\begin{align*} 
		\vecY_{m+i} &\leftarrow f_H( \vecX_{m+i}; \hat{\bftheta}_m ), \qquad i \ge 1,  \\
		\hat{\vecv}_m &\leftarrow \hat{\bfbeta}_m,
\end{align*} 
in order to detect a change in the regressors pointing to the necessity to retrain the network. 

\subsection{Financial portfolio optimization}

In financial portfolio optimization one aims at determining a portfolio vector $ \vecw $ of the proportions of money to invest in each asset, such that $ \sum_j w_j = 1 $.  $ w_j > 0 $ is a long position and $ w_j < 0 $ a short position in asset $j$. If a portfolio has many short positions, the gross--exposure $ \|\vecw \|_1 $ may be large. It is now well understood that portfolios should have bounded gross--exposure $ \|\vecw \|_1$. However, the celebrated mean--variance Markowitz theory, see \cite{Markowitz1952, Markowitz1959} and \cite{Sharpe1964}, which establishes the well known relationship between investment risk measured by the variance and mean return, does not guarantee bounded gross-exposure. It selects $ \vecw$ by minimizing the variance $ \Var( \vecw^\top \vecY ) $ of the portfolio return under the constraint  $ \vecw^\top \bfmu = \mu_0 $ of a target portfolio return $ \mu_0 $, where $ \vecY $ is a $d$--vector of asset returns with mean $ \bfmu $ and covariance matrix $ \bfSigma > 0 $. The optimal solution, $ \vecw^* $, is the unique portfolio realizing the target return $ \mu_0 $ and having the smallest variance possible, typically with short positions. It is a linear combination of $ \bfSigma^{-1} \bfmu $ and $ \bfSigma^{-1} \bm{1} $, so that
the results of our Section~\ref{Sec:EstimationOfV} can be applied.

The oracle risk, $ R( \vecw^* ) $,  the empirical risk, $ R_m( \vecw^* ) $ and the associated estimation error $ | R_m( \vecw^* ) - R(\vecw^*) |$, where $ R( \vecx ) = \vecx^\top \bfSigma \vecx $ and $ R_m( \vecx  ) = \vecx^\top \hat{\bfSigma}_m \vecx  $ with an arbitrary estimator $ \hat{\bfSigma}_m $ of $ \bfSigma $,  are, however, not controlled by the Markowitz solution. For large dimension the error of the risk estimate may suffer from substantial cumulation of estimation errors, especially due to the need to estimate $ d(d-1)/2 $ covariances. One popular and theoretically sound approach to mitigate this problem is to minimize the Markowitz mean--variance function $ M = \vecw^\top \bfmu + \lambda \vecw^\top \bfSigma \vecw $, for some risk--aversion parameter $ \lambda > 0 $, under a gross--exposure constraint $ \| \vecw \|_1 \le c $. Then, for any portfolio $ \vecw $, the sensitivity of $M$ can be bounded by $ \| \hat{\bfmu}_m - \bfmu \|_\infty \|\vecw \|_1 + \lambda \| \hat{\bfSigma}_m - \bfSigma \|_\infty \| \vecw \|_1^2 $, and this bound holds for arbitrary estimators of $ \bfmu $ and $ \bfSigma $. It follows that if one can estimate all entries of $ \bfmu $ and $ \bfSigma $ equally well, as guaranteed by our Corollary~\ref{ConvRate}, the sensitivity is controlled by the maximum estimation error and the gross--error exposure $ \|\vecw\|_1 $. This phenomenon has been observed in the finance literature by \cite{JagMa2003} for the choice $ c = 1 $ yielding the optimal no--short--sales (i.e., long--only) portfolio and has been recently discussed in greater detail by \cite{FanZhangYu2012}. Further, constraining $ \| \vecw \|_1 $ is equivalent to shrinking of the estimated covariance matrix, see \cite{FanZhangYu2012} and \cite{ZhaoLedoitJiang2021} for a recent discussion and further references. 

Even if $ \vecw $ is determined by another approach, for example by tracking an index by a sparse portfolio leading to a sparse portfolio vector $ \vecw $, i.e., with $ \| \vecw \|_1 \le c $ or even with finite $ \| \cdot \|_0 $--norm, as done by many exchange traded funds (ETFs), the portfolio variance $ \vecw^\top \bfSigma \vecw $ is a crucial risk functional, and one is faced with the following problem: If the portfolio risk changes, the weights should be updated.  For a large universe of assets, however, the frequent precise (re--) calculations of $ \vecw $ is a challenging and costly issue from a numerical and statistical viewpoint as well as from a practical viewpoint, since this may lead to massive trading costs and can even be infeasible for highly captitalize funds, which therefore  tend to update weights only twice a year. If the ETF implements a certain strategy, e.g.  momentum, more frequent updates triggered by a suitable monitoring procedure could, however, lead to more profitable short--term positions. 

Given a sequence of returns $ \vecY_1, \vecY_2, \ldots $ denote by $ \hat{\vecv}_m $ an estimator of an unknown true (preferably bounded--gross--exposure) portfolio $ \vecv $, estimated from the learning sample $ \vecY_1, \ldots, \vecY_m $. Recall that the risk of the unknown portfolio $ \vecv $ is given by $ R(\vecv) = \vecv^\top \bfSigma \vecv $. As long as the distribution of the returns changes without altering the true risk $ R(\vecv) $, there is no need to update the portfolio,  exactly the changes affecting $ R(\vecv) $ are of interest. Under Assumption A--iii and the common zero mean assumption for financial returns under normal market conditions, unbiased estimates of $ R(\vecv) $ are given by $ \vecv^\top \vecY_i \vecY_i^\top \vecv $, $ 1 \le i \le m$. A change of the expectation of $ \vecv^\top \vecY_{m+k} \vecY_{m+k}^\top \vecv $, $ k \ge 1 $, indicates the necessity to update the porfolio. This motivates to apply the proposed detector, $D$, which monitors the sequence $ \hat{\vecv}_m^\top \vecY_{m+k} \vecY_{m+k}^\top  \hat{\vecv}_m $, $ k \ge 1 $.

\subsection{Experimental support and simulations}

\textbf{Experiment:} To check whether the proposed approach works we applied it to synthetic data following the model
\[
  Y_t = \left\{ 
   \begin{array}{ll}
   	\cos( 10 \pi x_t   ) + e_t, \qquad & 1 \le t < 2,000, \\
   	\cos( 10 \pi x_t^*  ) + e_t, \qquad & 2,000 \le t < ,4000, \\
   	\cos( 4 \pi x_t^*  ) + e_t^*, \qquad & 4,000 \le t,
   	\end{array}
   \right.
\]
where $ x_t \sim U(-1/2,1/2) $, $ x_t^* \sim U(0,1) $, $ e_t  \sim N(0,0.1) $ and $ e_t^* \sim N(0,0.05 ) $. At $ k^*_1 = 2,000 $ the distribution of the regressor changes and at $ k_2^* = 4,000 $ the regression coefficient changes. Two independent uniformly distributed regressors were added to mimic the situation that the regression function depends only on a subset of the regressors. For this data a neural network needs to be trained three times: Once at the beginning and then after the two change--points. 

A neural network with three fully--connected hidden layers consisting of $4 $, $2$, and $1$ neurons, respectively, relu activations for the hidden layers and a linear output layer was trained using the first $ m = 1,000 $ observations using ADAM for $ 100 $ epochs with batch size 32, and a validation split of the training data of  $0.2$. For the simulation study below training was based on $ 30 $ epochs. The fitted network was then used for prediction and only retrained if a signal was given. The procedure was rolled over after a signal as follows: All data before the signal were discarded, and the next $m$ observations were used as a fresh training sample from which the neural network was retrained. Then, the detector was restarted after the training sample. To avoid non--termination of the detector, a sample of length $ T=50,000 $ was simulated and the stopping time was set to $ \infty $ if there was no signal up to $T$. The experiment was implemented in R using tensorflow and keras on a Mac M1 Max with GPU support. 

\begin{center}
	\begin{figure}
		\begin{center}	\includegraphics[width=8cm, bb=0 0 500 400]{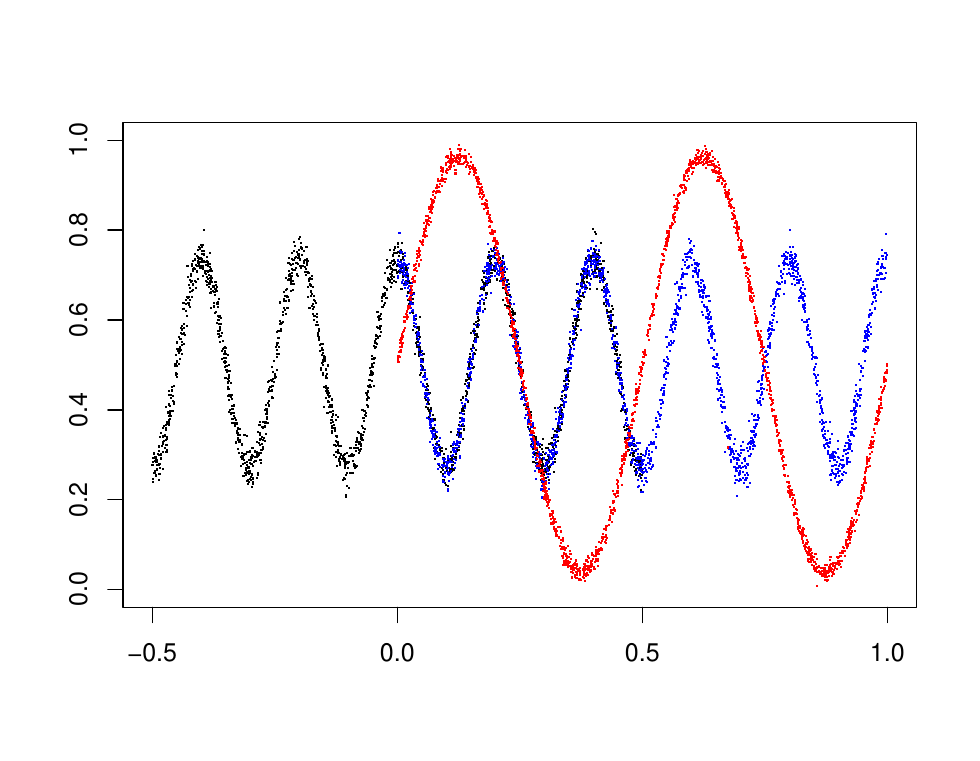} \end{center}
		\includegraphics[width=7.5cm, bb=0 0 500 400]{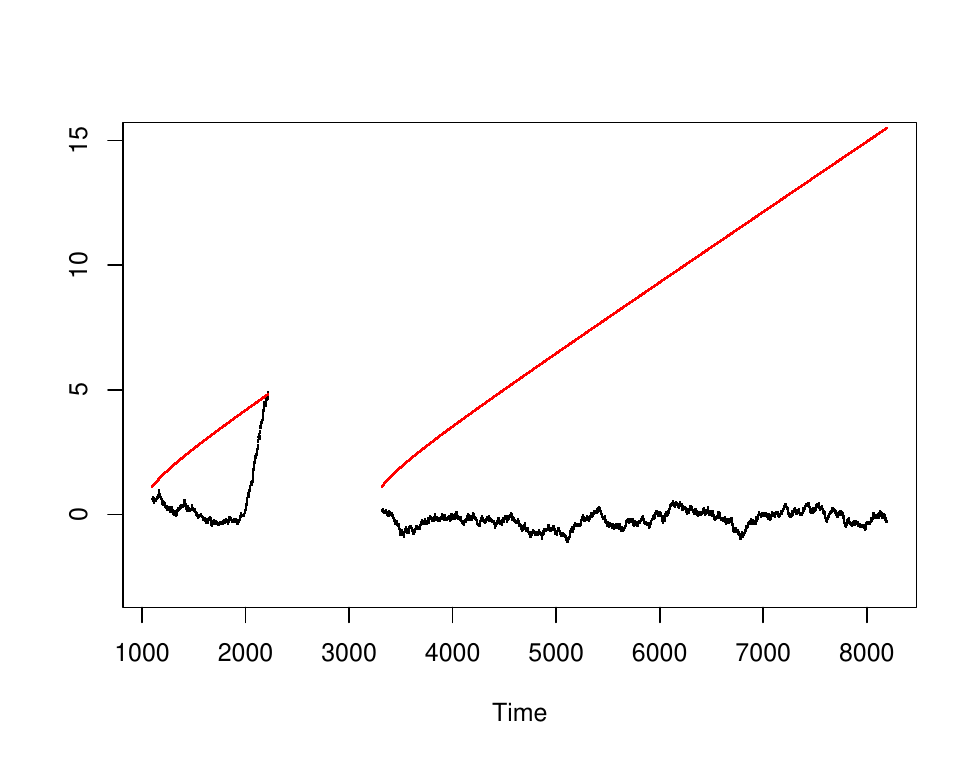}
		\includegraphics[width=7.5cm, bb=0 0 500 400]{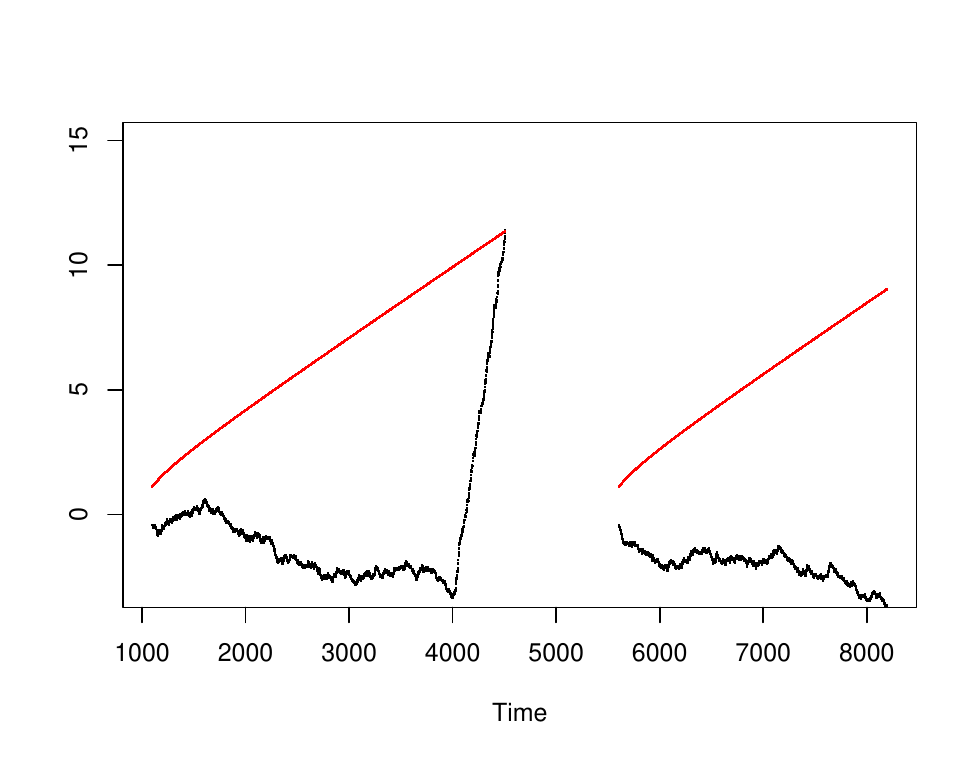}
		\caption{Experiments: The top panel shows the synthetic data set. Before the first change data points are black, after the first and before the second one blue, and data after both changes are in red. The  left (right) panel depicts the detector statistic $ D_{proj,k}$ ($D_{res,k}$) and the boundary $ B_k $. The method reacts quickly and there are no false alarms.}
		\label{Fig1}
	\end{figure}
\end{center}

Figure~\ref{Fig1} depicts the synthethic data ($Y_t$ scaled to the unit interval), simulated according to the above model, and, for the first 10,000 data points, the scaled detector statistic, $D_{proj,k} = \frac{1}{\sqrt{m}} Q(m,k) $, and the associated scaled  boundary, $ B_{k} = c_\alpha g(m,k) $, where the constant $ c_\alpha $ is selected to ensure an asymptotic false detection rate $ \alpha $. In additiona, a residuals detector, $ D_{res,k} $, was run obtained from  $ D_{proj,k} $ by replacing all $ \hat{\bfbeta}_m^\top \matS_t \hat{\bfbeta}_m  $ by $ Y_{m+\ell} - \hat{\bfbeta}_m^\top \matS_t \hat{\bfbeta}_m  $, $ t \ge 1 $. This detector is designed to detect a change in the regression coefficients, e.g. when using a different neural network. For the simulated synthetic data the detector $D_{proj,k}$ gives a signal at obs. $2,117$ corresponding to a delay of $ 117 $ and no further signal. The detector $ D_{resid,k} $ signals a change at obs. $4,508$ (delay: $508$) and gave no further signal as well. Both detectors react quickly to the changes and do not suffer from false signals. We can summarize that  for this example the proposed detection algorithm operates optimally and  minimizes the required computational costs of network training by invoking it only three times. 

\textbf{Simulations:} A simulation study was conducted to investigate the statistical behaviour of the proposed detectors in terms of the probability to get i) a signal (to reject $H_0$) at all and ii) to get a correct signals (i.e., after the change--point). A further characteristic of interest is the delay when a signal was given. For that purpose, both detectors were applied to samples simulated according to the above model for varying tuning parameter $ \gamma $ which controls the weighting function. The simulation model has one change--point of interest, $ k_1^*$, which should be detected by the proposed prediction detector $D$, and one change--point for which $D$ should be agnostic. 
 Analyzing the simulation results reveiled that indeed the detector $D$ only reacts for the change--point $ k_1^* $ and gives a signal at $ k_2^* $ only in rare cases. Therefore results for $ k_2^* $ are not reported. 
 


The empirical rejection rates and average delays were simulated under the null hypothesis and under  parameterized alternatives $ H_1(f) $, for $ f \in \{ 0.05, 0.1, 0.15 \} $, where the concept drift change--point $ k_1^*$,  at which the regressor distribution changes from a uniform law on $ [-1/2,1/2] $ to a uniform law on $ [-1/2+f,1/2+f] $, is located at $ k_1^* = (1+f) m $. That is, samples 
\[
Y_t = \left\{ 
\begin{array}{ll}
	\cos( 10 \pi x_t   ) + e_t, \qquad & 1 \le t < (1+f) 1,000, \\
	\cos( 10 \pi x_t^*  ) + e_t, \qquad & (1+f) 1,000 \le t < ,4000, \\
	\cos( 4 \pi x_t^*  ) + e_t^*, \qquad & 4,000 \le t,
\end{array}
\right.
\]
were generated, where $ x_t \sim U(-1/2,1/2) $, $ x_t^* \sim U(-1/2+f,1/2+f) $, $ e_t  \sim N(0,0.1) $ and $ e_t^* \sim N(0,0.05 ) $. Again, two additional uniform regressors were added and fed into the neural network. At $ k^*_1 = 2,000 $ the distribution of the regressor changes and at $ k_2^* = 4,000 $ the regression coefficient changes. 

The results are provided in Table~\ref{Tab2}. All methods perform quite well in terms of the type I error rate, but the Robbins--Siegmund detector turned out to result in quite conservative type I error rates, which leads to a loss of power to detect alternatives. This observation is in agreement with simulation results obtained in \cite{ChuStinchcombeWhite1996}, althought the boundary was used there in a quite different setting. To allow better comparisons, the value $a$ was selected to match the type I error rate of the $ \gamma $--detector, studied for and $ \gamma \in \{ 0.25, 0.45, 0.49 \} $, when $ \gamma = 0.49 $. In terms of the power to detect the concept drifts parameterized by $f$, the results also indicate that for larger changes the methods are quite similar. But for small changes the $ \gamma $--detector with $ \gamma = 0.49 $ (i.e., close to $1/2$) outperforms the competitors. Except rare cases, the signals are correct, i.e., after the change--point. However, in terms of the delay the procedures differ. Both the Erd\"os--Kac detector and the Robbins--Siegmund detector have much larger delays than the $ \gamma $--detector, and this holds even for large changes where the detection rate is around $ 80\% $. One can summarize the findings as follows: The Robbins--Siegmund detector is easy to apply but operates on a smaller--than--nominal type I error rate when using the theoretical value for $a$. The $ \gamma$--detector with $ \gamma $ close to $ 1/2 $ has high power and leads to quick signals compared to its competitors.



\begin{table}
	\begin{center}
	\begin{tabular}{|c|c|cccc|} \hline
		$ \gamma $ & Level & $ f $   & Power &  Correct Signal & Delay \\ \hline
		\multicolumn{6}{|c|}{Erd\"os--Kac detector}\\ \hline
		& 7.6 & 0.05 &  54.7 & 1 & 1041 \\
 				&& 0.10 & 75.2 & 1 & 662   \\
				&& 0.15 &  82.6  & 1 &  626 \\ \hline
		\multicolumn{6}{|c|}{$\gamma$--detector}\\ \hline
		0.25& 10.1 & 0.05 &  58.3 & 1 & 995 \\
				&& 0.10 & 76.0 & 1 & 594   \\
				&& 0.15 &  82  & 99.9 &  463 \\
		0.45 &9.7 & 0.05 &  53.5 & 1 & 809 \\
		         & & 0.10 & 75.4 & 99.3 & 464   \\
		         & & 0.15 &  82  & 98.3 &  394 \\
		0.49& 8.7 & 0.05 & 63.6 &1 &  609 \\
		        && 0.10 &  76.5 & 98.9 &  514 \\
		        & & 0.15 & 80.6 & 98.3 &  385  \\ \hline
		 \multicolumn{6}{|c|}{Robbins-Siegmund detector} \\ \hline
		 &  8.7  & 0.05 & 56.5 & 1  &  937 \\
		 &     & 0.10 & 74.8 & 1  &  688 \\
		 &     & 0.15 & 83.5 & 1  &  559 \\ \hline
	\end{tabular}
	\end{center}
	\caption{Simulated level and power under alternative $ H_1(f) $,  the conditional probabilties that a signal is correct, i.e., given after the change $ k_1^* $, and the mean delay.}
	\label{Tab2}
\end{table}

\section{Proofs}
\label{Sec:Proofs}

\begin{myproof}{of Proposition~\ref{Construction}} Assumption (\ref{HamburgerCondition}) yields the existence of a probability measure $ Q=Q_{t_1, \ldots, t_k} $ on $ \R^{kd_m} $, not necessarily unique, such that $ \int x_{11}^{n_{11}} \cdot \ldots \cdot x_{kd_m}^{n_{kd_m}} \, d Q = \mu_{\bm t, \bm n} $, for any $ \bm t = (t_1, \ldots, t_k)  \in \N^k $ with pairwise different entries and all $ \bm n \in \N_0^{kd_m} $.  By Kolmogorov's extension theorem, one can construct on a new probability space a stochastic process $ \tilde{\vecX}_t $, $ t \ge 1 $, whose finite dimensional distributions with respect to $ t_1, \ldots, t_k $, are given by $ Q_{t_1, \ldots, t_k} $, for all time points $ 1 \le t_1< \cdots < t_k $ and all $k \in \N $. The explicit construction of an equivalent process on the original probability space is as follows. Suppose we already have constructed $ \sigma = (\vecX_1, \ldots \vecX_{t-1}) $ with $ \vecX_i = \vecG_i(\bfeps_i) $, $ 1 \le i \le t-1$. Apply \cite[Lem.~1, p.212]{Billingsley1999} with  $ \nu = Q_{1, \ldots, t} $, $ \mu = \calL(\sigma) = Q_{1, \ldots, t-1} $ and $  U=\epsilon_t \sim U(0,1)$ to obain 
	a random variable $ \tau $ defined on $ (\Omega, \calF, \P) $, which is a function of $ (\vecX_1, \ldots, \vecX_{t-1}) $ and $ \epsilon_{t} $, such that $ (\vecX_1, \ldots, \vecX_{t-1}, \tau) \sim \nu $. Thus, we can take $ \vecX_{t} = \tau( \vecX_1(\bfeps_1), \ldots, \vecX_{t-1}(\bfeps_{t-1}), \epsilon_t) =: \vecG_t( \bfeps_t ) $. 
	It remains to construct $ \vecX_1 = (X_{11}, \ldots, X_{1d_m})$. Use $ \epsilon_1 $ as a seed to define i.i.d. uniform random variables $ \epsilon_{11}, \ldots, \epsilon_{1d_m} $ (which are functions of $ \epsilon_1 $). Put $ X_{11} =  X_{11}( \epsilon_{11}) :=  F_{Q_{11}}^{-1}( \epsilon_{11} ) $, where $ F_{Q_{11}}(x) = F_{Q_{1}}( x, \infty, \ldots, \infty) $ is the 1st marginal d.f. of $ Q_{1} $. Next, apply iteratively \cite[Lem.~1, p.212]{Billingsley1999}:  Having already constructed $ \sigma = (X_{11}, \ldots, X_{1,i-1}) $, put $ \nu = Q_{11, \ldots, 1i} $ (the marginal distribution of $ Q_1 $ with respect to the first $i$ coordinates) and $ U = \epsilon_{1i} $, to obtain $ \tau $, a function of $  (X_{11}, \ldots, X_{1,i-1})  $ and $ \epsilon_{1i}$ with $ (\sigma,\tau) = (X_{11}, \ldots, X_{1,i-1},\tau) \sim \nu = Q_{11, \ldots, qi} $. Thus, put $ X_{1i} = \tau( X_{11}(\epsilon_{i1}), \ldots, X_{1,i-1}(\epsilon_{1,i-1}), \epsilon_{1i} ) $. This gives us $ \vecX_1 = (X_{11}, \ldots, X_{1d_m}) $, a function $ \vecG_1( \bfeps_1 ) $ of $ \bfeps_1 = (\epsilon_1, \epsilon_0, \ldots ) $. 
\end{myproof}
\color{black}

The proofs of Theorems~\ref{ThIP} and \ref{ThH0} rely on several auxiliary results.

\begin{Lemma} 
	\label{lemma1}
	Suppose that $ \xi_i' $, $ i \ge 1 $, are mean zero and unit variance random variables defined on a probability space  $ (\Omega', \calF', \P') $ on which there exists a standard Brownian motion $ B_1' $, such that with $ S_n' = \sum_{i=1}^n \xi_i' $ the strong approximation
	\begin{equation}
		\label{SApproxBasic}
		| S_n' - B_1'(n) |  = o( n^{1/\nu} ), \qquad a.s.,
	\end{equation}
	holds for some $ \nu > 2 $. Define the standard Brownian motion $ B_{2,m}'(t) = B_1'(m+t) - B_1'(m) $, $ t \ge 0 $, which is independent of $ B'_1(t) $, $ 0 \le t \le m $. Then for any $ \delta > 0 $
	\[
	\sup_{k \ge \delta m} k^{-1/\nu} | S'_{m+k} - S'_m - B_{2,m}'(k) | = O_\P(1),  \qquad m \to \infty.
	\]
\end{Lemma}

\begin{myproof}{of Lemma~\ref{lemma1}}
	Observe that 
	\begin{align*}
		k^{-\frac{1}{\nu}} | S'_{m+k} - S'_m - B_{2,m}'(k) |
		& \le 
		\left( \frac{m+k}{k} \right)^{\frac{1}{\nu}} \frac{1}{(m+k)^{\frac{1}{\nu}}} | S'_{m+k} - B'_1(m+k) | \\
		& \qquad + \left( \frac{m}{k} \right)^{\frac{1}{\nu}} \frac{1}{m^{\frac{1}{\nu}}} | S'_{m} - B'_1(m) |. 
	\end{align*}
	The first term can be bounded by
	\begin{equation}
		\label{BoundBySup}
		\frac{1+ \delta}{\delta} \sup_{k \ge \delta m} (m+k)^{-\frac{1}{\nu}} | S'_{m+k} - B'_1(m+k) | \le \frac{1+\delta}{\delta}  \sup_{N \ge 1}  N^{-\frac{1}{\nu}} | S'_N - B'_1(N) |. 
	\end{equation}
	By assumption, there exists an event $ A $ with $P(A)=1 $ such that for each $ \omega \in A $ the random variable  $ Z_N(\omega) =  N^{-\frac{1}{\nu}} | S'_N - B'_1(N) |(\omega) $ is $o(1) $ and hence bounded. Thus, $ C(\omega) = \sup_{N \ge 1} Z_N(\omega) < \infty $. We obtain with probability one
	\[
	\sup_{k \ge \delta m} \left( \frac{m+k}{k} \right)^{\frac{1}{\nu}} \frac{1}{(m+k)^{\frac{1}{\nu}}} | S'_{m+k} - B'_1(m+k) | \le \frac{1+\delta}{\delta} C
	\]
	Since the random variable $ \frac{1+\delta}{\delta} C $ is stochastically bounded, the first term is $ O_\P(1) $.
	The second term can be bounded by
	\[
		\sup_{k \ge \delta m} \left( \frac{m}{k} \right)^{\frac{1}{\nu}} 			\frac{1}{m^{\frac{1}{\nu}}} | S'_{m} - B'_1(m) | \le \frac{1}{\delta^\frac{1}{\nu} }  \frac{1}{m^{\frac{1}{\nu}}} | S'_{m} - B'_1(m) | = o_\P(1).
	\] 
	This proves the assertion. 
\end{myproof}

The following lemma checks for several series derived from $ \vecY_t $ that they satisfy Assumptions A-i and A-ii with $V(2)$ and determines appropriate powers and constants.

\begin{Lemma} 
	\label{WeakDepLeftProj}
	Suppose that Assumptions A--i and A--ii hold. 
	\begin{itemize}
		\item[(i)] The time series $ \vecv^\top \matS_i \vecv - \vecv^\top \matM_i \vecv  $, $ i \ge 1 $, satisfies Assumption A--i with power $ q/2 $ instead of $ q $,
		 	$4 K^2 \|\vecv\|_1^2 $ instead of $K$,  and $ 4 K \|\vecv \|_1^2 C_1 $ instead of $ C_1 $. A--ii holds with $ 4K^2 \|\vecv \|_1^2 L $ instead of $K L $ for $ V(2) $.
		\item[(ii)]  The vector time series $ \vecv^\top \vecY_t \vecY_t^\top $, $ t \ge 1 $, satisfies Assumption A--i with power $ q/2 $ instead of $q$,  $ 2 K^2 \|\vecv\|_1 $ instead of $K$, and $ 2K \| \vecv \|_1 (\| \vecv \|_1 + 1 ) C_1 \sqrt{d} $ instead of $C_1$. Especially,
		\begin{equation}
			\label{PhysDepMeas}
			\left( \E \left\| \vecG_t(\bfeps_t) \vecv^\top \vecG_t(\bfeps_t ) - \vecG_t(\bfeps_{t,t-j} ) \vecv^\top \vecG_t(\bfeps_{t,t-j} ) \right\|_2^{\frac{q}{2}} \right)^{\frac{2}{q}} 
			\le 2K \| \vecv \|_1 (\| \vecv \|_1 + 1 ) C_1\sqrt{d} j^{-\beta}.
		\end{equation}
		Assumption A-ii holds  for $ V(2) $ with $ 4 \|\vecv\|_1^2 K^2L $ instead of $ KL $. 
		\item[(iii)] $ \left( \E \| \vecY_t \|_2^q \right)^{\frac{1}{q}} \le  K^{\frac{1}{q}} \sqrt{d} $. 
		\item[(iv)] $ \left( \E \| \vecG_t( \bfeps_t ) - \vecG_t( \bfeps_{t,t-j} ) \|_2^q \right)^\frac{1}{q} \le C_1^{\frac{1}{q}} \sqrt{d} j^{-\beta} $, $ j \ge 0 $.
		\item[(v)] The (vectorized) series $ \vecY_t \vecY_t^\top $ satisfies  Assumption A--i with power $ q/2 $ and constants $Kd $ and $ C_1 d $.
	\end{itemize}
\end{Lemma}

\begin{myproof}{Lemma~\ref{WeakDepLeftProj}} 
	As a preparation, let $ X, Y \in L_q $ be random variables with expectations $ \mu_X $ and $ \mu_Y $. The decomposition
	\[
		X Y = (X-\mu_X)(Y- \mu_Y) +  (X-\mu_X ) \mu_Y + \mu_X ( Y - \mu_Y) + \mu_X \mu_Y 
	\] implies the bound $ \| X Y \|_{L_{r}} \le \| X - \E X \|_{L_q} \|Y - \E Y \|_{L_q}  + |\E Y | \| X - \E X \|_{L_q} + | \E X |\| Y - \E Y \|_{L_q} + | \E X | | \E Y | $ for $1 \le r \le q $. Therefore,  $ | \E Y_{tj} Y_{tk} | \le 4 K^2 $ under Assumption A-i, such that   $\| \matM_t \|_\infty \le 4 K^2 $, since $ q > 4 $. Further $ \| \bfSigma_t \|_\infty = \max_{1 \le j,k \le d} | \E(Y_{tj}-\mu_{tj})(Y_{tk}-\mu_{tk}) |\le K^2 $, and $ \max_{1 \le j \le d} \| Y_{tj} \|_{L_r} \le \max_{1 \le j \le d} \| Y_{tj} - \mu_{tj} \|_{L_r} + K  \le 2 K $, for any $ 2 \le r \le q $.
		
	Let us first show (ii). Rearranging the above decomposition and denoting $ \Sigma_{XY} = \Cov(X,Y) $ and  $ M_{XY} = \E(XY) = \Sigma_{XY} + \mu_X \mu_Y $, one has the representation
	\[
	 \overline{XY} = XY - M_{XY} = (X-\mu_X)(Y-\mu_Y) + (X-\mu_X) \mu_Y + \mu_X(Y-\mu_Y) - \Sigma_{XY}.
	\]
	Thus,
	\begin{align*}
		\| \overline{XY} \|_{L_{q/2}} 
		& \le \| X -\mu_X \|_{L_{q}} \| Y - \mu_Y \|_{L_{q}} + | \mu_Y | \| X - \mu_X \|_{L_{q/2}} + |\mu_X| \| Y-\mu_Y \|_{L_{q/2}} + | \Sigma_{XY} |.
	\end{align*}
	Applying this inequality with $ X = Y_{tk} $ and $ Y = \vecY_t^\top \vecv $ gives for the $k$th coordinate $ Y_{tk} \vecY_t^\top \vecv $ of $ \vecY_t \vecY_t^\top \vecv  $, $ 1 \le k \le d $,
	\begin{align*}
		\|Y_{tk} \vecY_t^\top \vecv - \E(Y_{tk} \vecY_t^\top \vecv ) \|_{L_{q/2}}
		& \le 
		\|Y_{tk} - \mu_{tk} \|_{L_{q}} \| (\vecY_t - \bfmu_t)^\top \vecv \|_{L_{q}} 
		+ \|\bfmu_t \|_\infty \| \vecv \|_1 \| Y_{tk} - \mu_{tk} \|_{L_{q/2}} \\
		& \qquad + \| \bfmu_t \|_\infty  
		\| \vecv^\top (\vecY_t - \bfmu_t ) \|_{L_{q/2}} + | \Cov( Y_{tk}, \vecv^\top \vecY_t ) |.
	\end{align*}
 	To estimate the terms on the right side, notice that for any $ 2 \le r \le q $ 
	\[
	\left\| (\vecY_t - \bfmu_t)^\top \vecv  \right\|_{L_r}
	\le \sum_{\ell=1}^d | v_\ell | \| Y_{t\ell} - \mu_{t\ell} \|_{L_{r}}  \le K \| \vecv \|_1
	\]
	and
	\[
	| \Cov( Y_{tk}, \vecv^\top \vecY_t ) | \le \sqrt{ \E(Y_{tk}-\mu_{tk})^2 } 
	\sqrt{ \E\left( \vecv^\top (\vecY_t - \bfmu_t) \right)^2 }  \le  K^2 \| \vecv \|_1.
	\]
	Therefore,
	\[
	\| Y_{tk} \vecY_t^\top \vecv - \E(Y_{tk} \vecY_t^\top \vecv ) \|_{L_{q/2}}
	\le  4 K^2 \|\vecv\|_1, 
	\]
	Observe that by the triangle inequality and Jensen's inequality, for any $ 1 \le r \le q $,
	\begin{align*}
		\E |\vecv^\top \vecG_t( \bfeps_t ) |^r 
		& \le \|\vecv \|_1^r \E \left| \sum_{j=1}^d \frac{|v_j|}{\| \vecv \|_1 } | G_{tj}( \bfeps_t ) |   \right|^r \\
		& \le \|\vecv \|_1^r \sum_{j=1}^d  \frac{|v_j|}{\| \vecv \|_1 }  \E | G_{tj}(\bfeps_t ) |^r \\
		& \le \|\vecv \|_1^r \max_{1 \le j \le d} \E | G_{tj}( \bfeps_t ) |^r,
	\end{align*}
	so that
	\begin{equation}
	\label{BoundOnInnerProduct}
		  \| \vecv^\top \vecG_t( \bfeps_t )  \|_{L_r} \le \| \vecv \|_1 \max_{1 \le j \le d} \| G_{tj}( \bfeps_t ) \|_{L_r} \le 2 K \| \vecv \|_1.
	\end{equation}
    Using this fact, the lag $j$ (coordinate--wise) physical weak dependence measure can now be estimated by
	\begin{align*}
		& \left( \E | G_{tk}(\bfeps_t) \vecv^\top \vecG_t(\bfeps_t ) - G_{tk}(\bfeps_{t,t-j}) \vecv^\top \vecG_t(\bfeps_{t,t-j} ) |^\frac{q}{2} \right)^{\frac{2}{q}}  \\
		& \quad \le \left( \E | (G_{tk}(\bfeps_t) - G_{tk}(\bfeps_{t,t-j})) \vecv^\top\vecG_t(\bfeps_t ) |^\frac{q}{2} \right)^{\frac{2}{q}}
		+ \left(  \E  | G_{tk}( \bfeps_{t,t-j} )( \vecv^\top \vecG_t(\bfeps_{t} ) - \vecv^\top \vecG_t(\bfeps_{t,t-j} ) ) |^{\frac{q}{2}} \right)^{\frac{2}{q}} \\
		& \quad \le \left( \E |\vecv^\top \vecG_{t}( \bfeps_t ) |^q \right)^{\frac{1}{q}} 
		\left( \E | G_{tk}( \bfeps_t ) - G_{tk}( \bfeps_{t,t-j} ) |^q \right)^{\frac{1}{q}}  
		\\
		& \qquad + \left( \E | G_{tk}( \bfeps_{t,t-j} ) |^q \right)^{\frac{1}{q}}
		\left( \E | \vecv^\top \vecG_t( \bfeps_t ) - \vecv^\top \vecG_t( \bfeps_{t,t-j} ) |^q \right)^{\frac{1}{q}} \\
		& \quad \le 2 K \| \vecv \|_1  C_1 j^{-\beta} + 2 K \| \vecv \|_1^2 \max_{1 \le k \le d} 
		\left( \E | G_{tk}(\bfeps_t) - G_{tk}(\bfeps_{t,t-j}) |^q \right)^{\frac{1}{q}}  \\
		& \quad \le 2K \| \vecv \|_1 (\| \vecv \|_1 + 1 ) C_1 j^{-\beta}.
	\end{align*}
    This also  yields for the associated vector version 
    \begin{align*}
      & 	\left( \E \left\| \vecG_t(\bfeps_t) \vecv^\top \vecG_t(\bfeps_t ) - \vecG_t(\bfeps_{t,t-j} ) \vecv^\top \vecG_t(\bfeps_{t,t-j} ) \right\|_2^{\frac{q}{2}} \right)^{\frac{2}{q}} \\
      & \qquad =  \left( \E \left( \sum_{k=1}^d (G_{tk}(\bfeps_t) \vecv^\top \vecG_t(\bfeps_t ) - G_{tk}(\bfeps_{t,t-j} ) \vecv^\top \vecG_t(\bfeps_{t,t-j} ) )^2 \right)^{\frac{q}{4}}  \right)^{\frac{2}{q}} \\
      & \qquad \le \left(  d^{\frac{q}{4}-1} \sum_{k=1}^d \E (G_{tk}(\bfeps_t) \vecv^\top \vecG_t(\bfeps_t ) - G_{tk}(\bfeps_{t,t-j} ) \vecv^\top \vecG_t(\bfeps_{t,t-j} )  )^{\frac{q}{2}}  \right)^{\frac{2}{q}} \\
      & \qquad \le  2K \| \vecv \|_1 (\| \vecv \|_1 + 1 ) C_1 \sqrt{d} j^{-\beta},
    \end{align*}
verifying (\ref{PhysDepMeas}).  Assumption A--ii can be shown as follows. We have
\begin{align*}
  \sum_{t=1}^n \| \vecv^\top( \vecG_{t}( \bfeps_0 ) - \vecG_{t-1}( \bfeps_0 ) ) \|_{L_4}
  &\le \sum_{j=1}^d |v_j | \sum_{t=1}^n \| G_{tj}( \bfeps_0 ) -  G_{t-1,j}( \bfeps_0 ) \|_{L_4}\\
  & \le \| \vecv \|_1 KL.
\end{align*}
This gives, using $\| \vecv^\top \vecY_t \|_{L_4} \le 2 K \| \vecv \|_1$,
\begin{align*}
	& \| G_{tk}( \bfeps_0 ) \vecv^\top \vecG_t( \bfeps_0 ) - G_{t-1,k}( \bfeps_0 ) \vecv^\top \vecG_{t-1}( \bfeps_0 ) \|_{L_2} \\
	& \le \| \vecv^\top \vecG_t( \bfeps_0 ) \|_{L_4} \| G_{tk}( \bfeps_0 ) - G_{t-1,k}( \bfeps_0 ) \|_{L_4} +
	\| G_{t-1,k}( \bfeps_0 ) \|_{L_4} \| \vecv^\top( \vecG_t( \bfeps_0 ) - \vecG_{t-1}( \bfeps_0 ) ) \|_{L_4} \\
	& \le 2 K \| \vecv \|_1 \| G_{tk}( \bfeps_0 ) - G_{t-1,k}( \bfeps_0 ) \|_{L_4} + 2 K \| \vecv^\top( \vecG_t( \bfeps_0 ) - \vecG_{t-1}( \bfeps_0 ) ) \|_{L_4},
\end{align*}
and therefore
\[
  \sum_{t=1}^n \| G_{tk}( \bfeps_0 ) \vecv^\top \vecG_t( \bfeps_0 ) - G_{t-1,k}( \bfeps_0 ) \vecv^\top \vecG_{t-1}( \bfeps_0 ) \|_{L_2} \le 4 \| \vecv \|_1 K^2 L.
\]
This verifies A--ii for $ V(2) $ with constant $ 4 \| \vecv \|_1 K^2 L$ instead of $ KL $.  

To show (i) use the fact that 
\[
  \vecY_t \vecY_t^\top - \matM_t = (\vecY_t - \bfmu_t)(\vecY_t - \bfmu_t)^\top + (\vecY_t - \bfmu_t) \bfmu_t^\top  + \bfmu_t(\vecY_t - \bfmu_t)^\top - \bfSigma_t
\]
yielding
\begin{align*}
  \| \vecv^\top(\vecY_t\vecY_t^\top - \matM_t) \vecv \|_{L_{q/2}} 
  &\le \| \vecv^\top(\vecY_t-\bfmu_t) \|_{L_q}^2 + 2 |\bfmu_t^\top \vecv | \| \vecv^\top (\vecY_t - \bfmu_t) \|_{L_q} + \max_{1 \le j, k \le d} | \bfSigma_{t,jk} | \| \vecv \|_1^2 \\
  & \le \| \vecv \|_1^2 \max_{1 \le j \le d} \| Y_{tj} - \mu_{tj} \|_{L_q}^2 
  + 2 K \| \vecv \|_1^2  \max_{1 \le j \le d} \| Y_{tj} - \mu_{tj} \|_{L_q}  + K^2 \| \vecv \|_1^2 \\
  & \le 4 K^2 \| \vecv \|_1^2.
\end{align*}
Further, observing that $   \vecv^\top \vecY_t \vecY_t^\top \vecv - \vecv^\top \vecG_t( \bfeps_{t,t-j} )  \vecG_t(\bfeps_{t,t-j} )^\top \vecv  $ can be written as the sum of
$
  \vecv^\top( \vecG_t( \bfeps_t ) - \vecG_t(\bfeps_{t,t-j} ) ) \vecG_t( \bfeps_t )^\top\vecv  $ and 
  $ \vecv^\top\vecG_t( \bfeps_{t,t-j} ) ( \vecG_t( \bfeps_t ) - \vecG_t( \bfeps_{t,t-j} ) )^\top \vecv $, we obtain 
\begin{align*}
	   \| \vecv^\top \vecY_t \vecY_t^\top \vecv - \vecv^\top \vecG_t( \bfeps_{t,t-j} )  \vecG_t(\bfeps_{t,t-j} )^\top \vecv \|_{L_{q/2}} 
	  &  \le 2 \| \vecv^\top( \vecG_t(\bfeps_t) - \vecG_t(\bfeps_{t,t-j}) ) \vecv^\top \vecG_t(\bfeps_t ) \|_{L_{q/2}} \\
	& \le 2 \| \vecv^\top \vecY_t \|_{L_q} \| \vecv^\top( \vecG_t(\bfeps_t) - \vecG_t(\bfeps_{t,t-j}) ) \|_{L_q} \\
	& \le 4 K \|\vecv\|_1^2 \max_{1 \le k \le d} \| G_{tk}( \bfeps_t ) - G_{tk}( \bfeps_{t,t-j} ) \|_{L_q} \\
	& \le 4 K \|\vecv \|_1^2 C_1 j^{-\beta},
\end{align*}
which verifies A--i with constant $ 4 K \| \vecv \|_1^2 C_1 $.
Similarly, one shows that A--ii holds for $ V(2) $ with  $KL$ replaced by $ 4 \|\vecv\|_1^2 K^2L$: By symmetry,
\begin{align*}
  \vecv^\top \vecY_t \vecY_t^\top \vecv - \vecv^\top \vecY_{t-1} \vecY_{t-1}^\top \vecv 
  &= \vecv^\top( \vecY_t - \vecY_{t-1}) \vecY_t^\top \vecv + \vecv^\top \vecY_{t-1} (\vecY_t - \vecY_{t-1})^\top \vecv \\
  &= 2 \vecv^\top( \vecY_t - \vecY_{t-1}) \vecY_t^\top \vecv,
\end{align*}
such that by A--ii
\begin{align*}
	\sum_{t=1}^n \| \vecv^\top \vecY_t \vecY_t^\top \vecv - \vecv^\top \vecY_{t-1} \vecY_{t-1}^\top \vecv  \|_{L_{2}} & \le 2 \| \vecv\|_1 \sum_{t=1}^n \max_{1 \le j \le d} \| Y_{tj} \|_{L_4} \| \vecv^\top( \vecY_t - \vecY_{t-1} ) \|_{L_4} \\
	& \le 4  \| \vecv \|_1 K \sum_{t=1}^n \| Y_{tj} - Y_{t-1,j} \|_{L_4} \\
	& \le 4 \| \vecv \|_1^2 K^2 L.
\end{align*}

  Assertions (iii) and (iv) follow by Jensen's inequality, 
\[ 
E \| \vecY_t \|_2^q \le d^\frac{q}{2} \E\left( \frac{1}{d} \sum_{j=1}^d Y_{tj}^2 \right)^{\frac{q}{2}}  \le d^\frac{q}{2} \max_{1 \le j \le d} \E(Y_{tj}^q), 
\]
and 
\begin{align*}
	\E \| \vecG_t( \bfeps_t ) - \vecG_t( \bfeps_{t,t-j} ) \|_2^q 
	&= \E \left( \sum_{k=1}^d \E (G_{tk}( \bfeps_t ) - G_{tk}(\bfeps_{t,t-j} ) )^2 \right)^{\frac{q}{2}}  \\
	& \le d^\frac{q}{2}  d^{-1} \sum_{k=1}^d \E \left| G_{tk}(\bfeps_t) - G_{tk}( \bfeps_{t,t-j} ) \right|^{q}  \\
	& \le  d^\frac{q}{2} \max_{1 \le k \le d} \E \left| G_{tk}(\bfeps_t) - G_{tk}( \bfeps_{t,t-j} ) \right|^{q} 
\end{align*}
	
	Let us now verify (v):  Write $ \vecH_t( \bfeps_t ) = ( H_{tk}( \bfeps_t) )_{k=1}^{\bar{d}} $ for the $\bar{d} = d^2 $--dimensional vectorized series. By Jensen's inequality we may estimate the $j$th physical dependence measure with respect to the vector--2 norm in $ L_q $ as follows: 
	\begin{align*}
		  \left(  \E \| \vecH_t( \bfeps_t ) - \vecH_t( \bfeps_{t,t-j} ) \|_2^q \right)^{\frac{1}{q}} 
		  & = \left(  \bar{d}^{\frac{q}{2}} \E \left( \frac{1}{\bar{d}}  \sum_{k=1}^{\bar{d}} [H_{tk}( \bfeps_t ) - H_{tk}(\bfeps_{t,t-j})]^2 \right)^{\frac{q}{2}} \right)^{\frac{1}{q}} \\
		  & \le\left( \bar{d}^{\frac{q}{2}-1} \sum_{k=1}^{\bar{d}} \E | H_{tk}(\bfeps_t) - H_{tk}(\bfeps_{t,t-j} ) |^q \right)^{\frac{1}{q}} \\
		  & \le \bar{d}^{\frac{1}{2} - \frac{1}{q}} \left( \bar{d} C_1^q j^{-\beta q} \right)^{\frac{1}{q}} \\
		  & = C_1 d j^{-\beta}.
	\end{align*}
    The moment bound $ \left( \E \| \vecH_t( \bfeps_t ) \|_2^q \right)^{\frac{1}{q}} \le K d $ follows analogously.
    
\end{myproof}

\begin{myproof}{Theorem~\ref{LemmaRatesMeanCovEst}} 
The product series $ Y_{tj}  Y_{tk} = G_{tj}( \bfeps_t) G_{tk}(\bfeps_t) $, $ t \ge 1 $, satisfies Assumption A--i with power $ q/2 $ instead of $q$: Firstly, $ \E | Y_{ij} Y_{ik} |^\frac{q}{2} \le \| Y_{ij} \|_{L_2} \| Y_{ik} \|_{L_q} \le 4 K^2 $, see the proof of Lemma~\ref{WeakDepLeftProj}. Secondly,  $ \E( G_{tj}( \bfeps_t ) G_{tk}( \bfeps_t) ) = \E( G_{tj}( \bfeps_{t,t-\ell} ) G_{tk}( \bfeps_{t,t-\ell} )  ) $, for $ \ell \ge 1 $, such that 
	\begin{align*}
		\| G_{tj}( \bfeps_t ) G_{tk}( \bfeps_t) - G_{tj}( \bfeps_{t,t-\ell} ) G_{tk}( \bfeps_{t,t-\ell} ) \|_{L_{\frac{q}{2}}}   
		&\le \| G_{tk}( \bfeps_0 ) \|_{L_q} \| G_{tj}( \bfeps_t ) - G_{tj}( \bfeps_{t,t-\ell} ) \|_{L_q} \\
		& \quad + \| G_{tj}( \bfeps_0 ) \|_{L_q} \| G_{tk}( \bfeps_t ) - G_{tk}( \bfeps_{t,t-\ell} ) \|_{L_q} \\
		& \le 4 K^2 C_1  \ell^{-\beta},
	\end{align*}
	Thus, condition (G.1) of \cite{MiesSteland2022} holds with $ \Theta = \max(4 K^2, 4 K C_1) $, such that Theorem~3.1 therein yields the existence of a universal constant $c_1 $ such that 	
	\[
	\max_{1 \le j, k \le d} \left( \E \left| \sum_{t=1}^m (Y_{tj} Y_{tk} - \E(Y_{tj} Y_{tk})  ) \right|^{\frac{q}{2}} \right)^{\frac{2}{q} }\le c_1 \sqrt{m} \zeta(\beta),
	\]
	i.e.,  by monotonicity of the norms,
	$ 
	\max_{1 \le j,k \le d} \left( \E \left| \hat{M}_{m,jk} - \bar M_{m,jk} \right|^{r/2} \right)^{2/r} \le \frac{c_1}{\sqrt{m}} \zeta(\beta)
	$, for $ 2 \le r \le q $.
	Similarly, one shows that $ \max_{1 \le j \le d} \left( \E | \hat{\mu}_{mj} - \bar \mu_{mj} |^{r} \right)^{1/r} \le \frac{c_2}{\sqrt{m}} \zeta(\beta)$ for some constant $ c_2 $, $ 2 \le r \le q $.	This also implies
	\begin{align*}
		\| \hat{\mu}_{mj} \hat{\mu}_{m,k} - \bar \mu_{mj} \bar \mu_{mk}  \|_{L_\frac{r}{2}}
		& \le \| \hat{\mu}_{mk} \|_{L_r} \| \hat{\mu}_{mj} - \bar \mu_{mj}\|_{L_{r}} 
		+ | \mu_j | \| \hat{\mu}_{mk} - \bar \mu_{mk} \|_{L_r} \\
		& \le \left( \frac{c_2}{\sqrt{m}} \zeta(\beta) + K \right) \frac{c_2}{\sqrt{m}} \zeta(\beta)+ \frac{K c_2}{\sqrt{m}} \zeta(\beta)\\
		& \le \frac{c_3}{\sqrt{m}} \zeta(\beta),
	\end{align*}
	with $ c_3 = (c_2 \zeta(\beta) + K ) c_2 + K c_2 $. 
	Noting that $ \bar\bfSigma_m = \bar\matM_m + \bar\bfmu_m \bar\bfmu_m^\top $ and $ \hat{\bfSigma}_m = \hat{\matM}_m + \hat{\bfmu}_m \hat{\bfmu}_m^\top $ we get from the triangle inequality $ \max_{1 \le j, k \le d} \left( \E | \hat{\Sigma}_{m,jk} - \bar{\Sigma}_{m,jk} |^{r/2} \right)^{2/r}  \le \frac{c_4}{\sqrt{m}} \zeta(\beta) $ for some constant $ c_4 $. Now the theorem follows with $c = \max(c_1, \ldots, c_4) $.
\end{myproof}

\begin{Lemma} 
\label{LemmaBoundEstimatedV}
	Suppose that Assumptions A--i and A--ii as well as A--iv hold, so that particularly $ \| \vecv \|_1 \le C_{\vecv} $ and \[ \| \hat{\vecv}_m - \vecv \| = O_\P\left( \sqrt{\frac{r_d }{m} } \right), \]
	with $ \sqrt{s r_d} = O( m^{\zeta} ) $ for some $ 0 \le \zeta < \frac{1}{2} $ under $s$--sparsity, whereas $  \sqrt{d r_d} = O( m^{\zeta} )  $ for non--$\ell_0$--sparse estimation.
	Then
	\[
	 R_m = \max_{k \le m} \left| \sum_{i \le k} \frac{\hat{\vecv}_m^\top (\vecY_i \vecY_i^\top - \matM_i ) \hat{\vecv}_m }{\sigma_{0,\infty}} - \sum_{i \le k} \frac{\vecv^\top ( \vecY_i \vecY_i^\top - \matM_i)\vecv }{\sigma_{0,\infty}} \right| = O_\P( \sqrt{d r_d} ),
	\]
	and under $s$--sparseness
	\[
R_m = \max_{k \le m} \left| \sum_{i \le k} \frac{\hat{\vecv}_m^\top (\vecY_i \vecY_i^\top - \matM_i ) \hat{\vecv}_m }{\sigma_{0,\infty}} - \sum_{i \le k} \frac{\vecv^\top ( \vecY_i \vecY_i^\top - \matM_i)\vecv }{\sigma_{0,\infty}} \right| = O_\P( \sqrt{s r_d} ).
\]
\end{Lemma}

\begin{myproof}{of Lemma~\ref{LemmaBoundEstimatedV}}  To simplify presentation assume $ \sigma_{0,\infty}^2 =1 $. 
	We have the decomposition, by symmetry of $ \vecY_i \vecY_i^\top - \matM_i $,
	\begin{align*}
	  R_m &\le 
	  \max_{k \le m} \left| 
\sum_{i \le k} (\hat{\vecv}_m-\vecv)^\top ( \vecY_i \vecY_i^\top - \matM_i ) (\hat{\vecv}_m-\vecv) \right| 
   +  2  \max_{k \le m} \left| 
\sum_{i \le k} \vecv^\top ( \vecY_i \vecY_i^\top - \matM_i ) (\hat{\vecv}_m-\vecv)
\right|.
	\end{align*}
We shall estimate the terms on the right hand side separately. By Lemma~\ref{WeakDepLeftProj}~(ii) the series $ \vecY_i \vecY_i^\top \vecv $ satisfies Assumptions A--i and A--ii. Consequently, since $ q > 4 $, $ d = O(\sqrt{m} ) $ and the physical dependence measure of $ \vecY_i \vecY_i^\top \vecv $  is  $O(\sqrt{d}) $, see (\ref{PhysDepMeas}), \cite[Th.~3.2]{MiesSteland2022} with $ r = 2 $ yields
\[
  \E \max_{k \le m} \left\|  \sum_{i \le k} \vecv^\top (\vecY_i \vecY_i ^\top - \matM_i) \right\|_2^2 \le c_1 d m 
\]
for some constant $c_1$. By Assumption A--iv, for any given $ \varepsilon > 0 $ we can find a constant $ C > 0 $ such that 
\begin{equation}
	\label{VEstOP1}
	\P\left( \| \hat{\vecv}_m - \vecv \|_2 > C \sqrt{\frac{r_d}{m} } \right)  < \frac{\varepsilon}{2}.
\end{equation}
Using the elementary inequality $ | \vecx^\top \matA \vecx | \le \| \matA \vecx \|_2 \| \vecx \|_2 $, $ \matA \in \R^{d \times d} $, $ \vecx \in \R^d $, we obtain for the second term of the above decomposition, for any $ B > 0 $,
\begin{align*}
	& \P\left(  \frac{2}{\sqrt{d r_d}} \max_{k \le m}  \left| \sum_{i \le k}\vecv^\top (\vecY_i \vecY_i ^\top - \matM_i) ( \hat{\vecv}_m - \vecv) \right| > B \right)  \\
	& \quad \le \P\left( 2C \sqrt{\frac{1}{dm} } \max_{k \le m} \left\|  \sum_{i \le k} \vecv^\top (\vecY_i \vecY_i ^\top - \matM_i)  \right\|_2 > B \right) + 
	\P\left( \| \hat{\vecv}_m - \vecv \|_2 > C \sqrt{\frac{r_d}{m}} \right) \\
	& \quad \le \frac{ 4C^2 \E \max_{k \le m} \left\|  \sum_{i \le k} \vecv^\top (\vecY_i \vecY_i ^\top - \matM_i) \right\|_2^2 }{B^2 dm}
	+ \frac{\varepsilon}{2}
\end{align*}
Now choose $B $ large enough to ensure that the last bound is less than $ \varepsilon $. Therefore, 
\[
  \max_{k \le m} \left\| 2  \vecv^\top\sum_{i \le k}  ( \vecY_i \vecY_i^\top - \matM_i ) ( \hat{\vecv}_m - \vecv)  \right\|_2 = O_\P( \sqrt{d r_d} ).
\]
Further, by Lemma~\ref{WeakDepLeftProj}~(ii) and \cite[Th.~3.2]{MiesSteland2022}, for some constant $ c_ 2$,
\[
\E \max_{k \le m} \left\| \sum_{i \le k} \vecY_i \vecY_i^\top - \matM_i  \right\|_2^2  \le c_2 d^2 m,
\] 
since $ d = O( \sqrt{m} ) $ by assumption. Hence,
\begin{align*}
	& \P\left(  \frac{\sqrt{m}}{d r_d} \max_{k \le m}  \left| \sum_{i \le k} ( \hat{\vecv}_m - \vecv)^\top (\vecY_i \vecY_i ^\top - \matM_i) ( \hat{\vecv}_m - \vecv) \right| > B \right)  \\
	& \quad \le \P\left( \frac{C^2}{d \sqrt{m}} \max_{k \le m} \left\|  \sum_{i \le k} \vecY_i \vecY_i ^\top - \matM_i \right\|_2 > B \right) + 
	\frac{\varepsilon}{2} \\
	& \quad \le \frac{ 4C^2 \E \max_{k \le m} \left\|  \sum_{i \le k} \vecY_i \vecY_i ^\top - \matM_i \right\|_2^2 }{Bd^2 m}
	+ \frac{\varepsilon}{2}
\end{align*}
By increasing $B$, if necessary, this upper bound is less than $ \varepsilon $. Thus,
\[
\max_{k \le m} \left| 
\sum_{i \le k} (\hat{\vecv}_m-\vecv)^\top ( \vecY_i \vecY_i^\top - \matM_i ) (\hat{\vecv}_m-\vecv) \right| = O_\P\left( \frac{d r_d}{\sqrt{m}} \right)
= O_\P( \sqrt{d r_d} ),
\]
since $ \sqrt{d r_d} = o( \sqrt{m} ) $ by assumption. Let us now consider the $\ell_0$--sparse setting. Denote $ A = \{ j : v_j \not= 0 \}$, $ \hat{A} = \{ j: \hat{v}_j \not= 0 \} $ and let $ \vecY_{iA} = ( Y_{ij} )_{j \in A}, \vecv_A = (v_j)_{j \in A} $ as well as $ \hat{\vecv}_A = ( \hat{v}_j )_{j \in A} $. Notice that on the event $ \{ \hat{A} = A \} $ quadratic forms in the remainder $R_m $ as well as in the terms of the upper bound we need to estimate now collapse to quadratic forms in terms of the $s$--dimensional (random) vectors $\hat{\vecv}_A $, $ \vecv_A $ and $ \vecY_{iA} $. The arguments given for the non--$\ell_0$--sparse case apply with $ d = s $. For example, now we have
\[
\E \max_{k \le m} \left\| \sum_{i \le k} \vecY_{iA} \vecY_{iA}^\top - \matM_{iA}  \right\|_2^2  \le c_2 s^2 m,
\] 
since $ s = O(\sqrt{m} ) $ by assumption.  On the event $ \{ \hat{A} = A  \} $ it holds
$ ( \hat{\vecv}_m - \vecv)^\top (\vecY_i \vecY_i ^\top - \matM_i) ( \hat{\vecv}_m - \vecv)
= ( \hat{\vecv}_A - \vecv_A)^\top (\vecY_{iA} \vecY_{iA} ^\top - \matM_i) ( \hat{\vecv}_A - \vecv_A) $, and therefore 
\begin{align*}
	&\P\left(  \frac{\sqrt{m}}{s r_d} \max_{k \le m}  \left| \sum_{i \le k} ( \hat{\vecv}_m - \vecv)^\top (\vecY_i \vecY_i ^\top - \matM_i) ( \hat{\vecv}_m - \vecv) \right| > B \right) \\
	& \le \P\left(  \frac{\sqrt{m}}{s r_d} \max_{k \le m}  \left| \sum_{i \le k} ( \hat{\vecv}_A - \vecv_A)^\top (\vecY_{iA} \vecY_{iA} ^\top - \matM_{iA}) ( \hat{\vecv}_A - \vecv_A) \right| > B \right)  + \P( A \not= \hat{A})
\end{align*}
yielding
\[
\max_{k \le m} \left| 
\sum_{i \le k} (\hat{\vecv}_m-\vecv)^\top ( \vecY_i \vecY_i^\top - \matM_i ) (\hat{\vecv}_m-\vecv) \right| = O_\P\left( \frac{s r_d}{\sqrt{m}} \right)
= O_\P( \sqrt{s r_d} ),
\]
if $ \sqrt{s r_d} = o( \sqrt{m} ) $. The first term is treated analogously. This completes the proof.
\end{myproof}

\begin{Lemma}
\label{LemmaFiniteLRV}
	Under Assumptions A--i and A--ii it holds for all $ 1 \le j, k \le d $, $ t, s \ge 1 $ and all $d, m $
	\[
	  \Cov( G_{tj}(\bfeps_t), G_{sk}( \bfeps_{s+h} ) ) \le c \zeta( \beta )
	\]
	for some constant $ c $. Especially, if $ \beta > 2 $ then $ \sum_{h \in Z} | \Cov( G_{tj}(\bfeps_t), G_{sk}( \bfeps_{s+h} ) )  | < \infty $.
\end{Lemma}

\begin{myproof}{of Lemma~\ref{LemmaFiniteLRV}} The proof goes along the lines of the proof of \cite[Prop.~5.4]{MiesSteland2022}. For completeness, we provide the details. Let $ 
	\bar {\bfeps}_{h,0} = (\epsilon_h, \ldots, \epsilon_{1},\epsilon_0', \epsilon_1', \ldots ) $. 
  By independence of $ G_{tj}(\bfeps_0) $ and $ G_{sk}(\bar{\bfeps}_{h,0}) $,
  \begin{align*}
  	  \Cov( G_{tj}(\bfeps_t), G_{sk}( \bfeps_{s+h} ) ) 
  	  & = \Cov( G_{tj}(\bfeps_0), G_{sk}( \bfeps_{h} ) )  \\
  	  & = \Cov( G_{tj}( \bfeps_0 ), G_{sk}(\bfeps_h) - G_{sk}( \bar\bfeps_{h,0} ) )
  	  + \Cov( G_{tj}( \bfeps_0 ), G_{sk}( \bar\bfeps_{h,0}  )) \\
  	  & = \Cov( G_{tj}(\bfeps_0), G_{sk}( \bfeps_{h} ) - G_{sk}(\bar\bfeps_{h,0} ) ).
  \end{align*}
  Using $ \| G_{tj}( \bfeps_0 ) - \E G_{tj}( \bfeps_0 ) \|_{L_2} \le K $ and $ \|G_{sk}( \bfeps_{h} ) - G_{sk}(\bar\bfeps_{0,h} ) \|_{L_2} = O( \sum_{\ell=h}^\infty \ell^{-\beta} )$ (by a telescoping argument and applying Assumption A-ii to each term), the assertion follows from the Cauchy-Schwarz inequality.  
\end{myproof}

\begin{Lemma} 
\label{LemmaDiffLRV}
	Under Assumptions A--i - A--iii it holds 
	\[
	  | \sigma_{i,\infty}^2 - \sigma_{0,\infty}^2 | \le 
	  \frac{C K^2 \| \vecv \|_1^2}{m} \zeta(\beta-1),
	\]
	uniformly in $ i \ge 1 $, for some constant $ C$.
\end{Lemma}

\begin{myproof}{of Lemma~\ref{LemmaDiffLRV}} Observe that by Assumption A--iii there exists some constant $C_2 $ such that for all $ i \ge 1 $ and $ 1 \le j, k \le d $ 
	\[
	| \Cov( G_{ij}(\bfeps_0), G_{ik}(\bfeps_h)) - \Cov( G_{0j}(\bfeps_0), G_{0k}( \bfeps_h )) | \le \frac{C_2}{m} | \Cov( G_{0j}(\bfeps_0),  G_{0k}(\bfeps_{h} ) ) |.
	\]
	Therefore, using Lemma~\ref{LemmaFiniteLRV}
	\begin{align*}
		| \sigma_{i,\infty}^2 - \sigma_{0,\infty}^2 | 
		& \le \sum_{h \in \Z} |  \Cov( \vecv^\top\vecG_i( \bfeps_0 ), \vecv^\top\vecG_i(\bfeps_h ) ) - \Cov( \vecv^\top\vecG_0( \bfeps_0 ), \vecv^\top\vecG_0( \bfeps_h ) )  | \\
		& \le  \sum_{h \in \Z} \sum_{j,k=1}^d |v_j | | v_k || \Cov( G_{ij}(\bfeps_0), G_{ik}(\bfeps_h)) - \Cov( G_{0j}(\bfeps_0), G_{0k}( \bfeps_h )) | \\
		& \le \frac{C_2 \| \vecv \|_1^2}{m} \sum_{h \in \Z} \max_{1\le j,k \le d} | \Cov( G_{0j}(\bfeps_0),  G_{0k}(\bfeps_{h} ) ) |  \\
		& \le \frac{C_2  K^2  \| \vecv \|_1^2}{m} \sum_{h \in \Z}  \sum_{\ell=h}^\infty \ell^{-\beta}   \\
		& \le \frac{C K^2\| \vecv \|_1^2}{m} \sum_{\ell=0}^\infty \ell^{-\beta+1}
	\end{align*}
uniformly in $ i \ge 1 $, for some constant $ C$.
\end{myproof}

\begin{myproof}{of Theorem~\ref{ThIP}}
	 (i): Observe that by Assumption A--iii \[ \inf_{i \ge 1} \sigma_{i,\infty}^2 \ge \inf_{i \ge 1}  \vecv^\top \bfSigma_{i,\infty} \vecv > 0. \] Define the random variables
	\[
	\xi_i = \vecv^\top (\vecY_i  \vecY_i^\top - \matM_i)\vecv  / \sigma_{i,\infty} \quad \text{with} \quad
	\sigma_{i,\infty}^2 = \sum_{h \in \Z} \vecv^\top \Cov( \vecG_i( \bfeps_0 ),  \vecG_i( \bfeps_h ) ) \vecv, \ i \ge 1,
	\]
	and let $ S_n = \sum_{i=1}^n \xi_i $. 
		Denote the $j$th physical dependence coefficients of $ \{ \xi_{i} \} $  by $ \vartheta_{j,q} $, $ j \ge 1 $, and the tail dependence measure by $ \Theta_{i,q} = \sum_{j \ge i} \vartheta_{j,q} $. By Lemma~\ref{WeakDepLeftProj}~(i) the random variables $ \xi_i $ have a summable tail dependence measure, since 
	$\Theta_{j,q} \le 4 K \| \vecv \|_1^2 C_1  j^{-\beta+1} $, and satisfy the assumptions of \cite[Th.~3.1]{MiesSteland2022}. Hence, after redefining the vector time series on a new probability space, there exist i.i.d. random variables $ G_i \sim N(0,1) $, such that for fixed $ 0  < \varepsilon < \xi(q,\beta)  $
	\begin{equation}
		\label{GaussianApprox}
		\max_{k\le m} \left| S_k - \sum_{i=1}^k G_i \right| = O_\P( \sqrt{\log(m)} m^{\frac{1}{2} - \xi(q,\beta)} ) = o_\P( m^{\frac{1}{2} - \xi(q,\beta)  + \varepsilon } )
		= o_\P( m^{\frac{1}{\nu_1}} ), 
	\end{equation}
	where $ \nu_1 = \frac{1}{\frac{1}{2} - \xi(q,\beta) + \varepsilon} > 2 $. Next we show that we can use $ \sigma_{0,\infty}^2 = \vecv^\top \bfSigma_{0,\infty} \vecv $ instead of $ \sigma_{i,\infty}^2 = \vecv^\top \bfSigma_{i,\infty} \vecv $ for standardization: Using $ | \sqrt{x} - \sqrt{y} | \le \sqrt{|x-y|} $, we obtain by Assumption A--iii, Lemma~\ref{WeakDepLeftProj}~(i), Lemma~\ref{LemmaDiffLRV}  and \cite[Th.~3.2]{MiesSteland2022} 
		\begin{align*}
		& \E \max_{k \le m} \left| \sum_{i=1}^k \frac{\vecv^\top ( \vecY_i \vecY_i^\top - \matM_i)\vecv}{ \sigma_{0,\infty} } - \sum_{i=1}^k \frac{\vecv^\top (\vecY_i \vecY_i^\top - \matM_i)\vecv}{ \sigma_{i,\infty} } \right|^2 \\
		&\qquad \le \E \max_{k \le m} 
		\sum_{i=1}^k \left| \frac{\sqrt{|\sigma_{0,\infty}^2 - \sigma_{i,\infty}^2|} } {\sigma_{0,\infty} \sigma_{i,\infty}} \right|^2 | \vecv^\top  (\vecY_i \vecY_i^\top - \matM_i )\vecv |^2 \\
		& \qquad \le \sup_{i \ge 1} \frac{{|\sigma_{0,\infty}^2 - \sigma_{i,\infty}^2| }} {\sigma_{0,\infty}^2 \inf_{\ell \ge 1} \sigma_{\ell,\infty}^2}
		\E \max_{k \le m} \sum_{i=1}^k | \vecv^\top  (\vecY_i \vecY_i^\top - \matM_i )\vecv |^2\\
		& \qquad 
		=  O\left( 1\right),
	\end{align*}
	Therefore,
	\[
	\max_{k \le m} \left| \sum_{i=1}^k \frac{\vecv^\top (\vecY_i \vecY_i^\top - \matM_i )\vecv}{ \sigma_{0,\infty} } - S_k \right| = O_\P( 1 ) = o_\P( m^{\frac{1}{\nu}} ).
	\]
	We can assume that $ G_i = B_{i} - B_{i-1} $, $ i \in \N $, for some standard Brownian motion $ B_{t} $, $ t \ge 0 $, such that
	\[
	\max_{k \le m} | S_k - B_k | = o_\P( m^{\frac{1}{\nu_1}} ),
	\]
	which verifies (i) with $ \nu = \nu_1$.
	
	(ii): By Lemma~\ref{LemmaBoundEstimatedV}~(iii), under Assumption A--iv, which provides $ d = O(\sqrt{m} ) $ and $ \sqrt{d r_d} = O( m^\zeta ) $ for some $ 0 \le \zeta < \frac{1}{2} $,  one can replace $ \vecv $ by its estimator $ \hat{\vecv}_m $, because 
	\begin{align}
		\label{OrderWithEst}
		\max_{k \le m}  \left| \sum_{i \le k} \frac{ \hat{\vecv}_m^\top (\vecY_i\vecY_i^\top - \matM_i ) \hat{\vecv}_m}{\sigma_{0,\infty}} - \sum_{i \le k} \frac{\vecv^\top ( \vecY_i \vecY_i^\top - \matM_i ) \vecv}{\sigma_{0,\infty}} \right|
		& = O_\P\left(  \sqrt{d r_d} \right) 
		= o_\P\left( m^{\frac{1}{\nu_2}} \right), 
	\end{align}
	where $ \nu_2= \frac{1}{\zeta+\varepsilon} > 2 $, if $ \varepsilon $ is decreased, when necessary, to ensure $ \zeta + \varepsilon < \frac{1}{2} $. 
		Hence, by the triangle inequality, \eqref{GaussianApprox} remains true when replacing the $ \xi_i $'s by
	\[ 
	\hat{\xi}_i =  \frac{\hat{\vecv}_m^\top (\vecY_i \vecY_i^\top - \matM_i ) \hat{\vecv}_m}{\sigma_{0,\infty}}, \qquad 1 \le i \le m,
	\]
	and $ \nu_1 $ by $ \nu = \min( \nu_1, \nu_2 ) > 2 $. Thus, if $ \hat{S}_k = \sum_{i \le k} \hat{\xi}_i $, (\ref{OrderWithEst}) shows that 
	$
	\max_{k \le m} | \hat{S}_k - S_k | = o_\P\left( m^\frac{1}{\nu} \right),
	$    
	 and we obtain
	\[
	\max_{k\le m} \left| \hat{S}_k - B_k \right| = o_\P( m^{\frac{1}{\nu}} ).
	\]
	It remains to consider the case of $ s $--sparseness. Under Assumption A--iv the quadratic forms appearing in (\ref{OrderWithEst}) collapse to quadratic forms in terms of $ \hat{\vecv}_A, \vecv_A $ and $ \vecY_{iA} $. Thus, all properties derived from the vector time series can assume that the dimension is $s$. Note that the assumption of \cite[Th.~3.1]{MiesSteland2022} that the dimension is $ O(m) $ is void, since we apply it to a univariate series, and the dimension only appears in the constants.  Lemma~\ref{LemmaBoundEstimatedV}~(iii) provides the bound $ O_\P( \sqrt{s r_d} ) $ which is of the order $ o_\P( m^\frac{1}{\nu_2} ) $, 
    if $ \sqrt{s r_d} = O( m^\zeta ) $ as assumed in Assumption A--iv. Now the proof can be completed as for non--$\ell_0 $--sparse estimation.
\end{myproof}

\begin{myproof}{of Theorem~\ref{ThH0}}
	Under the null hypothesis   $ \vecv^\top \matM_i \vecv = \theta_ 0 $, $ i \ge 1 $, we have
	\[
	  \E_0( \vecv^\top \matS_i \vecv ) = \vecv^\top \E_0( \matS_i ) \vecv = \theta_0 
	\]
	where $ \E_0 $ indicates that the expectation is taken under $ H_0 $. Therefore, 
	centering all summands under $ \E_0 $ yields
	\[
	\sum_{m < i \le m+k} {\vecv}^\top \matS_i  {\vecv} - \frac{k}{m} \sum_{j=1}^m {\vecv}^\top \matS_j  {\vecv}
	= 
	\sum_{m < i \le m+k} {\vecv}^\top (\matS_i - \matM_i) {\vecv} - \frac{k}{m} \sum_{j=1}^m {\vecv}^\top (\matS_j - \matM_j) {\vecv}.
	\]
	Therefore,  for known $ \vecv $ we can be write
	\[
	\max_{\delta m \le k} \frac{Q(m,k)}{\hat{\sigma}_{0,\infty} g(m,k) } = \max_{\delta m \le k < \infty} \frac{\left| \sum_{m < i \le m+k} {\vecv}^\top (\matS_i - \matM_i) {\vecv} - \frac{k}{m} \sum_{j=1}^m {\vecv}^\top (\matS_j - \matM_j) {\vecv}  \right|}{\hat{\sigma}_{0,\infty} g(m,k)}.
	\]
	As in the proof of Theorem~\ref{ThIP}, denote $ \xi_i ={\vecv}^\top (\matS_i - \matM_i) {\vecv} / \sigma_{0,\infty}  $, 
	 $ \hat{\xi}_i = \hat{\vecv}_m^\top (\matS_i - \matM_i) \hat{\vecv}_m / \sigma_{0,\infty} $, $ i \ge 1 $, and  $ S_k = \sum_{i=1}^k \xi_i $,  $ \hat{S}_k = \sum_{i=1}^k \hat{\xi}_i $, $ k \ge 1 $.	By Theorem~\ref{ThIP}, after redefining the vector time series on a richer probability space, there exists a standard Brownian motion $ \{B_t : t \ge 0 \} $ such that
	\[
		\max_{k\le m} \left| S_k - B_k \right| = o_\P( m^{\frac{1}{\nu}} ) \qquad \text{and} \qquad \max_{k\le m} \left| \hat{S}_k - B_k \right| = o_\P( m^{\frac{1}{\nu}} )
	\]
	with $ \nu $ as specified in the theorem. Next, invoking the Skorohod representation theorem shows that one can define, on a new probability space, random variables $ \hat{\xi}_i', i \ge 1 $,  together with a standard Brownian motion $ B_m' $ such that
	$ \{ \hat{\xi}_i : i \ge 1 \} \stackrel{d}{=} \{ \hat{\xi}_i' : i \ge 1 \} $ and with 
	$ \hat{S}_k' = \sum_{i=1}^k \hat{\xi}_i' $
	\[
		m^{-\frac{1}{\nu}} \left| \hat{S}_m' - B_m' \right| = o(1),  \qquad a.s..	 
	\]
	Therefore,  the assumptions of Lemma~\ref{lemma1} are satisfied, and we can assume that there exist, for each $m$, two independent standard Brownian motions $ B_{1,m} $ and $ B_{2,m} $, such that
	\[
	| \hat{S}_m' - B_{1,m} | = o( m^{ \frac{1}{\nu} } )
	\]
	and
	\[
	\sup_{k \ge \delta m} k^{-\frac{1}{\nu}} |  \hat{S}_{m+k}' - \hat{S}_m' - B_{2,m} | = O_\P(1).
	\]
	The rest of the proof can be carried out along the lines of \cite{HHKS2004} which some modifications. For sake of clarity, we provide details and repeat some arguments. Consider 
	\begin{align*}
		&\sup_{k \ge m\delta} 
		\frac{ \left| \sum_{m < i \le m + k} \hat{\xi}_i' - \frac{k}{m} \sum_{i=1}^m \hat{\xi}_i' -  \left( B_{1,m}(k) - \frac{k}{m} B_{2,m}(m) \right)  \right| }{ g(m,k) } \\
		& = \sup_{k \ge  m\delta} 
		\frac{ O_\P( k^{\frac{1}{\nu}} ) + \frac{k}{m} O_\P( m^{\frac{1}{\nu}} ) }{ m^\frac{1}{2} \left( 1 + \frac{k}{m} \right) \left( \frac{k}{m+k} \right)^\gamma } \\
		& = O_\P(1) \sup_{k \ge \delta m} \frac{ k^{\frac{1}{\nu}} + \frac{k}{m} m^{\frac{1}{\nu}} }{   m^\frac{1}{2} \left( 1 + \frac{k}{m} \right) \left( \frac{k}{m+k} \right)^\gamma } \\
		& = o_\P(1),
	\end{align*}
Let $ \hat{Q}'(m,k) = \left| \sum_{m < i\le m + k} \hat{\xi}_i' - \frac{k}{m} \sum_{j=1}^m \hat{\xi}_j' \right| $ and define 
\[
	  \tilde{Q}(m,k) = \left| B_{1,m}(k) - \frac{k}{m} B_{2,m}(m) \right|.
\]
We have shown
\begin{equation}
	\label{Approx1}
	\left|\sup_{k \ge m\delta}  \frac{\hat{Q}'(m,k)}{\sigma_{0,\infty} g(m,k)}  - \sup_{k \ge m\delta} \frac{\tilde{Q}(m,k)}{g(m,k)} \right| = o_\P(1).
\end{equation}
Let us clarify the  distributional limit of $ \sup_{k \ge m\delta} \frac{\tilde{Q}(m,k)}{g(m,k)}  $. 
We obtain for each fixed $ T> 0 $
%
%
\begin{align*}
	\max_{\delta m \le k \le Tm} \frac{\tilde{Q}(m,k)}{g(m,k)} 
	& = \max_{t \in \left\{ \frac{k}{m} : \delta m \le k \le Tm\right\}} \frac{\left| \frac{1}{\sqrt{m}} B_{1,m}\left(t m\right) - \frac{k}{m} \frac{1}{\sqrt{m}} B_{2,m}\left( m \right) \right|}{\left(1+t \right) \left( \frac{t}{t+1} \right)^\gamma} \\
	& =\sup_{t \in [\delta, T]} \frac{|\frac{1}{\sqrt{m}} B_{1,m}(t m) - t \frac{1}{\sqrt{m}} B_{2,m}(m)|}{(1+t)\left(\frac{t}{1+t} \right)^\gamma} + o(1),
\end{align*}
as $ m \to \infty $, a.s., by Levy's continuity modulus. Since $ \frac{1}{\sqrt{ m}}(B_{1,m}(t m), B_{2,m}(m)) \stackrel{d}{=} (B_1(t), B_2(1))$ for all $m \ge 1 $,
where $ B_1, B_2$ are independent standard Brownian motions, we may conclude that
\[
 \P\left(  \left| \max_{\delta m \le k \le Tm} \frac{\tilde{Q}(m,k)}{g(m,k)}  -  \sup_{t \in [\delta, T]} \frac{| B_{1}(t ) - t  B_{2}(1)|}{(1+t)\left(\frac{t}{1+t} \right)^\gamma} \right| \to 0 \right) = 1.
\]
The law of the iterated logarithm yields
\[
  \sup_{k > mT} \frac{\left| B_{1}\left(\frac{k}{m} \right) - \frac{k}{m} B_{2}\left( \frac{k}{m} \right) \right|}{\left(1+\frac{k}{m}\right) \left( \frac{k}{k+m} \right)^\gamma}  = o_\P(1),
\]
and therefore we arrive at
\[
	\max_{\delta m \le k < \infty} \frac{\tilde{Q}(m,k)}{g(m,k)} =
	\sup_{\delta \le t < \infty}  \frac{| B_{1}(t ) - t  B_{2}(1)|}{(1+t)\left(\frac{t}{1+t} \right)^\gamma} + o_\P(1).
\]
Combining that results with \eqref{Approx1} gives
\[
	\max_{\delta m \le k < \infty} \frac{{Q}(m,k)}{\sigma_{0,\infty} g(m,k)} =
\sup_{\delta \le t < \infty}  \frac{| B_{1}(t ) - t  B_{2}(1)|}{(1+t)\left(\frac{t}{1+t} \right)^\gamma} + o_\P(1),
\]
and an application of the second half of the Skorohod representation theorem, yields, on the original probability space,
\[
\max_{\delta m \le k < \infty} \frac{Q(m,k)}{\sigma_{0,\infty}g(m,k)} \stackrel{d}{\to} \sup_{\delta \le t < \infty}\frac{|B_1(t) - t B_2(1)|}{(1+t)\left(\frac{t}{1+t} \right)^\gamma},
\]
as $ m \to \infty $. Since $ \{ B_1(t) - t B_2(1) : t \ge 0 \} \stackrel{d}{=} \{ (1+t) B(t/(1+t)) : t \ge 0 \} $, we may conclude that 
\[
\sup_{\delta \le t < \infty}\frac{|B_1(t) - t B_2(1)|}{(1+t)\left(\frac{t}{1+t} \right)^\gamma} \stackrel{d}{=} \sup_{ \frac{\delta}{1+\delta} \le t \le 1 } \frac{|B(t)|}{t^\gamma},
\]
where $B(t) $, $ 0 \le t \le 1 $, is a standard Brownian motion. Therefore,
\[
\sup_{k \ge m\delta}  \frac{Q(m,k)}{\sigma_{0,\infty} g(m,k)}  \stackrel{d}{\to} \sup_{ \frac{\delta}{1+\delta} \le t \le 1 } \frac{|B(t)|}{t^\gamma},
\]
and if $ \hat{\sigma}_{0,\infty}^2 \stackrel{\P}{\to} \sigma_{0,\infty}^2 $, $m \to \infty $, then we also have
\[
	\sup_{k \ge m\delta}  \frac{Q(m,k)}{\hat{\sigma}_{0,\infty} g(m,k)}  \stackrel{d}{\to} \sup_{ \frac{\delta}{1+\delta} \le t \le 1 } \frac{|B(t)|}{t^\gamma}.
\]
This completes the proof.
\end{myproof}

\begin{myproof}{of Theorem~\ref{ErdoesKacDetector}}
	Define
	\[
	  \calQ_m(t) = \frac{Q(m, \lfloor mt \rfloor )}{\sqrt{m} \hat \sigma_{0,\infty} (1+t)}, \qquad t \ge 0.
	\]
	It follows from the proof of Theorem~\ref{ThH0} that 
	\[
	  \left\{ \frac{\calQ_m(t) } {g_a( t/(1+t) )}  : t \ge \delta \right\} \Rightarrow \left\{ \frac{B(t/(1+t) )}{g_a( t/(1+t) )}   : t \ge \delta \right\},
	\]	
	as $ m \to \infty $, for some standard Brownian motion $ B$, where $ g_a(x) = \sqrt{(x+1) (a^2+ \log(x+1)} $. Therefore, by transforming time using the inverse $ h^{-1}: [\delta/(1+\delta),1) \to [\delta,\infty) $, $ t \mapsto t/(1-t) $, of $ h : [\delta,\infty) \to [\delta/(1+\delta),1) $, $  t \mapsto t/(1+t) $, we obtain
	\[
	  \left\{ \frac{\calQ_m( t/(1-t)) }{g_a(t)} :  t \in [\delta/(1+\delta),1) \right\} \Rightarrow \left\{ \frac{B(t)}{g_a(t)} : t \ge \delta \right\},
	\]
	such that we can conclude
	\begin{align*}
	 \lim_{m \to \infty} \P( \calQ_m(k/(m-k)) > g_a(k/m), \exists \delta m \le k < m  ) &= \P( B(t) > g_a(t), \exists t \ge \delta )\\
	 & \nearrow \P( B(t) > g_a(t), \exists t > 0 )\\
	 & = e^{-a^2/2},
	\end{align*}
	 as $ \delta \searrow 0 $. Note that $ 1+ t/(1-t) = 1/(1-t) $, so that 
	\begin{align*}
	  \calQ_m\left(  \frac{k}{m-k} \right) &= \left(1 - \frac{k}{m} \right)\frac{1}{\sqrt{m}\hat \sigma_{0,\infty}} \left( \sum_{i=m+1}^{m+\lfloor m \frac{k}{m-k} \rfloor} \xi_i - \frac{\lfloor m \frac{k}{m-k} \rfloor}{m} \sum_{j=1}^m \xi_j  \right) \\
	  & = \frac{m-k}{m} \frac{1}{\sqrt{m}\hat\sigma_{0,\infty}} \left( \sum_{i=m+1}^{m+\lfloor m \frac{k}{m-k} \rfloor} \xi_i - \frac{\lfloor m \frac{k}{m-k} \rfloor}{m} \sum_{j=1}^m \xi_j  \right)
	\end{align*}
	Hence, the event that an alarm is raised at physical time $k \ge \delta m $ where
	\[
	  \left|\calQ_m\left(  \frac{k}{m-k} \right) \right| > \sqrt{\frac{m+k}{m}\left(a^2 + \log(\frac{m+k}{m})\right)}
	\]
	is equivalent to raising an alarm at transformed time $ \delta m \le k < m $ where
	\begin{align*}
	\left|   \sum_{i=m+1}^{m+ \lfloor m \frac{k}{m-k} \rfloor} \xi_i - \frac{\lfloor m \frac{k}{m-k} \rfloor}{m} \sum_{j=1}^m \xi_j \right| &> 
	\hat \sigma_{0,\infty} \sqrt{m} \frac{m}{m-k} \sqrt{ \frac{m+k}{m} } \sqrt{a^2+ \log \frac{m+k}{m}},
	\end{align*}
	i.e., after rearranging terms,
	\begin{align*}
	\left|   \sum_{i=m+1}^{m+ \lfloor k \frac{m}{m-k} \rfloor} \xi_i - \frac{\lfloor k \frac{m}{m-k} \rfloor}{m} \sum_{j=1}^m \xi_j \right| &> 
	\hat\sigma_{0,\infty} \sqrt{m+k} \frac{m}{m-k} \sqrt{a^2+ \log \frac{m+k}{m} }.
	\end{align*}
	This means, when for (transformed  time) $k$ the threshold is exceeded for the first time, this occurs when the centered partial sum using the first $ \ell = \lfloor m \frac{k}{m-k} \rfloor \ge k $  observations after $m$ crosses the boundary. Vice versa, when a signal is given when the centered partial sum using the first $ k \ge \delta m$  observations (physical time) crosses the boundary, this corresponds to $k^* = \lceil k \frac{m}{m+k} \rceil $ (solve $ k \le k^* \frac{m}{m-k^*} < k + 1 $ for $k^*$), i.e.,
	\[
	\left|   \sum_{i=m+1}^{m+k} \xi_i - \frac{k}{m} \sum_{j=1}^m \xi_j \right| > 
	\hat\sigma_{0,\infty} \sqrt{m+k^*} \frac{m}{m-k^*} \sqrt{a^2+ \log \frac{m+k^*}{m} },
	\]
	which verifies \eqref{DetectorED}.
\end{myproof}
\color{black}

\begin{myproof}{of Theorem~\ref{ThH1}}

	 By Assumption A--iii we have $ \matM_{i} = \matM_i^0 $ for $ 1 \le i \le m $, where $ \matM_i^0 $ guarantees $ \vecv^\top \matM_i^0 \vecv = \theta_0 $, and, 	under the sequence of local alternatives for the moment functional,
	 \[
	   \vecv^\top \matM_{m+k} \vecv = \theta_0 +  \frac{\Delta\left(\frac{k-k^*}{m}\right)}{\sqrt{m}} \eins_{ k \ge k^* }, \qquad k \ge 1.
	 \]
	 Therefore
	 \[
	   \sum_{i=m+1}^{m+k} (\vecv^\top \matS_i \vecv - \theta_0) 
	   = \sum_{i=m+1}^{m+k} \vecv^\top (\matS_i - \matM_i) \vecv + \eins_{ k \ge k^* } \frac{1}{\sqrt{m}} \sum_{i=k^*}^{k} \Delta\left(\frac{i-k^*}{m}\right)
	 \]
	 and 
	 \begin{align*}
	   \sum_{i=m+1}^{m+k} \vecv^\top \matS_i \vecv - \frac{k}{m} \sum_{j=1}^m \vecv^\top \matS_j \vecv &=
	    \sum_{i=m+1}^{m+k} (\vecv^\top \matS_i  \vecv - \theta_0 ) - \frac{k}{m} \sum_{j=1}^m (\vecv^\top  \matS_j \vecv - \theta_0)\\
	   & = \sum_{i=m+1}^k \vecv^\top (\matS_i - \matM_i) \vecv  - \frac{k}{m} \sum_{j=1}^m \vecv^\top  (\matS_j - \matM_j) \vecv  \\
	   & \qquad + \eins_{k\ge k^*} \frac{1}{\sqrt{m}} \sum_{i=k^*}^k \Delta\left(\frac{i-k^*}{m}\right).
	 \end{align*}
    Since $ \Delta $ has bounded total variation $ TV( \Delta; [0,T] ) $ on $ [0,T] $, we have for $ k^* \le k $
 	\[
 	   \left| \frac{1}{\sqrt{m}} \sum_{i=k^*}^k \Delta \left(\frac{i-k^*}{m}\right) 
 	  - \sqrt{m}  \int_0^{\frac{k-k^*}{m}} \Delta(s) \, ds \right| 
 	  \le \frac{TV(\Delta; [0, T]) }{\sqrt{m}}
 	\]
 	and thus
 	\begin{align*}
 		\label{ApproxByIntegral}
 		  \max_{\delta m \le k \le Tm} \left| 
 		  \frac{\frac{\eins_{k\ge k^*} }{\sqrt{m}} \sum_{i=k^*}^k \Delta \left(\frac{i-k^*}{m}\right) 
 		  - \sqrt{m} \eins_{k\ge k^*} \int_0^{\frac{k-k^*}{m}} \Delta(s) \, ds}
 	  	{\sqrt{m}\left(1 + \frac{k}{m} \right) \left( \frac{k}{k+m} \right)^\gamma } \right| 
 		  & \le \frac{TV(\Delta; [0, T]) }{ (1+\delta) \frac{\delta}{T+1} m}.
	\end{align*}
	Arguing as in the proof of Theorem~\ref{ThH0}  
 	we obtain for $ \hat{Q}'(m,k) = \sum_{m < i \le m+k} \hat{\xi}_i '- \frac{k}{m} \sum_{j=1}^m \hat{\xi}_j' $
 	\[
 	  \left| \max_{\delta m \le k \le Tm} \frac{\hat{Q}'(m,k) }{\sigma_{0,\infty}^2g(m,k)} -
 	    \max_{\delta m \le k \le Tm} \frac{B_{1,m}(k) - \frac{k}{m} B_{2,m}(m) + \sqrt{m} \eins_{k\ge k^*} \int_0^{\frac{k-k^*}{m}} \Delta(s) \, ds }{g(m,k)}
 	  \right| = o_\P(1).
 	\]
 	Further, 
 	\begin{align*}
 		& \max_{\delta m \le k \le Tm} \frac{B_{1,m}(k) - \frac{k}{m} B_{2,m}(m) + \sqrt{m} \int_0^{\frac{k-k^*}{m}} \Delta(s) \, ds }{g(m,k)} \\
 		& \qquad  \stackrel{d}{=} \max_{\delta m \le k \le Tm}\frac{\left| B_{1}\left(\frac{k}{m} \right) - \frac{k}{m} B_{2}\left( 1 \right) + \eins_{k\ge k^*} \int_0^{(k-k^*)/m} \Delta(s) \, ds \right|}{\left(1+\frac{k}{m}\right) \left( \frac{k}{k+m} \right)^\gamma} \\
 		& \qquad \stackrel{d}{\to} \sup_{t \in [\delta, T]} \frac{|B_1(t) - t B_2(1) + \eins_{t \ge \vartheta } \int_0^{t-\vartheta} \Delta(s) \, ds |}{(1+t)\left(\frac{t}{1+t} \right)^\gamma},
	\end{align*}
 which completes the proof. 
\end{myproof}

The following proofs draw to some extent on \cite{BickelLevina2008}. 

\newcommand{\Sihat}{{\Gamma}_{ij}}
\newcommand{\Si}{{\Sigma_{ij}}}
 
\begin{myproof}{of Theorem~\ref{ThBoundThresholdOp}}
First observe that
\begin{align*}
  \|\calT_t \bfSigma  - \bfSigma \|_{op} 
  &\le \max_{1 \le i \le d} \sum_{j=1}^d | \Sigma_{ij} | \eins( | \Sigma_{ij} | \le t) 
   \le \max_{1 \le i \le d} \sum_{j=1}^d  | \Sigma_{ij} |^r t^{1-r} 
  \le t^{1-r} s_o \\
\end{align*}
Clearly, $ \| \calT_t  {\bfGamma}  - \calT_t  \bfSigma  \|_{op} \le \max_{1 \le i \le d} \sum_{j=1}^d  |\Sihat  \eins_{|\Sihat | \ge t} - \Sigma_{ij} \eins_{|\Sigma_{ij}| \ge t } | $. We have the decomposition
\begin{align*}
	  \Sihat \eins_{|\Sihat|\ge t} - \Si \eins_{|\Si| \ge t} 
	  & = ( \Sihat - \Si ) \eins_{|\Sihat| \ge t, |\Si| < t} + ( \Sihat - \Si ) \eins_{|\Sihat| \ge t, |\Si| \ge t} \\
	  & \qquad + \Si \left( \eins_{|\Sihat| \ge t, | \Si | < t}  - \eins_{|\Si| \ge t} +
	  \eins_{|\Sihat| \ge t, \Si \ge t } \right) \\
	  & = ( \Sihat - \Si ) \eins_{|\Sihat| \ge t, |\Si| < t} + ( \Sihat - \Si ) \eins_{|\Sihat| \ge t, |\Si| \ge t} \\
	& \qquad + \Si \left( \eins_{|\Sihat| \ge t, | \Si | < t}  - \eins_{|\Sihat| < t, |\Si| \ge t} \right) \\
	& = U_{ij}^1 + U_{ij}^2 + U_{ij}^3. 
\end{align*}
Fix $ 0 < \gamma < 1$. Since $ | \Sihat - \Si | > (1-\gamma)t $, if $ |\Sihat|\ge t $ and $ |\Si|< \gamma t$,
\begin{align*}
	\max_{1 \le i \le d} \sum_{j=1}^d | U_{ij}^1 |
	& \le \max_{1 \le i \le d} \sum_{j=1}^d  |\Sihat - \Si | \eins_{|\Sihat|\ge t, |\Si|< \gamma t}
	+ \max_{1 \le i \le d} \sum_{j=1}^d |\Sihat - \Si | \eins_{|\Sihat|\ge t, \gamma t \le |\Si | < t} \\
	& \le \| \bfGamma - \bfSigma \|_\infty \left( \#( |\Sihat - \Si | > (1-\gamma) t ) +
	\max_{1 \le i \le d} \sum_{j=1}^d  \frac{ |\Si |^r}{(\gamma t)^r} \eins( |\Si | \ge \gamma t) \right) \\
	& \le \| \bfGamma - \bfSigma \|_\infty \left( \#( |\Sihat - \Si | > (1-\gamma) t ) + \frac{1}{(\gamma t)^r} s_0  \right).
\end{align*}
Further,
\begin{align*}
	\max_{1 \le i \le d} \sum_{j=1}^d | U_{ij}^2 |
	& \le \| \bfGamma - \bfSigma \|_\infty \max_{1 \le i \le d} \sum_{j=1}^d \frac{|\Si|^r}{t^r} \eins_{|\Si| > t} 
	 \le \frac{s_0}{t^r} \| \bfGamma - \bfSigma \|_\infty.
\end{align*}
To bound $ \max_{1 \le i \le d} \sum_{j=1}^d | U_{ij}^3 | $ first observe that
\begin{align*}
  \max_{1 \le i \le d} \sum_{j=1}^d  |\Si| \eins_{|\Sihat| \ge t, | \Si | < t } 
  & \le \max_{1 \le i \le d} \sum_{j=1}^d | \Si |^r | \Si |^{1-r} \eins_{|\Si| < t} 
  \le t^{1-r} s_0.
\end{align*}
Lastly,
\begin{align*}
	\max_{1 \le i \le d} \sum_{j=1}^d | \Si | \eins_{|\Sihat| < t, | \Si | \ge t} 
	& \le \max_{1 \le i \le d} \sum_{j=1}^d | \Sihat - \Si | \eins_{|\Sihat| < t, | \Si | \ge t} 
	+ \max_{1 \le i \le d} \sum_{j=1}^d  | \Sihat | \eins_{|\Sihat| < t, |\Si | \ge t} \\
	& \le \| \bfSigma - \bfSigma \|_\infty \max_{1 \le i \le d} \sum_{j=1}^d \frac{|\Si|^r}{t^r} \eins_{|\Si| \ge t} + \max_{1 \le i \le d} \sum_{j=1}^d t \frac{|\Si|^r}{t^r} \eins_{|\Si| > t} \\
	& \le \| \bfGamma - \bfSigma \|_\infty  \frac{s_0}{t^r} + t^{1-r} s_0.
\end{align*}
Combining the above bounds and collecting term shows that
\begin{align*}
\| \calT_t {\bfGamma}- \bfSigma  \|_{op} 
& \le \| \calT_t \bfSigma  - \bfSigma \|_{op} + \| \calT_t \bfGamma - \calT_t \bfSigma  \|_{op} \\
& \le 2 t^{1-r} s_0 + \| {\bfGamma} - \bfSigma \|_\infty
\left( \#( |{\Gamma}_{ij} - \Sigma_{ij}| > (1-\gamma)t  ) + \frac{1}{(\gamma t)^r} s_0 + \frac{2}{t^r} s_0 \right).
\end{align*}
For the soft-thresholding operator we have
\[
  | \calS_t( \Gamma_{ij} ) - \calS_t( \Sigma_{ij} ) | 
  \le \begin{cases}
  	  | \calS_t( \Gamma_{ij} ) | \le | \Gamma_{ij} | \le | \Gamma_{ij} - \Sigma_{ij} | + | \Sigma_{ij} |, & | \Gamma_{ij} | \ge t, | \Sigma_{ij} | < t, \\
  	  | \Sigma_{ij} |, & |\Gamma_{ij} | < t, | \Sigma_{ij} | \ge t, \\
  	  2 t + | \Gamma_{ij} - \Sigma_{ij} |, & | \Gamma_{ij} | \ge t, | \Sigma_{ij} | \ge t,
     \end{cases}
\]
where in the last case we used the estimate
\begin{align*}
  | \calS( \Gamma_{ij} ) - \calS( \Sigma_{ij} ) | 
  &\le | \calS_t( \Gamma_{ij} ) - \Gamma_{ij} | + | \Gamma_{ij} - \Sigma_{ij} | + | \calS_t( \Sigma_{ij} ) - \Sigma_{ij} | \\
  & \le 2 t +  | \Gamma_{ij} - \Sigma_{ij} |.
\end{align*}
Therefore,
\begin{align*}
  \max_{1 \le i \le d} \sum_{j=1}^d | \calS_t( \Gamma_{ij} ) - \calS_t( \Sigma_{ij} ) |
  & \le \max_{1 \le i \le d} \sum_{j=1}^d | \Gamma_{ij} - \Sigma_{ij} | ( \eins_{|\Gamma_{ij}|\ge t, |\Sigma_{ij}| < t} + \eins_{|\Gamma_{ij}| \ge t, |\Sigma_{ij}| \ge t} ) \\
  & \qquad + \max_{1 \le i \le d} \sum_{j=1}^d | \Sigma_{ij} | ( \eins_{|\Gamma_{ij}| \ge t, |\Sigma_{ij}| < t} + \eins_{|\Gamma_{ij}| < t, |\Sigma_{ij}| \ge t}) \\
  & = U_{ij}^1 + U_{ij}^2 + U_{ij}^3,
\end{align*}
such that the same bound results. 
\end{myproof}

\begin{myproof}{of Theorem~\ref{ThConsThresholdedCov} and Corollary~\ref{ConvRate}}
	We need to control the upper bound on $ \| T_s \hat{\bfSigma}_m - \bar \bfSigma_m \|_{op} $. 
	By Lemma~\ref{LemmaRatesMeanCovEst}, Markov's inequality yields for $ t > 0 $
	\[
	  \P( | \hat{\Sigma}_{m,jk} - \bar \Sigma_{m,jk} | \ge t  ) \le  
	  \frac{\E | \hat{\Sigma}_{m,jk} - \bar \Sigma_{m,jk} |^{\frac{q}{2}} }{ t^{\frac{q}{2}} } \le \frac{ c m^{-\frac{q}{4}}} {t^{\frac{q}{2}} }
    \]
   such that the union bound yields
    \begin{align*}
    	\P\left( d^{-\frac{4}{q}} m^{\frac{1}{2}} 
    	 \max_{1  \le j,k \le d} | \hat{\Sigma}_{m,jk} - \bar \Sigma_{m,jk} | \ge t  \right) \le \frac{c}{t^{\frac{q}{2}} }.
    \end{align*}
    Therefore,
	\[
	  \| \hat{\bfSigma}_m - \bar \bfSigma_m \|_\infty = O_\P\left( \frac{d^{\frac{4}{q}}}{ \sqrt{m} } \right).
	\]
	Clearly, the proof carries over to $ \hat{\matM}_m $, and since $ \max_{1 \le j \le d} \E | \hat{\mu}_{mj} - \bar\mu_{mj} |^q = O( m^{q/2} )$, one gets $ \P( \frac{\sqrt{m}}{d^{2/q}} \| \hat{\bfmu}_m - \bar{\bfmu}_m \|_\infty  \ge t ) = O( t^{-q/2} ) $, so that $ \| \hat{\bfmu}_m - \bar{\bfmu}_m \|_\infty = O_\P( \frac{d^{2/q}}{\sqrt{m}} ) $ follows. Bounding the Frobenius norm by $ d $ times the supnorm gives
	\begin{align*}
	  \P\left( d^{-\frac{4+q}{q}} m^{\frac{1}{2}} \| \hat{\Sigma}_{m,jk} - \bar \Sigma_{m,jk} \|_F \ge t  \right) &\le 
	   \P\left(  d^{-\frac{4}{q}} m^{-\frac{1}{2}}  
	  \max_{1  \le j,k \le d} | \hat{\Sigma}_{m,jk} - \bar \Sigma_{m,jk} | \ge t \right) \\
	  & \le c \frac{d^2 m^{-\frac{q}{4}} }{t^{\frac{q}{2}} m^{-\frac{q}{4}} d^{2}} 
	  = \frac{c}{t^{\frac{q}{2}}},
	\end{align*}
	and $ \P(d^{\frac{2-q}{q}} m^\frac{1}{2} \| \hat \bfmu_m  - \bar \bfmu_m \|_2 \ge t )
	\le d \P( d^{2/q}  m^{1/2}  \| \hat \bfmu_m  - \bar \bfmu_m \|_\infty \ge t  ) = O(t^{-q/2} )$. Thus, Corollary~\ref{ConvRate} is shown. Next, we have
	\[
	  p_0 = \P\left( \#\left( | \hat{\Sigma}_{m,ij} - \bar\Sigma_{m,ij} | > (1-\gamma) t \right) > 0 \right) = \P( \| \hat{\bfSigma}_m - \bar\bfSigma_m \|_\infty > (1-\gamma)t )
	\]
	Thus, if one selects the threshold as a multiple of the convergence rate of $ \| \hat{\bfSigma}_m - \bar \bfSigma_m \|_\infty $, i.e.,
	$ t_{th} = C_{th} \frac{d^{\frac{4}{q}}}{ \sqrt{m}} $, for some constant $ C_{th} $, then $ p_0 $ is arbitrarily small if $ C_{th} $ is large enough.  Now, recalling that $ 0 \le r < 1 $, it easily follows by noting that $ \|  \hat{\bfSigma}_m - \bar \bfSigma_m \|_\infty	 = O_P( t_{th}^{1-r} ) $ that 
	\[
	  \| \calT_{th} \hat{\bfSigma}_m  - \bar \bfSigma_m \|_{op} = O( t_{th}^{1-r} ) + O_\P( \| \hat{\bfSigma}_m - \bar \bfSigma_m \|_\infty ) = O_\P(  t_{th}^{1-r}  )
	\]
  which establishes Theorem~\ref{ThConsThresholdedCov}.  
\end{myproof}

\color{\blau}
\begin{myproof}{of Theorem~\ref{ThConsPrecision}}
	Observe that the thresholded matrix $ \calS_{t_m} \hat{\bfSigma}_m  $, for some sequence $ t_m = o(1) $ of thresholds, has maximal eigenvalue bounded by $ \max_i \sum_j |\bar \Sigma_{m,ij} |\le M^{1-r} s_0 $ as well,  and that $ \lambda_{\min}( \calS_{t_m} \hat{\bfSigma}_m  ) \ge \frac{1}{2} \lambda_{\min}( \bar\bfSigma_m ) \ge \frac{\varepsilon_0}{2} $ if $ t_m $ is small enough to ensure that $ \| \calS_{t_m} - \calS_0 \|_{op} < \inf_{m \ge 1} \lambda_{\min}( \bar \bfSigma_m )  $, since the arguments of \cite[p.2580]{BickelLevina2008} carry over. Therefore $ \calS_{t_m} \hat{\bfSigma}_m \in \mathcal{U}( s, s_0, M, \varepsilon_0)  $. Consequently, the well known inequality
	\[
	\| \matA^{-1} - \matB^{-1} \|_{op} \le \frac{  \| \matA - \matB \|_{op} }{ \lambda_{\min}( \matA ) \lambda_{\min}( \matB )} 
	\]
	for regular square matrices $ \matA,  \matB $ of the same dimension,  yields for the inverses, under the assumptions of Theorem~\ref{ThConsThresholdedCov},
	\[
	\|( \calS_{t_m} \hat{\bfSigma}_m  )^{-1} - \bar\bfSigma_m^{-1} \|_{op}
	= O_\P\left(\left( \frac{d^{\frac{4}{q}}}{\sqrt{m}} \right)^{1-r}\right).
	\]
\end{myproof}
\color{black}

\begin{myproof}{of Proposition~\ref{LipschitzDeepLearner}}
	If $ \sigma : \R \to \R $ is $ \rho $--Lipschitz, then $ \vecx \mapsto \sigma( \vecx )$, $ \vecx  $ some vector, is $ \rho $--Lipschitz with respect to any $p$--vector norm. Hence, for a $ \rho $--Lipschitz activation $ \sigma $ with $ \sigma(0) = 0 $, $N \times r $ weighting matrices $ \matW, \wt{\matW} $ and $ \vecx \in \R^r $, $ N, r \in \N $,
	\begin{itemize}
		\item[(i)] $ | \sigma^{\wt{\matW}}(\vecx ) - \sigma^{\matW}(\vecx ) | \le \rho \| (\wt{\matW} - \matW)\vecx  \|_2 \le \rho \|\wt{\matW} - \matW \|_{op} \| \vecx \|_2 $,
		\item[(ii)] $ \| \sigma^{\matW} ( \wt{\vecx} ) - \sigma^{\matW} ( \vecx ) \|_2 \le \rho \| \matW (\wt{\vecx}- \vecx ) \|_2 \le \rho \| \matW \|_{op} \| \wt{\vecx} - \vecx \|_2$,
		\item[(iii)] $ \|\sigma^{\matW}( \vecx ) \|_2 = \| \sigma^{\matW} ( \vecx )  - \sigma(0) \|_2 \le \rho \| \matW \|_{op} \| \vecx \|_2 $, and therefore
		\item[(iv)]  $ \| \vecf_j( \vecx ) \|_2 = \| \sigma_j^{\matW_j} \circ \cdots \circ \sigma_1^{\matW_1}( \vecx ) \|_2 \le \prod_{i=1}^j \rho_i \| \matW_i \|_{op} \| \vecx \|_2 $.
	\end{itemize}
	We have for $ k <  j \le H $
	\begin{align*}
		& \| \vecf_j( \vecx; \bftheta(\wt{\matW}_k) ) - \vecf_j( \vecx; \bftheta( \wt{\matW}_{k-1} ) ) \|_2 \\
		& \quad = 
		\| \sigma_j^{\matW_j} \circ \cdots \circ \sigma_{k+1}^{\matW_{k+1}} 
		\circ \sigma_{k}^{\wt{\matW}_{k}} 
		\circ \cdots \circ \sigma_1^{\wt{\matW}_1}(\vecx)
		- \sigma_j^{\matW_j}\circ \cdots \circ \sigma_k^{\matW_k} \circ \sigma_{k-1}^{\wt{\matW}_{k-1}} \circ \cdots \circ \sigma_1^{\wt{\matW}_1} (\vecx) \|_2 \\
		& \quad  \stackrel{(ii)}{\le} \rho_j \| \matW_j \|_{op} 
		\| \sigma_{j-1}^{\matW_{j-1}} \circ \cdots \circ \sigma_{k+1}^{\matW_{k+1}} 
		\circ \sigma_{k}^{\wt{\matW}_{k}} 
		\circ \cdots \circ \sigma_1^{\wt{\matW}_1}(\vecx) \\
		& \qquad \qquad 
		- \sigma_{j-1}^{\matW_{j-1}}\circ \cdots \circ \sigma_k^{\matW_k} \circ \sigma_{k-1}^{\wt{\matW}_{k-1}} \circ \cdots \circ \sigma_1^{\wt{\matW}_1} (\vecx) \|_2 \\
		& \quad \le \cdots \le \prod_{i=k+1}^j  \rho_i \| \matW_i \|_{op} 
		\| \sigma_k^{\wt{\matW}_k} \circ \cdots \circ \sigma_1^{\wt{\matW}_1}( \vecx ) 
		- \sigma_k^{\matW_k} \circ \sigma_{k-1}^{\wt{\matW}_{k-1}} \cdots \circ \sigma_1^{\matW_1}( \vecx ) \|_2 \\
		& \quad \stackrel{(i)}{\le} \left( \prod_{i=k+1}^j  \rho_i \| \matW_i \|_{op} \right) \rho_k \| \wt{\matW}_k - \matW_k \|_{op} 
		\|  \sigma_{k-1}^{\wt{\matW}_{k-1}} \circ \cdots \circ \sigma_1^{\wt{\matW}_1}( \vecx ) \|_2 \\
		& \quad \stackrel{(iv)}{\le}  \prod_{i=1, i \not= k}^j \rho_i  \| \matW_i \|_{op}   \rho_k \| \wt{\matW}_k - \matW_k \|_{op} \| \vecx \|_2 \\
		& \quad = L_{jk} \| \wt{\matW}_k - \matW_k \|_{op} \| \vecx \|_2
	\end{align*}
	where $ L_{jk} = \prod_{i=1, i \not= k}^j \rho_i  \| \matW_i \|_{op}   \rho_k $. 
	Notice that $ \bftheta( \wt{\matW}_0 ) = \bftheta $ and $ \bftheta( \wt{\matW}_H ) = \wt{\bftheta} $. By telescoping and estimating all spectral norms by the corresponding  Frobenius matrix norms, we obtain
	\begin{align*}
		\| \vecf_H( \vecx; \wt{\bftheta} ) - \vecf_H(\vecx; \bftheta ) \|_2 
		& = \left\| \sum_{k=1}^H f_H( \vecx; \theta( \wt{\matW}_k )  ) - f_H( \vecx; \theta( \wt{\matW}_{k-1} ) ) \right\|_2 \\
		& \le \sum_{h=1}^H L_{Hk} \| \wt{\matW}_k - \matW_{k} \|_{op} \|\vecx \|_2 \\
		& \le H \max_{1 \le k \le H} L_{Hk}   \|\vecx \|_2 \frac{1}{H} \sum_{k=1}^H \| \wt{\matW}_k - \matW_k \|_F\\
		& \le H L_H \|\vecx \|_2 \sqrt{ \frac{1}{H} \sum_{k=1}^H \| \wt{\matW}_k - \matW_k \|_F^2 } \\
		& = L_H \sqrt{H}  \|\vecx \|_2 \sqrt{\sum_{k=1}^H \| \wt{\matW}_k - \matW_k \|_F^2 } \\
		& = L_H \sqrt{H}  \|\vecx \|_2 \| \wt{\bftheta} - \bftheta \|_F,
	\end{align*}
	where the last inequality follows from Jensen's inequality.
\end{myproof}

\section*{Acknowledgements}
The author acknowledges support from Deutsche Forschungsgemeinschaft (DFG, grant STE 1034/11--2). He thanks M.Sc. Florian Scholze for carefully proof--reading an earlier version of the manuscript.

\bibliography{lit}

@article {BickelLevina2008,
    AUTHOR = {Bickel, Peter J. and Levina, Elizaveta},
     TITLE = {Covariance regularization by thresholding},
   JOURNAL = {Ann. Statist.},
  FJOURNAL = {The Annals of Statistics},
    VOLUME = {36},
      YEAR = {2008},
    NUMBER = {6},
     PAGES = {2577--2604},
      ISSN = {0090-5364},
     CODEN = {ASTSC7},
   MRCLASS = {62H12 (62F12 62G09)},
  MRNUMBER = {2485008 (2010b:62197)},
MRREVIEWER = {M. Hu{\v{s}}kov{\'a}},
}

@article{RothmanLevinaZhu2009,
	author = {Rothman, Adam J. and Levina, Elizaveta and Zhu, Ji},
	title = {Generalized Thresholding of Large Covariance Matrices},
	journal = {Journal of the American Statistical Association},
	volume = {104},
	number = {485},
	pages = {177-186},
	year = {2009},
	publisher = {Taylor & Francis},
	doi = {10.1198/jasa.2009.0101},
}

@book{WhiteHalbert1992,
  author = {White, Halbert},
  title = {Artificial Neural Networks},
  year = {1992},
  publisher = {Blackwell Publishers},
}

@article {RafaSteland2014,
AUTHOR = {Steland, Ansgar and Rafaj{\l}owicz, Ewaryst},
TITLE = {Decoupling change-point detection based on characteristic
functions: methodology, asymptotics, subsampling and
application},
JOURNAL = {J. Statist. Plann. Inference},
FJOURNAL = {Journal of Statistical Planning and Inference},
VOLUME = {145},
YEAR = {2014},
PAGES = {49--73},
ISSN = {0378-3758},
MRCLASS = {62L10 (62G10 62G20 62M07 62M10)},
MRNUMBER = {3125349},
MRREVIEWER = {Antonio Cuevas},
DOI = {10.1016/j.jspi.2013.08.009},
URL = {https://doi.org/10.1016/j.jspi.2013.08.009},
}

@article{JagMa2003,
Author = {Jagannathan, R and Ma, TS},
Title = {Risk reduction in large portfolios: Why imposing the wrong constraints
helps},
Journal = {Journal of Finance},
Year = {2003},
Volume = {58},
Number = {4},
Pages = {1651-1683},
Month = {AUG},
Publisher = {WILEY-BLACKWELL},
Address = {111 RIVER ST, HOBOKEN 07030-5774, NJ USA},
Type = {Article},
Language = {English},
Affiliation = {Jagannathan, R (Corresponding Author), Univ Utah, David Eccles Sch Business, Salt Lake City, UT 84112 USA.
Univ Utah, David Eccles Sch Business, Salt Lake City, UT 84112 USA.
Northwestern Univ, Kellogg Sch Management, Evanston, IL 60208 USA.
Natl Bur Econ Res, Cambridge, MA 02138 USA.},
DOI = {10.1111/1540-6261.00580},
ISSN = {0022-1082},
EISSN = {1540-6261},
Keywords-Plus = {VARIANCE-EFFICIENT PORTFOLIOS; PERFORMANCE; SELECTION; RETURNS; MODEL},
Research-Areas = {Business \& Economics},
Web-of-Science-Categories  = {Business, Finance; Economics},
Number-of-Cited-References = {28},
Times-Cited = {647},
Usage-Count-Last-180-days = {1},
Usage-Count-Since-2013 = {33},
Journal-ISO = {J. Financ.},
Doc-Delivery-Number = {706UK},
Web-of-Science-Index = {Social Science Citation Index (SSCI)},
Unique-ID = {WOS:000184473300012},
OA = {Green Published},
DA = {2022-08-29},
}

@article {FanZhangYu2012,
AUTHOR = {Fan, Jianqing and Zhang, Jingjin and Yu, Ke},
TITLE = {Vast portfolio selection with gross-exposure constraints},
JOURNAL = {J. Amer. Statist. Assoc.},
FJOURNAL = {Journal of the American Statistical Association},
VOLUME = {107},
YEAR = {2012},
NUMBER = {498},
PAGES = {592--606},
ISSN = {0162-1459},
MRCLASS = {91G10 (62P05)},
MRNUMBER = {2980070},
DOI = {10.1080/01621459.2012.682825},
URL = {https://doi.org/10.1080/01621459.2012.682825},
}

@article{ ZhaoLedoitJiang2021,
Author = {Zhao, Zhao and Ledoit, Olivier and Jiang, Hui},
year = {2021},
Title = {Risk Reduction and Efficiency Increase in Large Portfolios:
Gross-Exposure Constraints and Shrinkage of the Covariance Matrix},
Journal = {Journal of Financial Econometrics},
Publisher = {OXFORD UNIV PRESS},
DOI = {10.1093/jjfinec/nbab001},
EarlyAccessDate = {FEB 2021},
}

@article {Markowitz1952,
AUTHOR = {Markowitz, H.M.},
year = {1952},
title = {Portfolio Selection},
journal = {Journal of Finance},
pages = {77-91},
vol = {7},
}

@article {Sharpe1964,
AUTHOR = {Sharpe, W.},
year = {1964},
title = {Capital Asset Prices: A Theory of Market Equilibrium Under Conditions of Risks},
journal = {Journal of Finance},
pages = {425-442},
vol = {19},
}

@book {Markowitz1959,
AUTHOR = {Markowitz, H.M.},
year = {1959},
title = {Portfolio Selection: Efficient Diversification of Investments},
publisher = {John Wiley \& Son},  
}

@article {SeqMonSqrdRes2010,
AUTHOR = {Chen, Zhan Shou and Tian, Zheng and Qin, Rui Bing and Leng,
Cheng Cai},
TITLE = {Sequential monitoring variance change in linear regression
model},
JOURNAL = {J. Math. Res. Exposition},
FJOURNAL = {Journal of Mathematical Research and Exposition},
VOLUME = {30},
YEAR = {2010},
NUMBER = {4},
PAGES = {610--618},
ISSN = {1000-341X},
MRCLASS = {62L10 (62J05)},
MRNUMBER = {2742115},
}

@article {HHKS2004,
AUTHOR = {Horv\'{a}th, Lajos and Hu\v{s}kov\'{a}, Marie and Kokoszka, Piotr and
Steinebach, Josef},
TITLE = {Monitoring changes in linear models},
JOURNAL = {J. Statist. Plann. Inference},
FJOURNAL = {Journal of Statistical Planning and Inference},
VOLUME = {126},
YEAR = {2004},
NUMBER = {1},
PAGES = {225--251},
ISSN = {0378-3758},
MRCLASS = {62L10 (62J05)},
MRNUMBER = {2090695},
MRREVIEWER = {Rainer Schwabe},
DOI = {10.1016/j.jspi.2003.07.014},
URL = {https://doi.org/10.1016/j.jspi.2003.07.014},
}

@article {Wu2005,
AUTHOR = {Wu, Wei Biao},
TITLE = {Nonlinear system theory: another look at dependence},
JOURNAL = {Proc. Natl. Acad. Sci. USA},
FJOURNAL = {Proceedings of the National Academy of Sciences of the United
States of America},
VOLUME = {102},
YEAR = {2005},
NUMBER = {40},
PAGES = {14150--14154},
ISSN = {0027-8424},
MRCLASS = {62M10},
MRNUMBER = {2172215},
DOI = {10.1073/pnas.0506715102},
URL = {https://doi.org/10.1073/pnas.0506715102},
}

@article {Haviland1936,
AUTHOR = {Haviland, E. K.},
TITLE = {On the {M}omentum {P}roblem for {D}istribution {F}unctions in
{M}ore {T}han {O}ne {D}imension. {II}},
JOURNAL = {Amer. J. Math.},
FJOURNAL = {American Journal of Mathematics},
VOLUME = {58},
YEAR = {1936},
NUMBER = {1},
PAGES = {164--168},
ISSN = {0002-9327,1080-6377},
MRCLASS = {99-04},
MRNUMBER = {1507139},
DOI = {10.2307/2371063},
URL = {https://doi.org/10.2307/2371063},
}

@article {GhasemiEtAl2016,
AUTHOR = {Ghasemi, Mehdi and Kuhlmann, Salma and Marshall, Murray},
TITLE = {Moment problem in infinitely many variables},
JOURNAL = {Israel J. Math.},
FJOURNAL = {Israel Journal of Mathematics},
VOLUME = {212},
YEAR = {2016},
NUMBER = {2},
PAGES = {989--1012},
ISSN = {0021-2172,1565-8511},
MRCLASS = {44A60 (47A57)},
MRNUMBER = {3505409},
MRREVIEWER = {Akio\ Arimoto},
DOI = {10.1007/s11856-016-1318-5},
URL = {https://doi.org/10.1007/s11856-016-1318-5},
}

@article{ErdoesKac1946,
author = {P. Erd{\"o}s and M. Kac},
title = {{On certain limit theorems of the theory of probability}},
volume = {52},
journal = {Bulletin of the American Mathematical Society},
number = {4},
publisher = {American Mathematical Society},
pages = {292 -- 302},
year = {1946},
}

@article{RobbinsSiegmund1970,
author = {Herbert Robbins and David Siegmund},
title = {{Boundary Crossing Probabilities for the Wiener Process and Sample Sums}},
volume = {41},
journal = {The Annals of Mathematical Statistics},
number = {5},
publisher = {Institute of Mathematical Statistics},
pages = {1410 -- 1429},
year = {1970},
doi = {10.1214/aoms/1177696787},
URL = {https://doi.org/10.1214/aoms/1177696787}
}

@article{erdogmus2004minimax,
title={Minimax mutual information approach for independent component analysis},
author={Erdogmus, Deniz and Hild, Kenneth E and Rao, Yadunandana N and Principe, Jose C},
journal={Neural Computation},
volume={16},
number={6},
pages={1235--1252},
year={2004},
publisher={MIT Press}
}

@article{John2007,
title = {Techniques for the reconstruction of a distribution from a finite number of its moments},
author = {V. John and I. Angelov and A.A. \"{O}nc\"{u}l and D. Th\'{e}venin},
journal = {Chemical Engineering Science},
volume = {62},
number = {11},
pages = {2890-2904},
year = {2007},
issn = {0009-2509},
doi = {https://doi.org/10.1016/j.ces.2007.02.041},
url = {https://www.sciencedirect.com/science/article/pii/S0009250907002072},
}

@article {MeadPapaniolaou1984,
AUTHOR = {Mead, Lawrence R. and Papanicolaou, N.},
TITLE = {Maximum entropy in the problem of moments},
JOURNAL = {J. Math. Phys.},
FJOURNAL = {Journal of Mathematical Physics},
VOLUME = {25},
YEAR = {1984},
NUMBER = {8},
PAGES = {2404--2417},
ISSN = {0022-2488,1089-7658},
MRCLASS = {82A05 (44A60 82A15)},
MRNUMBER = {751523},
DOI = {10.1063/1.526446},
URL = {https://doi.org/10.1063/1.526446},
}

@article{ChuStinchcombeWhite1996,
  AUTHOR = {Chu, Chia-Shang James and Stinchcombe, Maxwell and
  White, Halbert},
  TITLE = {Monitoring structural change},
  JOURNAL = {Econometrica},
  YEAR = {1996},
  VOLUME = {64},
  NUMBER = {5},
  PAGES = {1045--1065},
}

@article {StelandJMVA2020,
AUTHOR = {Steland, Ansgar},
TITLE = {Testing and estimating change-points in the covariance matrix
of a high-dimensional time series},
JOURNAL = {J. Multivariate Anal.},
FJOURNAL = {Journal of Multivariate Analysis},
VOLUME = {177},
YEAR = {2020},
PAGES = {104582, 24},
ISSN = {0047-259X},
MRCLASS = {62E20 (60F17 62H99 62L10 62L12 62M10)},
MRNUMBER = {4068598},
DOI = {10.1016/j.jmva.2019.104582},
URL = {https://doi.org/10.1016/j.jmva.2019.104582},
}

@article {Sutradhar1986,
AUTHOR = {Sutradhar, Brajendra C.},
TITLE = {On the characteristic function of multivariate {S}tudent
{$t$}-distribution},
JOURNAL = {Canad. J. Statist.},
FJOURNAL = {The Canadian Journal of Statistics. La Revue Canadienne de
Statistique},
VOLUME = {14},
YEAR = {1986},
NUMBER = {4},
PAGES = {329--337},
ISSN = {0319-5724,1708-945X},
MRCLASS = {62H10 (62E10)},
MRNUMBER = {876759},
MRREVIEWER = {D.\ N.\ Shanbhag},
DOI = {10.2307/3315191},
URL = {https://doi.org/10.2307/3315191},
}

@inproceedings{Kingma2015,
AUTHOR = {Kingma, D. and Ba, J.},
TITLE = {Adam: A method for stochatic optimization},
BOOKTITLE = {Proceedings of the {3}rd International Conference on Learning Representation (ICLR 2015), San Diego.},
YEAR = {2015},
}

@article {RobbbinsMonro1951,
AUTHOR = {Robbins, Herbert and Monro, Sutton},
TITLE = {A stochastic approximation method},
JOURNAL = {Ann. Math. Statistics},
FJOURNAL = {Annals of Mathematical Statistics},
VOLUME = {22},
YEAR = {1951},
PAGES = {400--407},
ISSN = {0003-4851},
MRCLASS = {62.0X},
MRNUMBER = {42668},
MRREVIEWER = {R. P. Peterson},
DOI = {10.1214/aoms/1177729586},
URL = {https://doi.org/10.1214/aoms/1177729586},
}

@article{MiesSteland2022,
author = {Fabian Mies and Ansgar Steland},
title = {{Sequential Gaussian approximation for nonstationary time series in high dimensions}},
volume = {29},
journal = {Bernoulli},
number = {4},
publisher = {Bernoulli Society for Mathematical Statistics and Probability},
pages = {3114 -- 3140},
year = {2023},
doi = {10.3150/22-BEJ1577},
URL = {https://doi.org/10.3150/22-BEJ1577}
}

@article{MiesSteland2023,
author = {Fabian Mies and Ansgar Steland},
title = {{Projection inference for high-dimensional covariance matrices with structured shrinkage targets}},
volume = {18},
journal = {Electronic Journal of Statistics},
number = {1},
publisher = {Institute of Mathematical Statistics and Bernoulli Society},
pages = {1643 -- 1676},
year = {2024},
doi = {10.1214/24-EJS2225},
URL = {https://doi.org/10.1214/24-EJS2225}
}

@article{BartlettFosterTelgarsky2017,
title = {Spectrally-normalized margin bounds for neural networks},
author = {Bartlett, P.L. and Foster, D.J. and Telgarsky, M.J.},
journal={Advances in Neural Information Processing Systems}, 
fjournal={Advances in Neural Information Processing Systems}, 
volume = {30},
year={2017}
}

@article {AueHorvath2004,
AUTHOR = {Aue, Alexander and Horv\'{a}th, Lajos},
TITLE = {Delay time in sequential detection of change},
JOURNAL = {Statist. Probab. Lett.},
FJOURNAL = {Statistics \& Probability Letters},
VOLUME = {67},
YEAR = {2004},
NUMBER = {3},
PAGES = {221--231},
ISSN = {0167-7152},
MRCLASS = {62L10},
MRNUMBER = {2053524},
MRREVIEWER = {B. L. S. Prakasa Rao},
DOI = {10.1016/j.spl.2004.01.002},
URL = {https://doi.org/10.1016/j.spl.2004.01.002},
}

@article {Carlstein1986,
AUTHOR = {Carlstein, Edward},
TITLE = {The use of subseries values for estimating the variance of a
general statistic from a stationary sequence},
JOURNAL = {Ann. Statist.},
FJOURNAL = {The Annals of Statistics},
VOLUME = {14},
YEAR = {1986},
NUMBER = {3},
PAGES = {1171--1179},
ISSN = {0090-5364},
MRCLASS = {62G05 (60G10 62M09)},
MRNUMBER = {856813},
MRREVIEWER = {Sh. A. Khashimov},
DOI = {10.1214/aos/1176350057},
URL = {https://doi.org/10.1214/aos/1176350057},
}

@article {LohWainwright2015,
AUTHOR = {Loh, Po-Ling and Wainwright, Martin J.},
TITLE = {Regularized {$M$}-estimators with nonconvexity: statistical
and algorithmic theory for local optima},
JOURNAL = {J. Mach. Learn. Res.},
FJOURNAL = {Journal of Machine Learning Research (JMLR)},
VOLUME = {16},
YEAR = {2015},
PAGES = {559--616},
ISSN = {1532-4435},
MRCLASS = {62H12 (62F07 62J07 90C90)},
MRNUMBER = {3335800},
}

@article {JacodSoerensen2018,
AUTHOR = {Jacod, Jean and Soerensen, Michael},
TITLE = {A review of asymptotic theory of estimating functions},
JOURNAL = {Stat. Inference Stoch. Process.},
FJOURNAL = {Statistical Inference for Stochastic Processes. An
International Journal Devoted to Time Series Analysis and the
Statistics of Continuous Time Processes and Dynamical Systems},
VOLUME = {21},
YEAR = {2018},
NUMBER = {2},
PAGES = {415--434},
ISSN = {1387-0874},
MRCLASS = {62F12 (60J60 62M05 62M09 62M10)},
MRNUMBER = {3824976},
MRREVIEWER = {Rosa Maria Mininni},
DOI = {10.1007/s11203-018-9178-8},
URL = {https://doi.org/10.1007/s11203-018-9178-8},
}

@article{ GoukEtAl2021,
  author = {Gouk, H. and Frank, E. and Pfahringer, B. and Cree, M.J.},
  title = {Regularisation of neural networks by enforcing Lipschitz continuity},
  journal = {Mach. Learning},
  fjournal = {Machine Learning},
  volume = {110},
  pages = {393-416},
  year = {2021},
}

@InProceedings{pmlr-v48-yangc16,
title = 	 {Sparse Nonlinear Regression: Parameter Estimation under Nonconvexity},
author = 	 {Yang, Zhuoran and Wang, Zhaoran and Liu, Han and Eldar, Yonina and Zhang, Tong},
booktitle = 	 {Proceedings of The 33rd International Conference on Machine Learning},
pages = 	 {2472--2481},
year = 	 {2016},
editor = 	 {Balcan, Maria Florina and Weinberger, Kilian Q.},
volume = 	 {48},
series = 	 {Proceedings of Machine Learning Research},
address = 	 {New York, New York, USA},
month = 	 {20--22 Jun},
publisher =    {PMLR},
url = 	 {https://proceedings.mlr.press/v48/yangc16.html}
}

@article{Lederer2023,
AUTHOR = {Lederer, Johannes},
TITLE = {Statistical guarantees for sparse deep learning},
JOURNAL = {AStA Advances in Statistical Analysis},
FJOURNAL = {AStA Advances in Statistical Analysis},
YEAR = {2024},
VOLUME = {108},
PAGES = {231-258},
DOI = {10.1007/s10182-022-00467-3},
}

@misc{BanerjeeEtAl2022,
doi = {10.48550/ARXIV.2209.15106},

url = {https://arxiv.org/abs/2209.15106},

author = {Banerjee, Arindam and Cisneros-Velarde, Pedro and Zhu, Libin and Belkin, Mikhail},

title = {Restricted Strong Convexity of Deep Learning Models with Smooth Activations},

publisher = {arXiv},

year = {2022},

copyright = {arXiv.org perpetual, non-exclusive license}
}

@article {PeligradShao1995,
AUTHOR = {Peligrad, Magda and Shao, Qi Man},
TITLE = {Estimation of the variance of partial sums for {$\rho$}-mixing
random variables},
JOURNAL = {J. Multivariate Anal.},
FJOURNAL = {Journal of Multivariate Analysis},
VOLUME = {52},
YEAR = {1995},
NUMBER = {1},
PAGES = {140--157},
ISSN = {0047-259X},
MRCLASS = {60F05 (62F12)},
MRNUMBER = {1325375},
MRREVIEWER = {Jerome Senturia},
DOI = {10.1006/jmva.1995.1008},
URL = {https://doi.org/10.1006/jmva.1995.1008},
}

@book {ShorackWellner1986,
AUTHOR = {Shorack, Galen R. and Wellner, Jon A.},
TITLE = {Empirical processes with applications to statistics},
SERIES = {Wiley Series in Probability and Mathematical Statistics:
Probability and Mathematical Statistics},
PUBLISHER = {John Wiley \& Sons, Inc., New York},
YEAR = {1986},
PAGES = {xxxviii+938},
ISBN = {0-471-86725-X},
MRCLASS = {60-02 (60F17 62-02 62G10 62G30)},
MRNUMBER = {838963},
MRREVIEWER = {S\'{a}ndor Cs\"{o}rg\H{o}},
}

@book {Billingsley1999,
AUTHOR = {Billingsley, Patrick},
TITLE = {Convergence of probability measures},
SERIES = {Wiley Series in Probability and Statistics: Probability and
Statistics},
EDITION = {Second},
NOTE = {A Wiley-Interscience Publication},
PUBLISHER = {John Wiley \& Sons, Inc., New York},
YEAR = {1999},
PAGES = {x+277},
ISBN = {0-471-19745-9},
MRCLASS = {60B10 (28A33 60F17)},
MRNUMBER = {1700749},
DOI = {10.1002/9780470316962},
URL = {http://dx.doi.org/10.1002/9780470316962},
}

\end{document}